# Frequency-Volume Optimization of Arch-Dams using Charged System Search Method

Mohamad-Tagi Aalami[1], Siamak Talatahari[1,2*], Reza Parsiavash[1]

[1] *Department of Civil Engineering, University of Tabriz, Tabriz, Iran*

[2] *Faculty of Engineering & IT, University of Technology Sydney, Ultimo, NSW 2007, Australia*

**Abstract**

This paper presents a novel multi-objective formulation to investigate optimal shape of double curvature arch dams in dam-reservoir systems considering its volume and natural frequencies. For this purpose, multi-objective charge system search algorithm in combination with multi-criteria tournament decision-making method is utilized. Concrete volume and natural frequencies of dam body are considered as two main objective functions. Numerical modeling and frequency analysis of dam-reservoir system are implemented by developing a parallel-working MATLAB-APDL code. A case study of a well-known arch dam is performed to show the effectiveness and reliability of proposed methodology.

## 1. Introduction

The design of double-curvature arch dams is classified as a complex problem due to existing many efficient parameters on its geometry, the difficulty on numerical modeling (because of interaction between structure and surroundings), and the necessity of high computational

---

*Corresponding author: *E-mail address*: Siamak.talat@gmail.com (S. Talatahari)

costs. Finding an optimum design of such difficult problems is more complicated where we need to perform many different designs. The nonlinear nature of the required objective function and constraints of arch dam optimization problem and necessity of considering various boundary conditions make designers and researchers unavoidable to use novel optimization methods to handle this problem. Up to now, many of studies in the field of arch dam optimization are employed deterministic solution techniques [1-10]. Due to drawbacks of these methods such as requiring a good starting point, probability of reaching an unsatisfactory local minimum, requiring continuity and differentiability of objective function and having poor performance encountering with large-scale practical and complicated problems, there is an increasing interest in meta-heuristic approaches in recent years. Nevertheless, only a few studies are conducted using these approaches in the field of arch dam optimization problems.

Shape optimization of arch dams are performed by some meta-heuristic algorithms such as an improved particle swarm optimization [11], the genetic algorithm with the perturbation stochastic approximation [12]. In these studies dam volume is considered as the objective function and some geometrical, behavior and stability conditions are included as constraints. In spite of considering stress state as behavior constraints in some structural optimization studies, more recently, several attempts have been made to deal with frequency conditions [13]. In the other words, optimum design of arch dam is performed by limiting some natural dam frequencies to specific values [14]. In order to simplify the optimization process, only some frequencies and geometrical constraints were considered in these studies. Limitations are imposed to first-to-fourth natural frequencies in order to reduce the domain of vibration of the arch dam.

Aside from the fact that some simplifications are taken to account in previous studies, all of them were employed single-objective optimization approaches and therefore, the rest conditions were considered as constraints. The main achievement of such optimization approaches is a single result which may have a small objective function and satisfies the constraints. It should be noted that the final design necessarily is not a global optimum and some constraints may not be satisfied and therefore, the design may become an unfeasible design. In order to cover these defaults, this paper presents a new formulation for the frequency-volume optimum design of arch dams. We modified the problem in a way that two different types of objective functions are utilized, simultaneously. As a result, multi-objective formula for design of arch dames is developed in this paper for the first time. The first objective function is the concrete volume of dam body that should be minimized and the other ones are the different natural frequencies of the arch dam that should be maximized. To solve this problem, we develop a new multi-objective optimization methodology using a recent and efficient meta-heuristic optimization method, charged system search (CSS)[15].

## 2. Formulation of Frequency-Volume Multi-objective Optimization problem

Multi-objective optimization is a framework in which multiple objective functions are desired to have equal treatments. The solution to this problem is a set of multiple sub-solutions, which optimizes simultaneously the objectives. In addition, the feasible solutions have to comply with certain constraint conditions. This problem can be formulated as follows:

*Minimize:* $$\{fit_1(\boldsymbol{X}), fit_2(\boldsymbol{X}), ..., fit_k(\boldsymbol{X}), ..., fit_N(\boldsymbol{X})\} \quad (1)$$

*Subject to:* $$\Phi_i(\mathbf{X}) > 0 \qquad i = 1,2,3,...,p \quad (2)$$

where $fit_k(\boldsymbol{X})$ and $\Phi_i(\mathbf{X})$ denote the different objectives and the constraint condition, respectively; $N$ and $p$ are the number of objectives and constraints, respectively and $\mathbf{X}$ is the feasible set of decision variables. In this paper, two types of objectives are considered:

**Objective 1**: The concrete volume of the arch dam that should be minimized. It could be determined by integration on the dam surfaces, as:

*Minimize:* $$fit_1(\boldsymbol{X}) = \iint_A \left| y_u(x,z) - y_d(x,z) \right| dA \quad (3)$$

in which, $y_u(x,z)$ and $y_d(x,z)$ are the upstream and downstream surfaces of the arch dam, respectively and $A$ is an area produced by projecting the dam body on a $xz$ plan.

**Objective 2**: The natural frequencies of the arch dam that should be maximized:

*Maximize:* $$fit_2(\boldsymbol{X}) = fr_n \qquad n = 1,2,\ldots,n_{fr} \quad (4)$$

where, $n_{fr}$ is number of natural frequencies that may be considered.

Here, design variables vector, $\boldsymbol{X}$, can be adopted as:

$$\boldsymbol{X} = \{s \quad \beta \quad tc_1 \ldots tc_{n+1} \quad ru_1 \ldots ru_{n+1} \quad rd_1 \ldots rd_{n+1}\} \tag{5}$$

$\boldsymbol{X}$ contains $2+3(n+1)$ shape parameters of arch dam where, $n$ is number of divisions along the dam height.

In this optimization problem, various constraints including geometrical and stability constraints are taken into account. To ensure that upstream and downstream faces of the dam do not pass through each other Eq. (6) should be satisfied. Also for constructing facilities and having smooth cantilevers over the height of the dam, the slope of overhang at the upstream and downstream faces of the dam should satisfy Eq. (7), as:

$$r_{di} \le r_{ui} \Rightarrow \frac{r_{di}}{r_{ui}} - 1 \le 0, i = 1,2,..6 \tag{6}$$

$$s \le s_{alw} \Rightarrow \frac{s}{s_{alw}} - 1 \le 0 \tag{7}$$

where $r_{di}$ and $r_{ui}$ are the radius of curvature at $i$th level in the $z$ direction and $s = \cot\alpha$ is the slope of overhang at the downstream and upstream faces of the dam. $s_{alw}$ is the allowable absolute value for the aforementioned parameter. To ensure the sliding stability of the dam, following equation should be satisfied:

$$\varphi_l \le \varphi \le \varphi_u \tag{8}$$

where $\varphi$ is the central angle of arch dam at $i$th level in the $z$ direction that varies from $90^\circ$ to $130^\circ$ throughout the dam height.

Some geometric conditions for Double-curvature arch dams are also required. The following equation shows one polynomial of second order that is used to determine the curve of upstream boundary as shown in Fig. 1(a), [16, 17]:

$$g(z) = \frac{\gamma z^2}{2\beta h} - \gamma z \tag{9}$$

where $h$ is dam height and the point where the slope of the upstream face equals to zero is $z = \beta h$. By dividing the dam height into $n$ segments, a polynomial function can be employed to determine the thickness of crown cantilever:

$$t_c(z) = \sum_{i=1}^{n+1} L_i(z)\, t_{ci} \tag{10}$$

where $L_i(z)$ is a Lagrange interpolation formula and $t_{ci}$ is the thickness of central vertical section at $i$th level. $L_i(z)$ can be expressed as:

$$L_i(z) = \prod_{\substack{m=1 \\ m \neq i}}^{n+1} \frac{z - z_m}{z_i - z_m} \tag{11}$$

in which, $z_i$ denotes the $z$ coordinate of $i$th level. As shown in Fig. 1(b), for the purpose of symmetrical canyon and arch thickening from crown to abutment, the shape of the horizontal section of a parabolic arch dam is determined by the following two parabolas [16, 17]:

$$y_u(x, z) = \frac{x^2}{2 r_u(z)} + g(z) \tag{12}$$

$$y_d(x, z) = \frac{x^2}{2 r_d(z)} + g(z) + t_c(z) \tag{13}$$

where $y_u(x, z)$ and $y_d(x, z)$ are parabolas of upstream and downstream faces, respectively. $r_u$ and $r_d$ are radius of curvature for upstream and downstream curves at z direction and can be interpolated by $L_i(z)$ according to following equations of the $n$th order:

$$r_u(z) = \sum_{i=1}^{n+1} L_i(z)\, r_{ui} \tag{14}$$

$$r_d(z) = \sum_{i=1}^{n+1} L_i(z)\, r_{di} \tag{15}$$

in whicfh $r_{ui}$ and $r_{di}$ are the related values of $r_u$ and $r_d$ at the controlling levels, respectively.

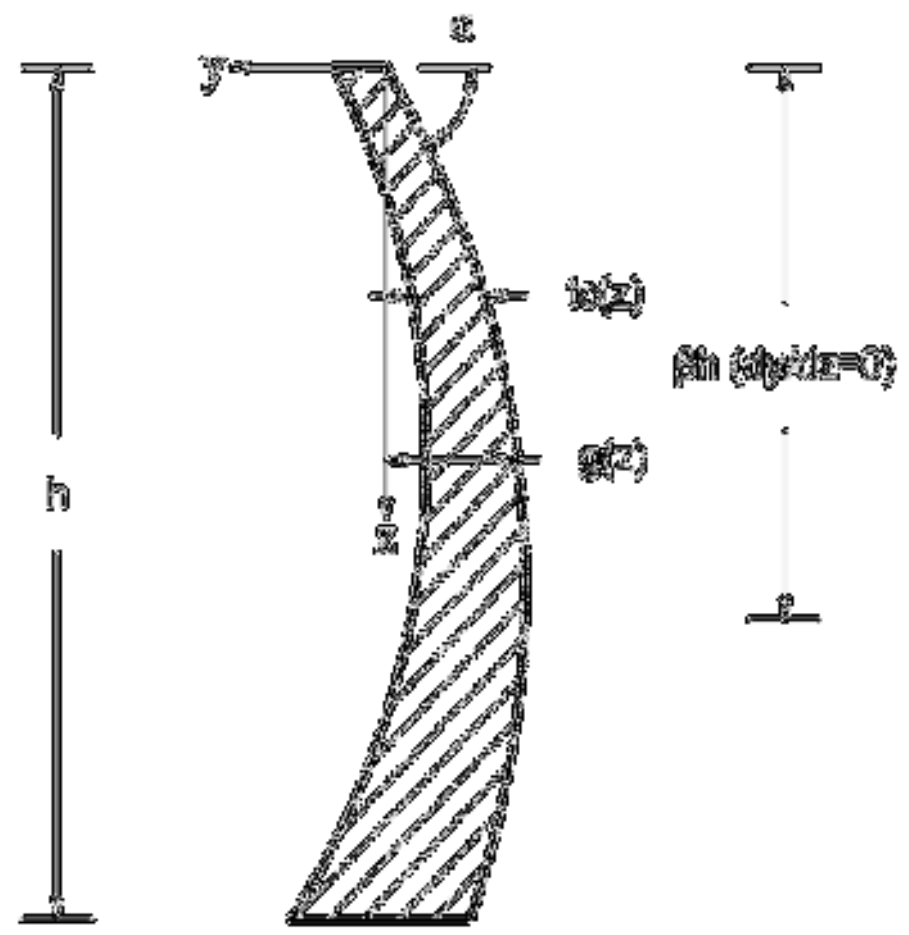

(b)

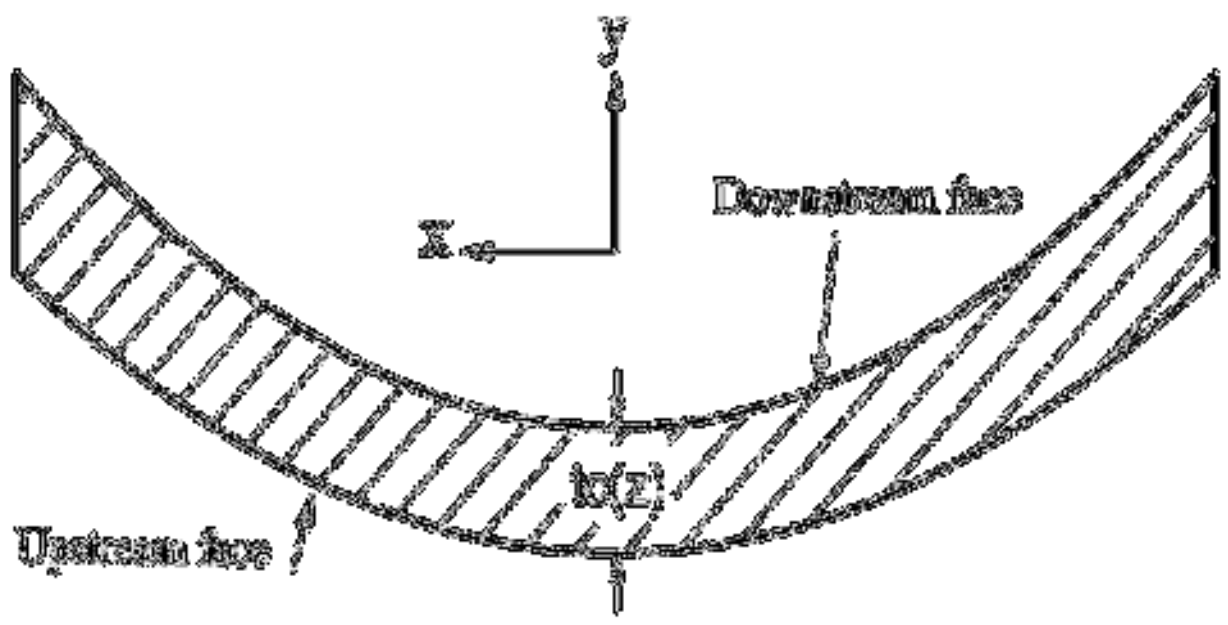

Fig. 1. Schematic of arch dam: a) Crown Cantilever profile b) parabolic shape of an elevation

## 3. Multi-objective Charged System Search

The Charged System Search (CSS) is inspired by the Coulomb law from electrostatics and the laws of motion from Newtonian mechanics [15]. In the CSS each agent (CP) is considered as a charged sphere with radius a, having a uniform volume charge density which can produce an electric force on the other CPs. The force magnitude for a CP located in the inside of the sphere is proportional to the separation distance between the CPs, while for a CP located outside the sphere, it is inversely proportional to the square of the separation distance between the CPs. Considering this resultant forces, each CP moves toward its new position and is this way the new solution is generated. Application of the CSS method for solving many benchmark and engineering problems shows that, it outperforms evolutionary algorithms and comparison of the results demonstrates the efficiency of the present algorithm [15].

Some modifications and additional steps are considered in the single-objective charged system search (CSS) algorithm for presenting a multi-objective CSS method (MoCSS). The pseudo-code for the MoCSS algorithm can be summarized as follows:

**Step 1:** *Initialization*. Similar to the standard CSS, the initial positions of CPs are determined randomly in the search space and the initial velocities of CPs are assumed to be zero. The values of the fitness functions for the CPs are determined and the CPs are ranked based on Pareto dominance criteria:

$$\forall i \neq j \qquad \exists fit_i \;|fit_i(\mathbf{X}A) \leq fit_i(\mathbf{X}B) \wedge fit_j(\mathbf{X}A) \leq fit_j(\mathbf{X}B) \tag{16}$$

**Step 2:** *definition of CM*. all CPs belong to first rank (according to the Pareto dominance criteria) are stored in a memory, so called charged memory (CM).

**Step 3:** *Determination of forces on CPs*. The force vector is calculated for each CP as:

$$\mathbf{F}_j = \sum_{i,i\neq j} \left( \frac{q_i}{a^3} r_{ij} \cdot i_1 + \frac{q_i}{r_{ij}^2} \cdot i_2 \right) ar_{ij} \, p_{ij} (\mathbf{X}_i - \mathbf{X}_j) \left\langle \begin{array}{l} j = 1,2,\ldots,N \\ i_1 = 1, i_2 = 0 \Leftrightarrow r_{ij} < a \\ i_1 = 0, i_2 = 1 \Leftrightarrow r_{ij} \geq a \end{array} \right. \tag{17}$$

where $\mathbf{F}_j$ is the resultant force acting on the $j$th CP; $N$ is the number of CPs. The magnitude of charge for each CP ($q_i$) is defined considering the quality of its solution as:

$$q_i = \prod_{k=1}^{N} \frac{fit_k(i) - fitworst_k}{fitbest_k - fitworst_k}, \qquad i = 1,2,\ldots,N \tag{18}$$

where operator $\Pi$ is used for repeated multiplications. $fitbest_k$ and $fitworst_k$ are the best and the worst fitness of the $k$th objective function for all CPs, respectively; $fit_k(i)$ represents the fitness of the agent $i$. $N$ is the number of objective functions. $P_{ij}$ is the probability of moving each CP towards the others and is obtained using the following function:

$$p_{ij} = \begin{cases} 1 & rank(i) > rank(j) \\ 0 & otherwise \end{cases} \tag{19}$$

In Eq. (17), $ar_{ij}$ indicates the kind of force and is defined as:

$$ar_{ij} = \begin{cases} +1 & rand < 0.8 \\ -1 & \text{otherwise} \end{cases} \tag{20}$$

where *rand* represents a random number.

**Step 4:** *Solution construction.* Each CP moves to the new position and the new velocity is calculated as:

$$\mathbf{X}_{j,new} = rand_{j1} \cdot k_a \cdot \mathbf{F}_j + rand_{j2} \cdot k_v \cdot \mathbf{V}_{j,old} + \mathbf{X}_{j,old} \tag{21}$$

$$\mathbf{V}_{j,new} = \mathbf{X}_{j,new} - \mathbf{X}_{j,old} \tag{22}$$

where $k_a$ is the acceleration coefficient; $k_v$ is the velocity coefficient to control the influence of the previous velocity; and $rand_{j1}$ and $rand_{j2}$ are two random numbers uniformly distributed in the range (0,1).

**Step 5:** *Position updating process.* If a new CP exits from the allowable search space, a harmony search-based handling approach is used to correct its position. According to this method, any variable of each solution that violates its corresponding boundary can be regenerated from CM as shown in Fig. 2. In this figure, "$w.p.$" is the abbreviation of "with the probability"; CMCR (the Charge Memory Considering Rate) varying between 0 and 1 sets the rate of selecting a value in the new vector from historic values stored in CM, and (1–CMCR) sets the rate of randomly choosing one value from a possible range of values. The value (1–PAR) sets the rate of doing nothing, and PAR sets the rate of choosing a value from neighboring the best CP.

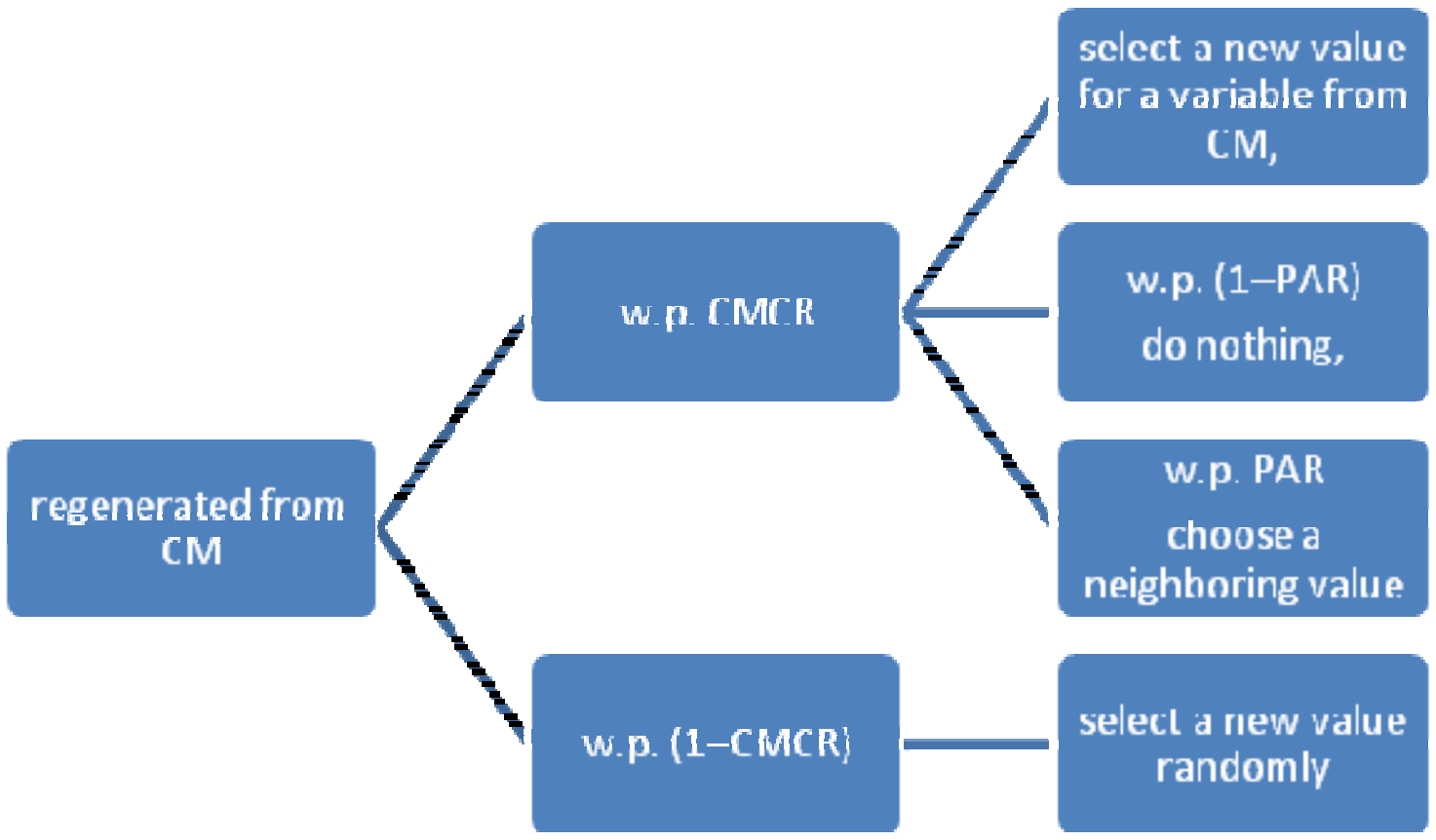

Fig. 2. Regenerating violated solutions from CM

**Step 6:** *CM updating.* All new CPs and available ones in the CM are ranked using Eq. (16) and then the first Pareto are saved in the CM; then the agents with small objective distances will be deleted from the memory. This distance is defined as [18, 19]:

$$d_{ij} = \sqrt{\sum_{k=1}^{N} \left( u_k \left( fit_k (\boldsymbol{X}_i) - fit_k (\boldsymbol{X}_j) \right) \right)^2} \qquad (23)$$

in which, $d_{ij}$ is the objective distance between solution $\boldsymbol{X}_i$ and $\boldsymbol{X}_j$ and $u_k$ is a constant chosen to make all $u_k fit_k$ close to each other. To determine the appropriate value of $u_k$, first we set all $u_k = 1$. Then by optimizing the problem, some near optimal solutions and $fitworst_k$ for each objective are obtained. Eventually, $u_k$ is obtained in the following way:

$$u_1 = \alpha > 0,\ u_k = u_{k-1} * \frac{fitworst_k}{fitworst_{k-1}} \qquad for\ fit_k, \qquad k > 1 \qquad (24)$$

**Step 7:** *Termination criterion control.* Steps 3-6 are repeated until a termination criterion is satisfied.

**Step 8:** *Multi Criteria Decision Making.* Since the Pareto optimal solutions cannot be ranked globally, the decision maker needs to define some extra information to reach the desired answers. Here, we use the multi-criteria tournament decision making method (MTDM) [20] because it is a simple method which uses the preferences of decision maker in the form of criterion importance weight to rank solutions from the best to the worst. The Preferences of DM can be expressed in terms of a global ranking function R. To fulfill this aim, tournament function $T_i\left(a,A\right)$, should be defined as following:

$$T_i\left(a,A\right)=\sum_{\forall b\in A,a\neq b}\frac{t_i\left(a,b\right)}{\left(\left|A\right|-1\right)}, \tag{25}$$

where:

$$t_i\left(a,b\right)=\begin{cases}1 & \textit{if } fit_i\left(b\right)-fit_i\left(a\right)>0\\ 0 & \textit{otherwise}\end{cases} \tag{26}$$

This function counts the ratio of times alternative $a$ wins the tournament against each other $b$ solution from $A$ in which, $a$ is a non-dominated point in the objective space.

The function $R$ is then defined in order to gather all individual criterion and their appropriate priority weights, $w_i$, into the global ranking function:

$$R\left(a\right)=\left(\prod_{i=1}^{N}T_i\left(a,A\right)^{w_i}\right)^{\frac{1}{N}} \tag{27}$$

in which $N$ is the number of objective functions and $w_i$ should be specified by the DM under following conditions:

$$w_i>0 \quad \textit{and} \quad \sum_{i=1}^{n}w_i=1. \tag{28}$$

The flowchart of proposed method is presented in Fig. 3.

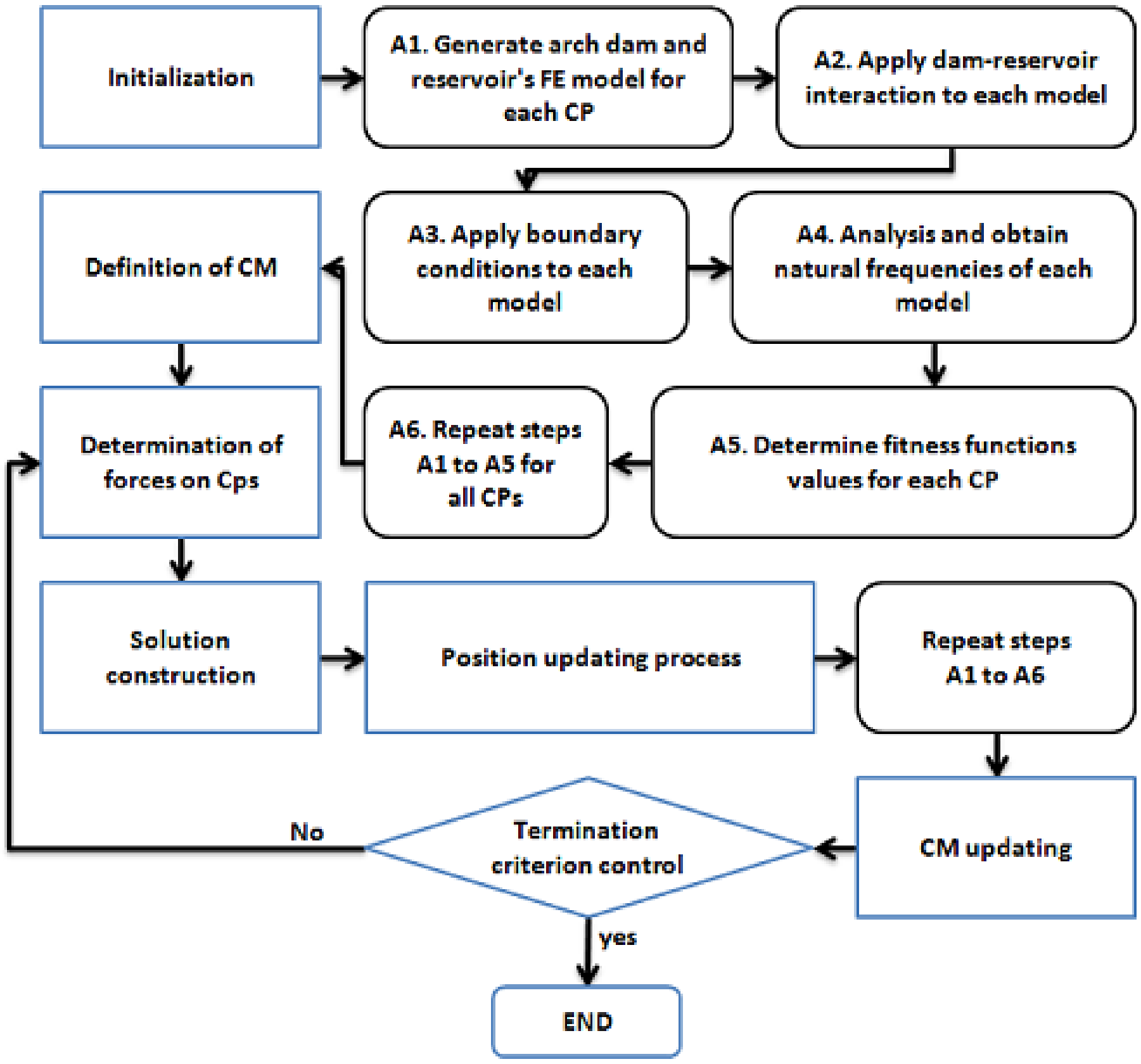

Fig. 3. Summarized flowchart of proposed method

## 4. Numerical Results

In this paper, in order to evaluate the proposed methodology, the Morrow Point arch dam as a real-world benchmark problem has been investigated. For this purpose, firstly, a numerical model of dam-reservoir system is developed and then to assure the reliability of proposed model, the presented model is verified with experimental and analytical results from literature. Then, the problem is optimized using the developed method.

For adapting the dam optimization problem with the developed method, the three types of problems are defined. For all these problems, the first objective function is the concrete volume of the arch dam as described in Eq. (3). However the second objective (natural frequencies) is different. At first, we exchange the maximization problem to a minimization ones as:

*Minimize:* $$fit_2(\boldsymbol{X}) = \frac{1}{fr_n} \qquad n = 1, 2, \ldots, n_{fr} \tag{28}$$

For the first problem (P1), Eq. (28) is used as the other objective functions; clearly, the number of these functions is related to the required frequencies to be maximized. The second problem (P2) just uses one function as the second objective:

*Minimize:* $$fit_2(\boldsymbol{X}) = \sum_{n=1}^{n_{fr}} \frac{1}{fr_n} \qquad n = 1, 2, \ldots, n_{fr} \tag{29}$$

Here, the summation of frequencies is utilized as the second objective. In the third problem (P3), the multiplication of inverse of frequencies are considered, as:

*Minimize:* $$fit_2(\boldsymbol{X}) = \prod_{n=1}^{n_{fr}} \frac{1}{fr_n} \qquad n = 1, 2, \ldots, n_{fr} \tag{30}$$

In this way for P1, we have many objectives while for P2 and P3, there are only two objective to be optimized. In this paper, we use first 10 frequents of structure.

Modeling and analyzing of the arch dam are handled using a combination of parallel working MATLAB and Ansys Parametric Design Language (MATLAB-APDL) codes that interfaces with the MoCSS to find the Pareto solutions. Reservoir is supposed to be full and the dam-reservoir interaction is taken to account in this example. The lower and upper bounds of design variables are obtained using classic design approaches [21] and presented in Table 2.

| Table 1: The lower and upper bounds of design variables |
| --- |

### *4.1. Arch Dam modeling*

To develop the finite element model, 8-node solid elements are employed for modeling the dam body. This element has three degrees of freedom at each node including translations in the nodal *x*, *y* and *z* directions. For modeling the reservoir, the 8-node fluid element with four degrees of freedom per node including translations in the *x*, *y* and *z* direction and pressure is used. The translation, however, are applicable only at nodes that are in the interface. The number of nodes and elements are not constant during the optimization process due to change in dimensions of the arch dam and mesh generation in each analysis, so it will be varied when needed. To prevent extra complexities, only the dam and reservoir interaction are considered and foundation rock is assumed rigid. The dam is treated as a 3D linear structure. Fluid medium is assumed as homogeneous, isotropic, irrotational and inviscid by linear compressibility. Fluid-structure interaction (FSI) effects are taken to account between reservoir and dam body as well as reservoir walls.

In order to verify analytical model, a finite element model (FEM) of morrow point arch dam-reservoir system has been developed (Fig. 4). This double curvature thin-arch concrete structure is located in the Gunnison River in west-central Colorado about 35 km east of Montrose. The dam has 142.65 m high and 220.68 m long along the crest also Its thickness varies from 3.66 m at the crest to 15.85 m at the base level [22]. Material properties of both dam body and water are presented in Table 3. The natural frequencies of arch dam with both empty and full reservoir are obtained using finite element analysis (FEA) and compared to experimental and analytical results from other works as shown in Table 4. It can be observed that good correspondence has been achieved between the results of present work with those reported in the literatures.

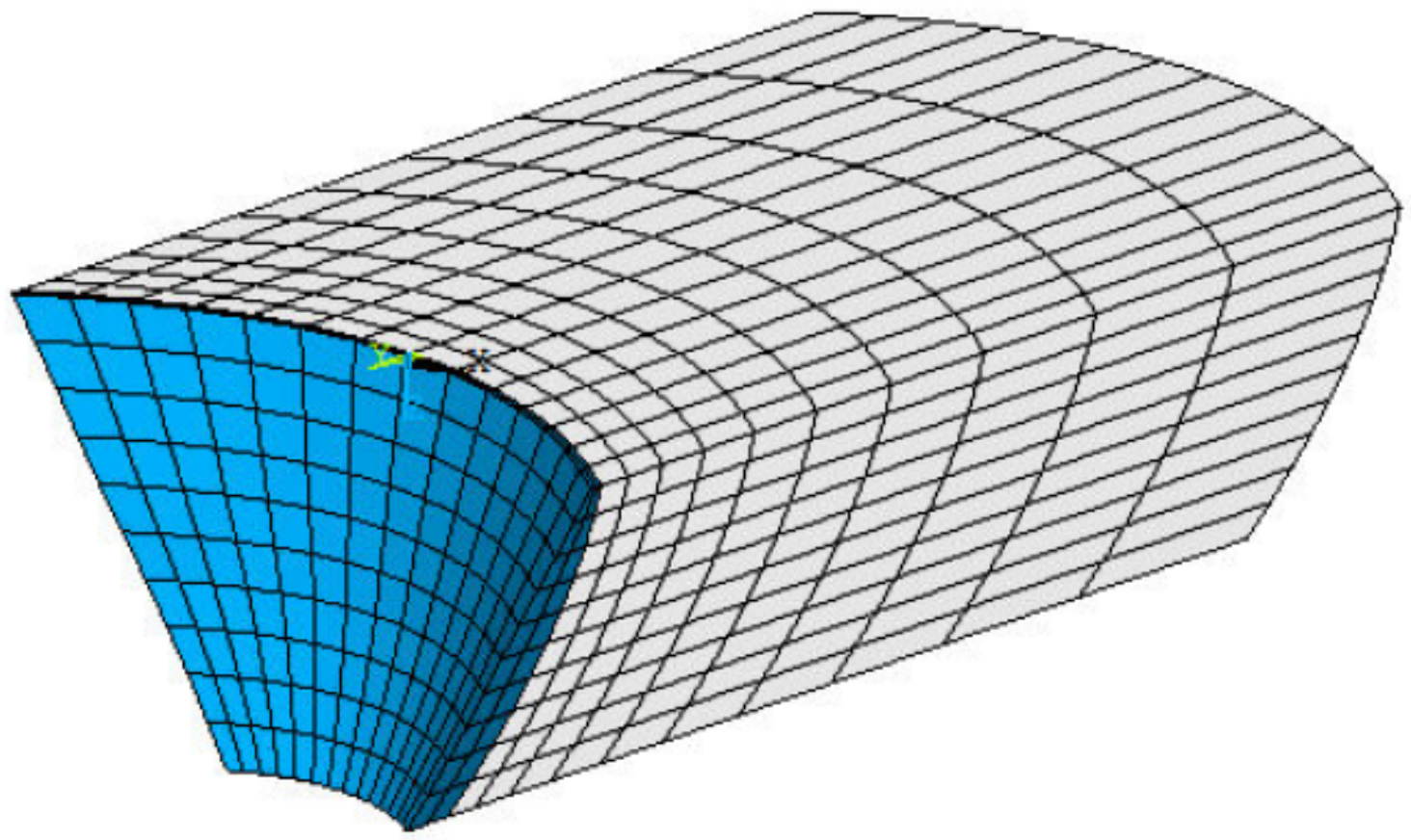

| |
|---|
| Fig. 4. Finite Element Model of Morrow Point Dam |
| Table 2: Material Property of Mass Concrete and Water |
| Table 3: Comparison the natural frequencies from the literature with the FEM |

### *4.2. Results and discussion*

In addition to the proposed MoCSS, the multi objective particle swarm optimization (MoPSO) [23] and a fast elitist non-dominated sorting genetic algorithm (NSGA-II) [18] are utilized as the well-known multi-objective optimization methods. For each problem, twenty runs are implemented with different random initial populations and best solutions of each algorithm are reported as the final design. The parameters of these methods are presented in Table 4 in which $N_{\text{var}}$ is the number of decision variables. In order to have equal computational effort for all algorithms, the same number of search iterations and population size are used for these algorithms. The Pareto fronts of the P1, P2 and P3 problems are shown in Figs 5, 6 and 7, respectively. The horizontal axes in these diagrams, are concrete volume of dam body (*fit1*) and the vertical axes are a combination of reversed frequencies, depended on the type of the problem. For the problem P1, there are eleven objective functions; therefore, we present eleven Pareto diagrams for this problem. In order to compare of these diagrams, the extreme values of Pareto solutions are shown in Tables 6, 6 and 8 for these problems, respectively. In addition, Figs 8, 9 and 10 present the related volume-frequency diagrams obtained by the three methods for problems P1, P2 and P3, respectively. Regarding the equality of employed criteria, comparing these diagrams is a suitable way to evaluate the efficiency of the used algorithms.

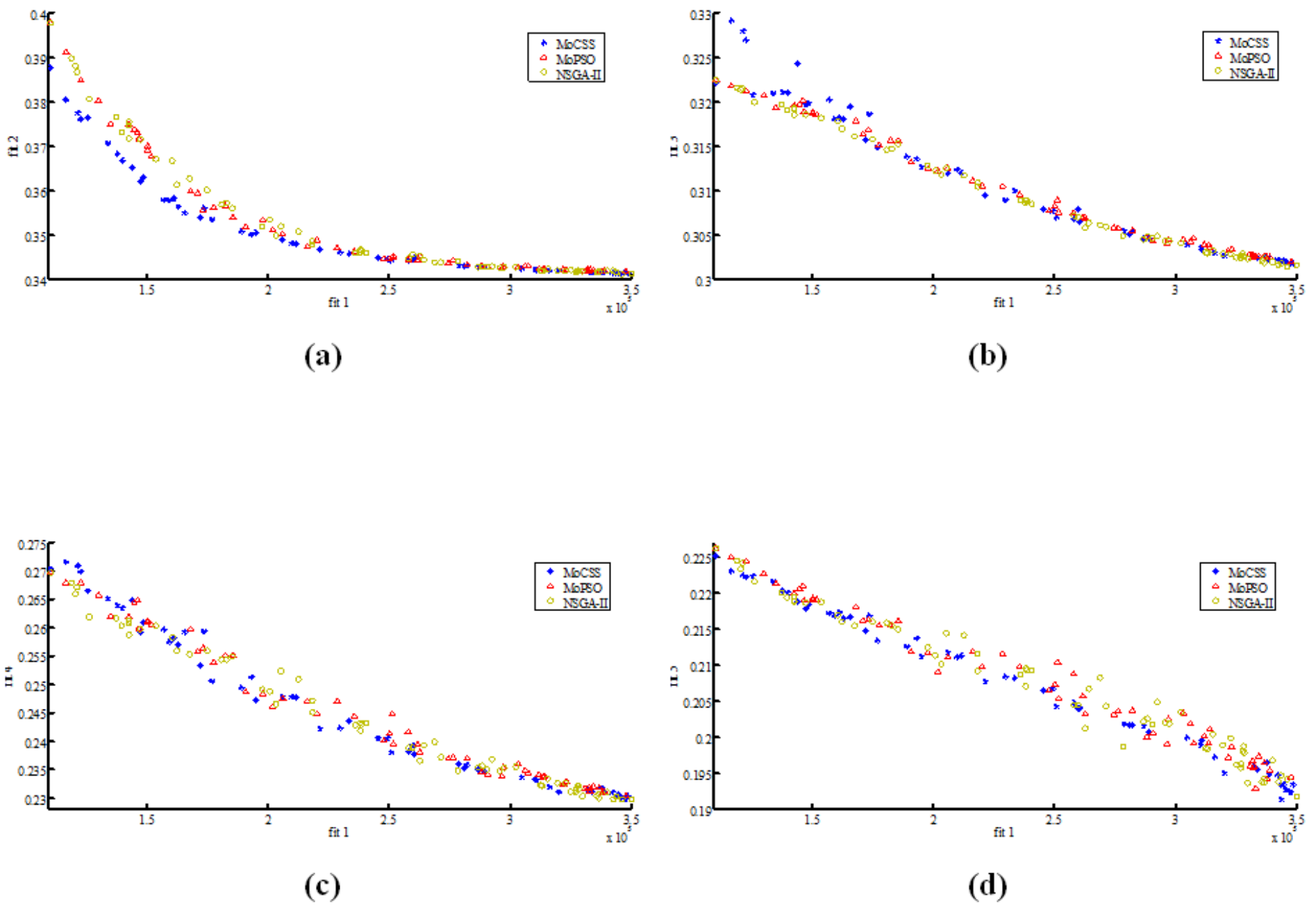

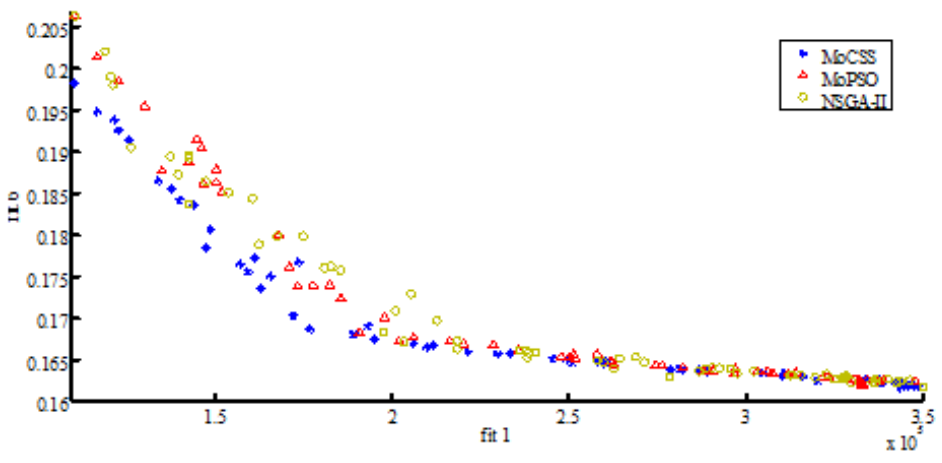

**(e)**

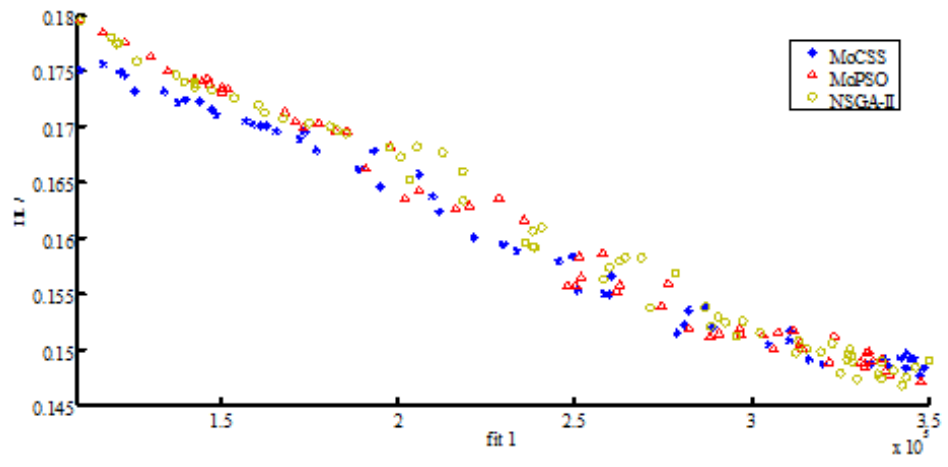

**(f)**

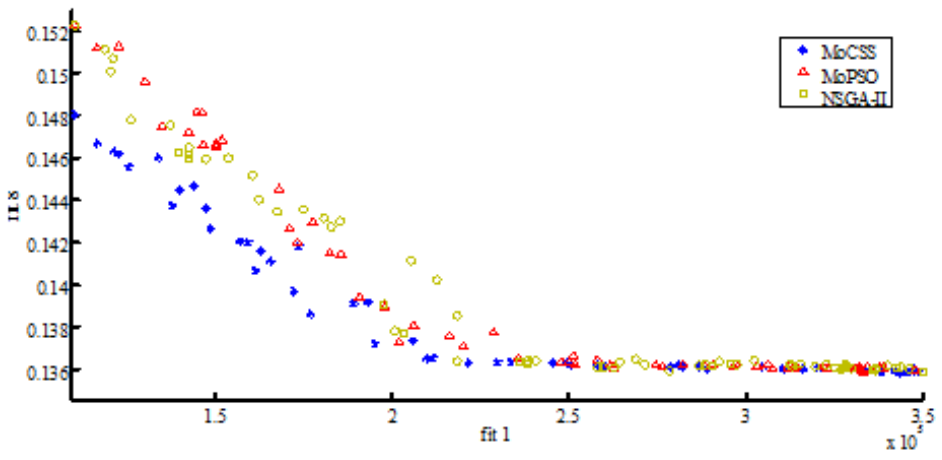

**(g)**

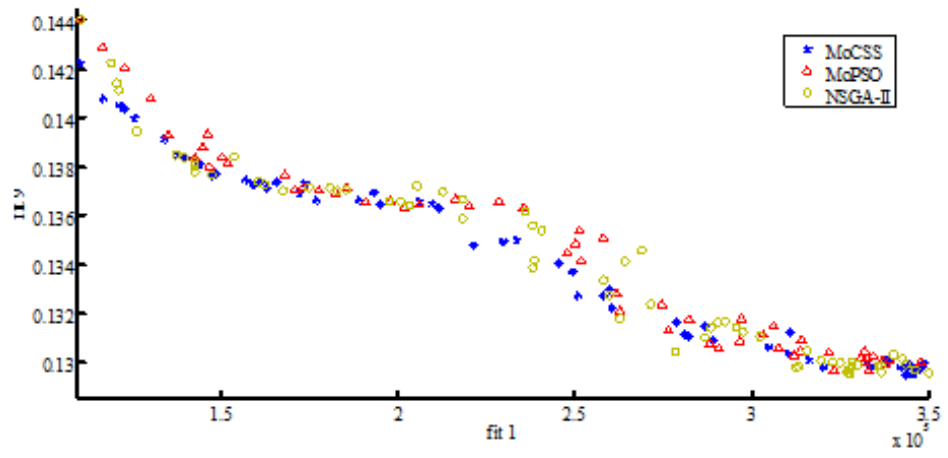

**(h)**

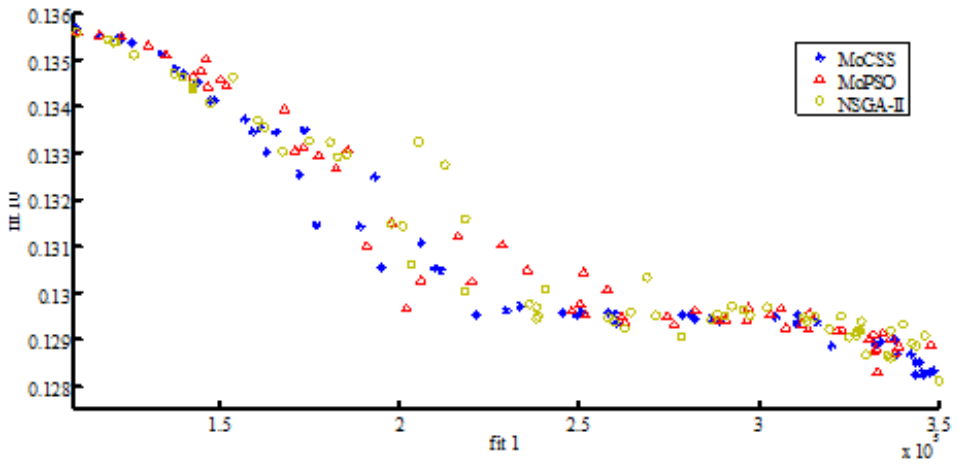

**(i)**

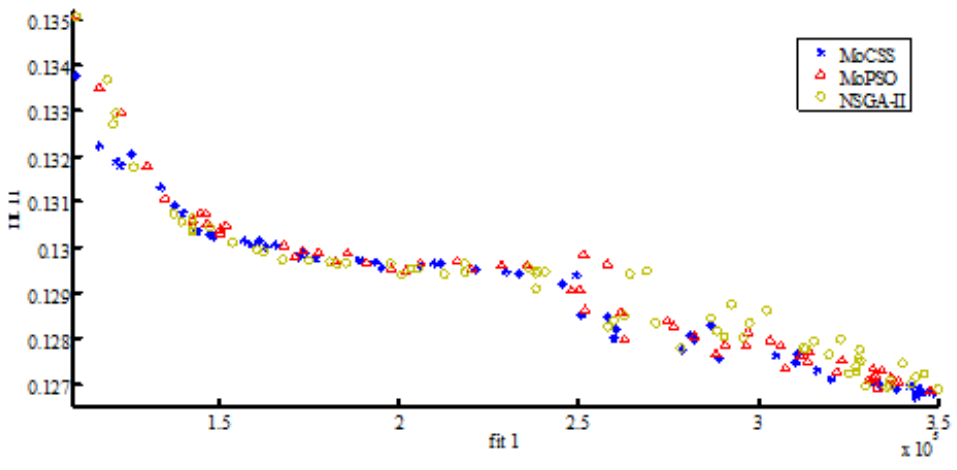

**(j)**

Fig. 5. Pareto fronts of three methods for P1

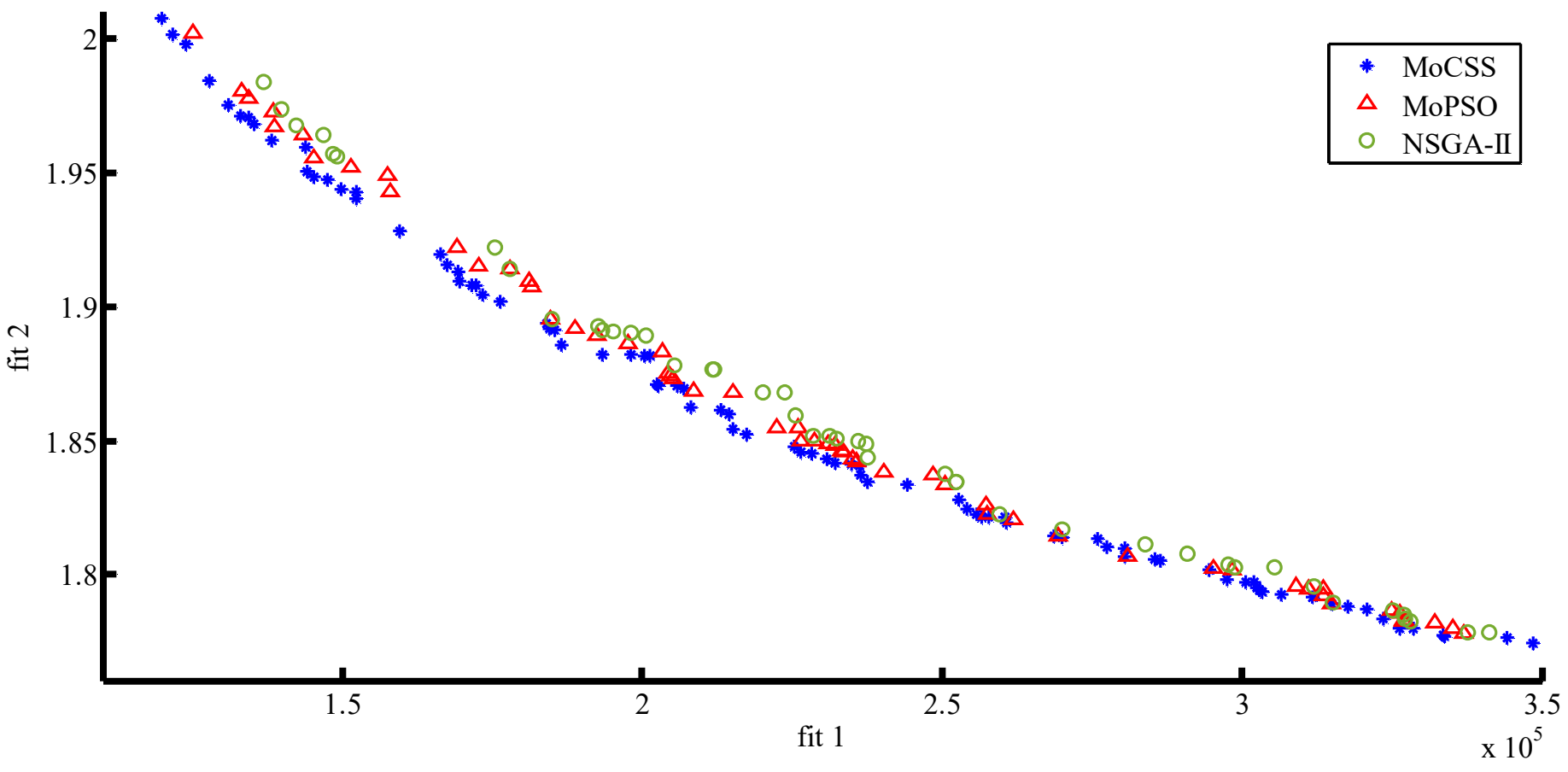

Fig. 6. Pareto front of three methods for P2

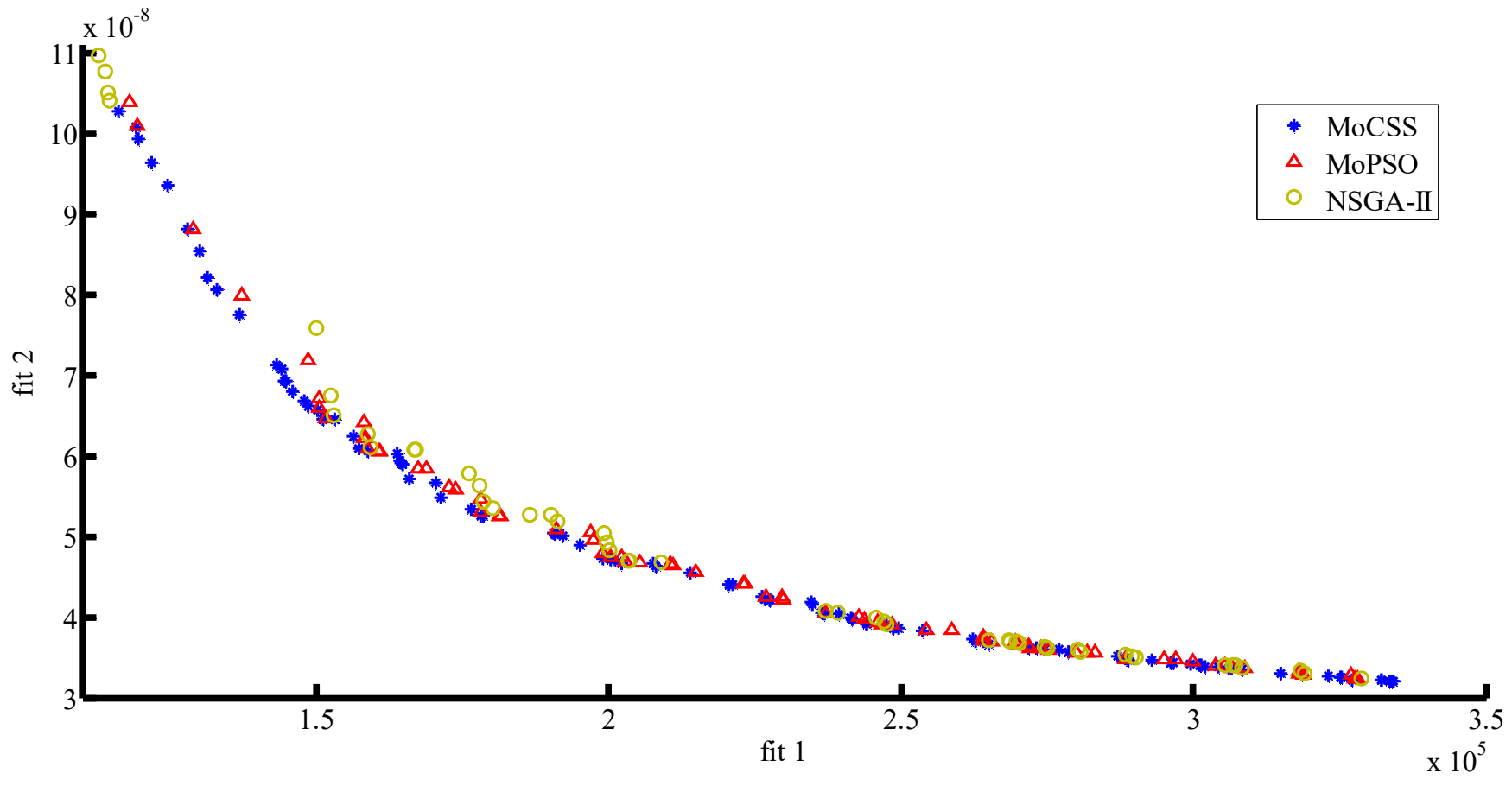

Fig. 7. Pareto front of three methods for P3

Table 4: Parameters of the MoCSS method

Table 5: Extreme values for the problem P1

Table 6: Comparison of the results for problem P2

Table 7: Comparison of the results for problem P3

As it can be seen clearly from Table 6, the MoCSS yields better extreme results comparing two other approaches. To illustrate, considering the extreme value of volume (*fit1*), however all algorithms have reached to an identical volume, MoCSS yields higher frequencies, which is more desirable. Besides, for extreme values of third to eleventh objective functions (second to tenth frequencies), the MoCSS finds better results and just in the case of 2nd, 6th, 8th and 9th objective functions, the other methods performs better. As it can be observed from Fig. 8, which compares the volume-frequency values of three methods for P1, the results of MoCSS are higher than those of other two approaches and have good dispersion among search space. It is worth to mention that the main objective is to reach maximum frequency in a definite volume. It is clear that, the MoCSS yields designs with higher frequencies in a certain volume. In the problem P2, the best extreme values are obtained by MoCSS, which are minimum volume for *fit1* and designs with maximum of third to seventh and tenth frequencies for *fit2*. Extreme value of MoPSO has maximum in eighth frequency and NSGA-II has maximum in ninth frequency (Table 7). Fig. 9 clearly reveals that the dispersion and quality of the MoCSS results are better than results of two other approaches. Although differences between the results are slight in some cases and much in others, by comparison ten volume-frequency diagrams (Fig. 9), it is apparent that the MoCSS yields better results than other algorithms. In the case of the problem P3, the MoCSS also performs better than two other methods as shown in Fig. 10. In comparison of volume-frequency diagrams, the reason of oscillation of results returns to nature of defined objective functions. In other words, when we investigate the Pareto fronts of these two problems (Figs. 6-7) the efficiency of the MoCSS becomes clear. Table 8 depicts the extreme values for problem P3.

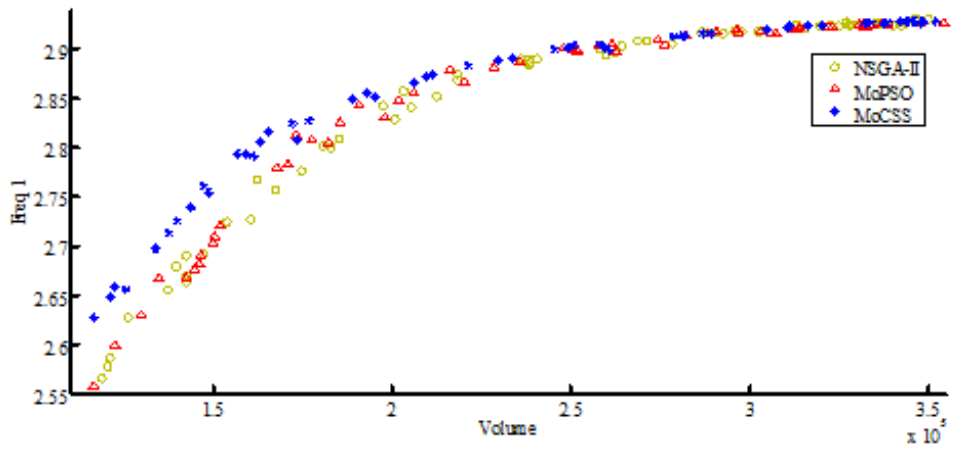

**(a)**

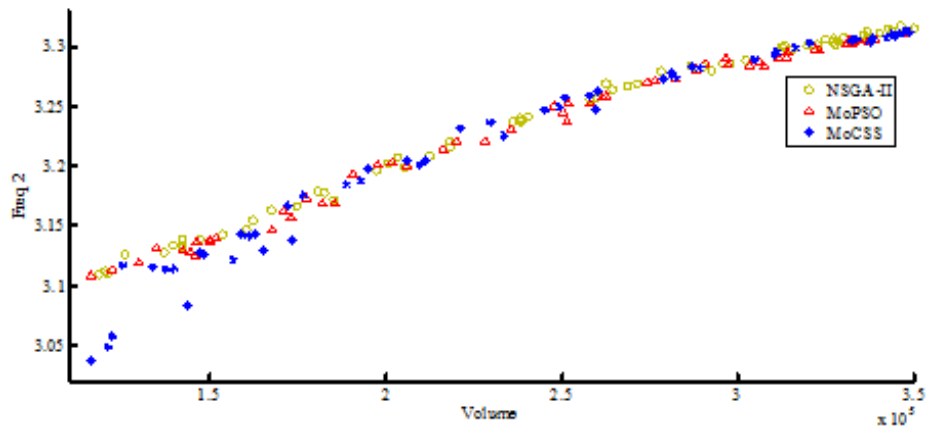

**(b)**

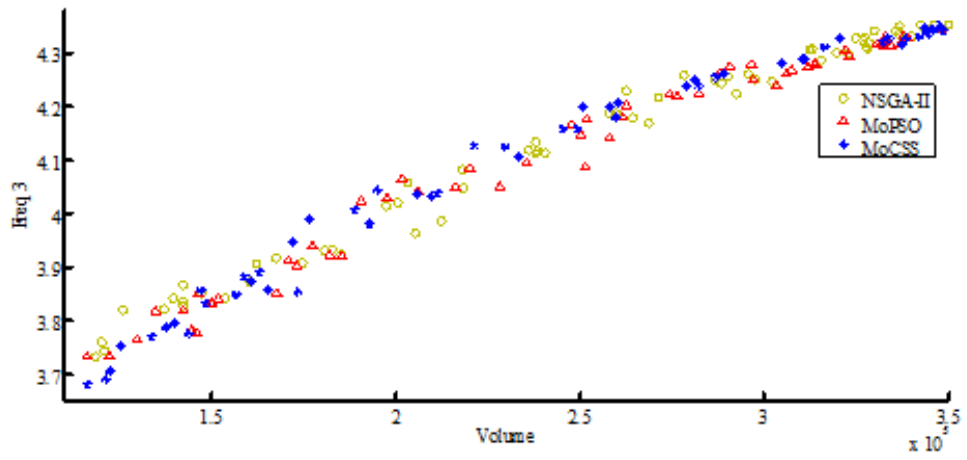

**(c)**

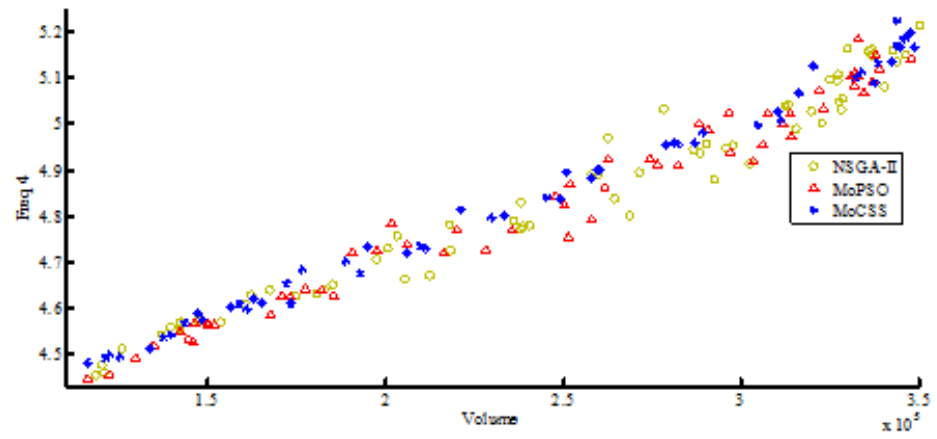

**(d)**

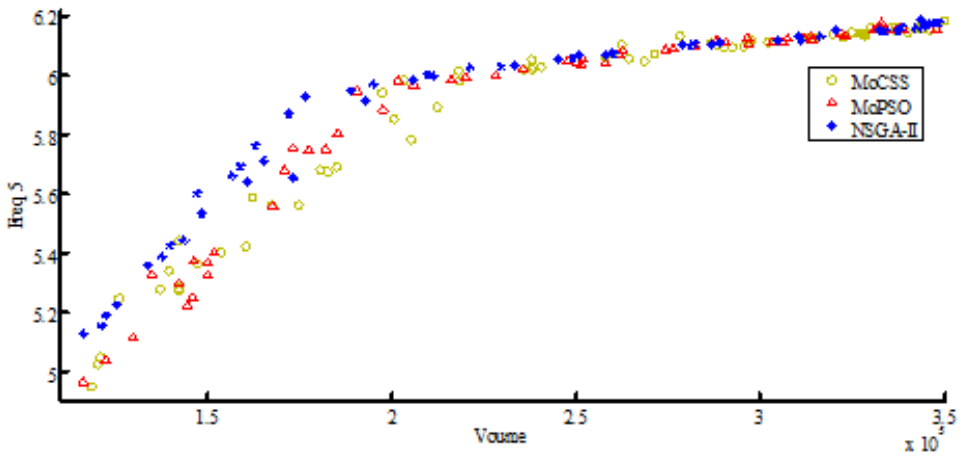

**(e)**

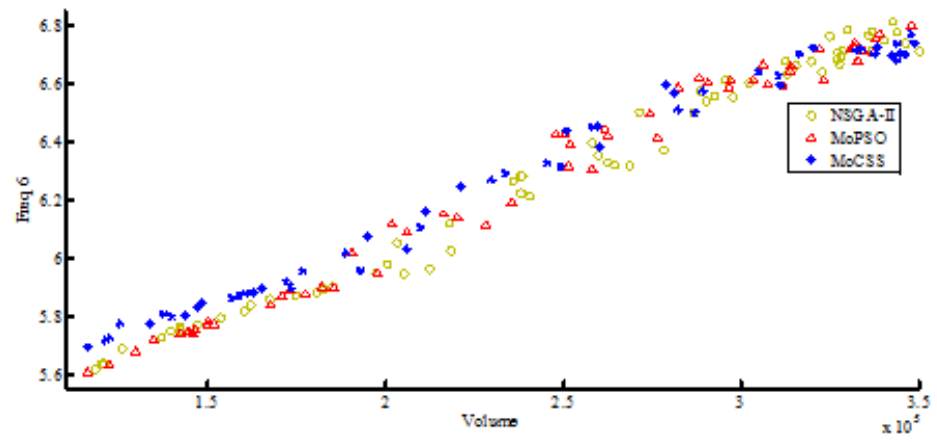

**(f)**

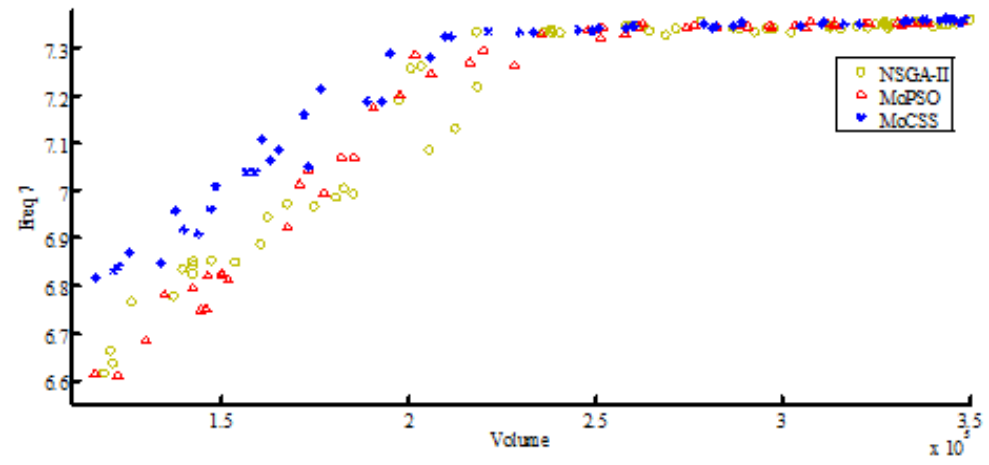

**(g)**

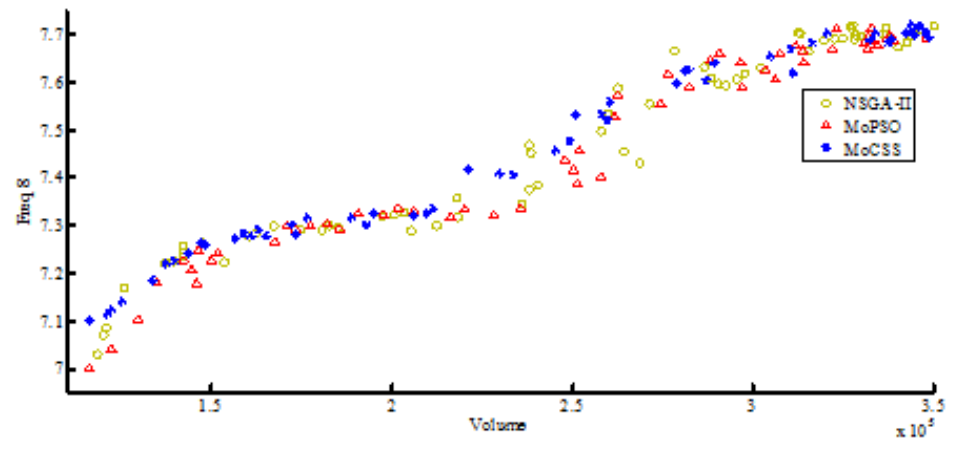

**(h)**

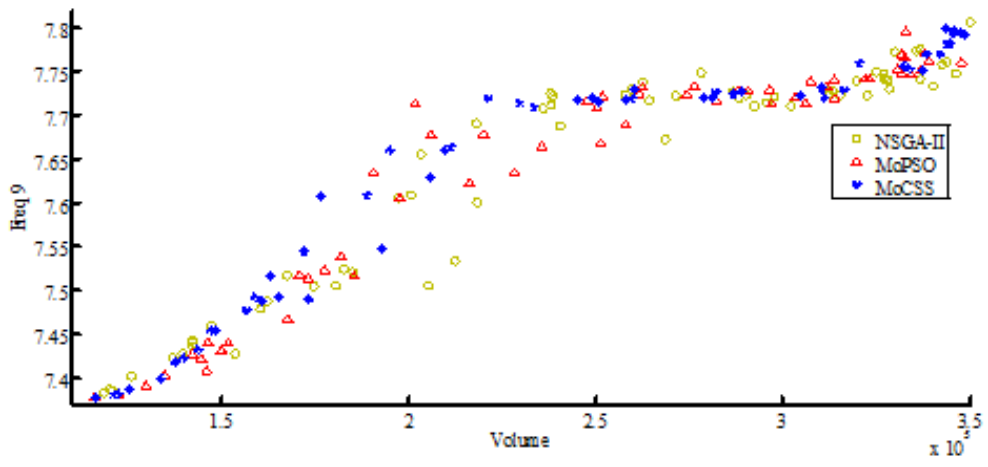

**(i)**

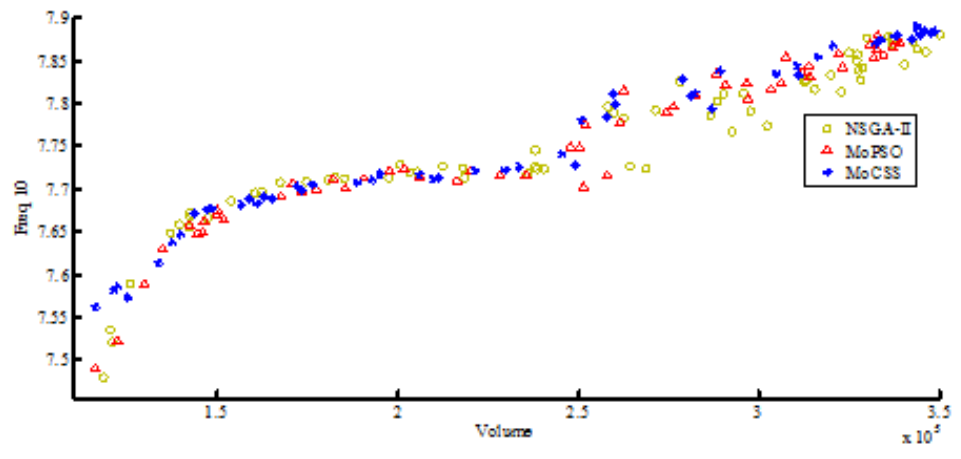

**(j)**

Fig. 8. The related Volume-Frequencies of three methods for P1

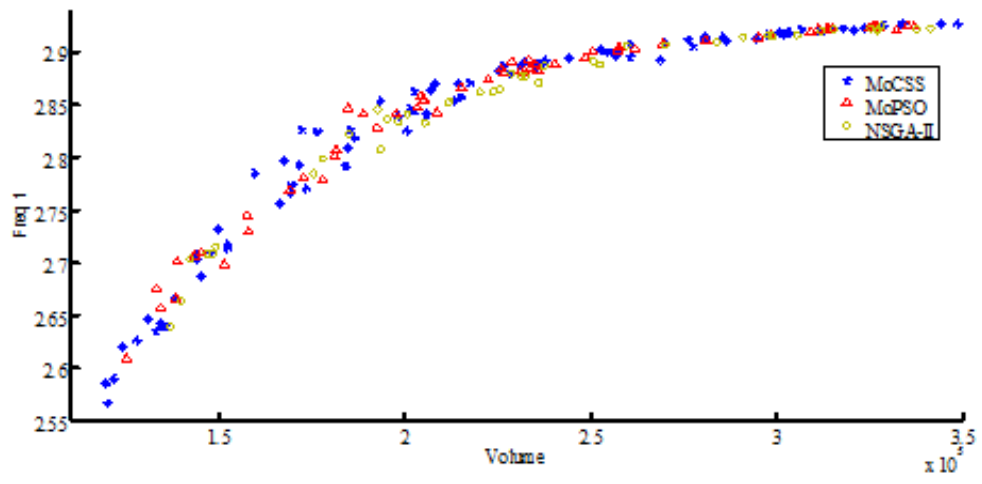

**(a)**

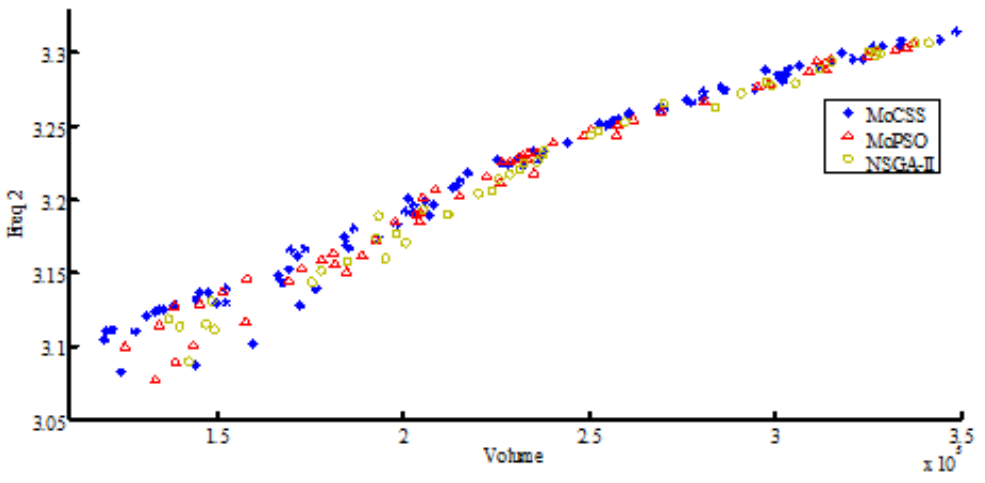

**(b)**

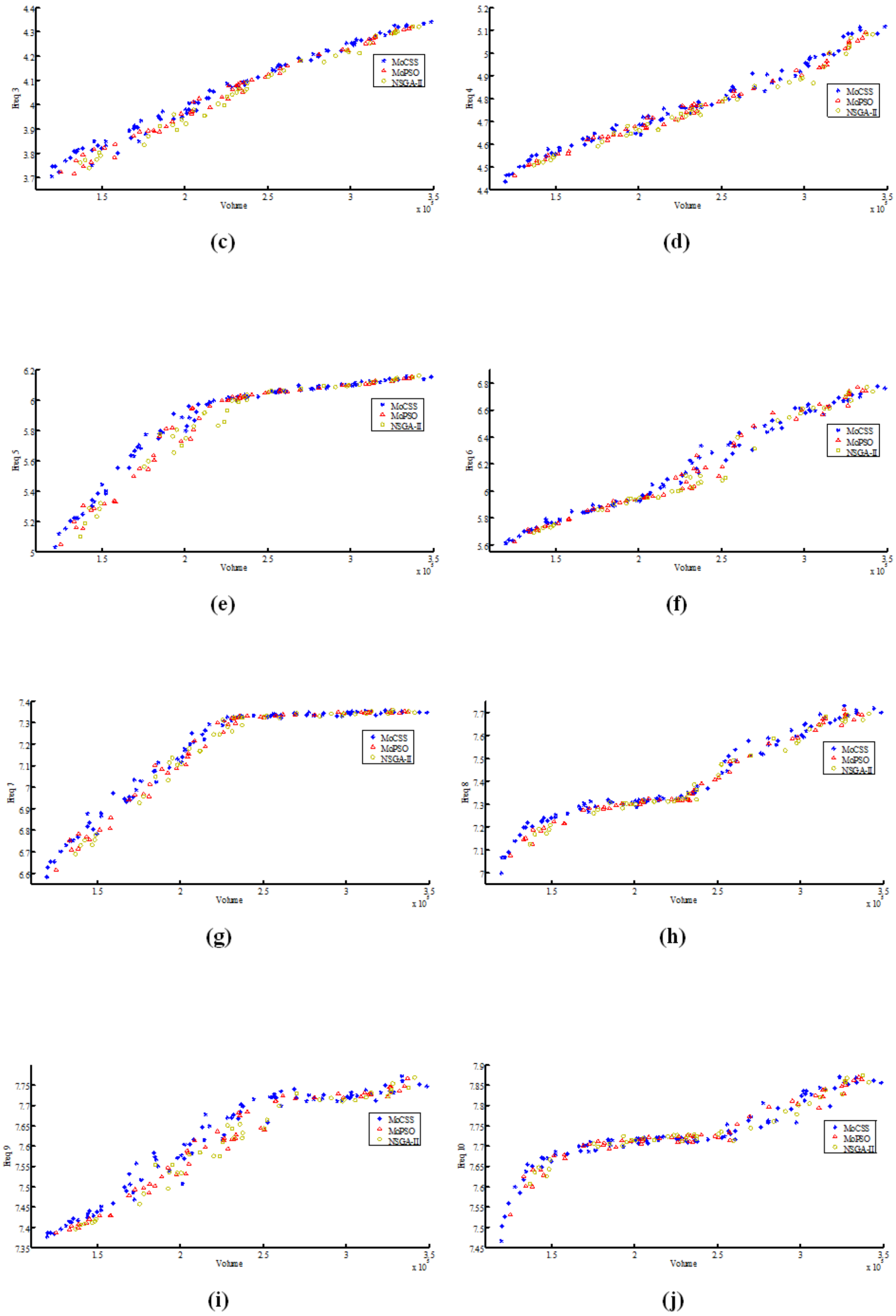

Fig. 9. The related Volume-Frequencies of three methods for P2

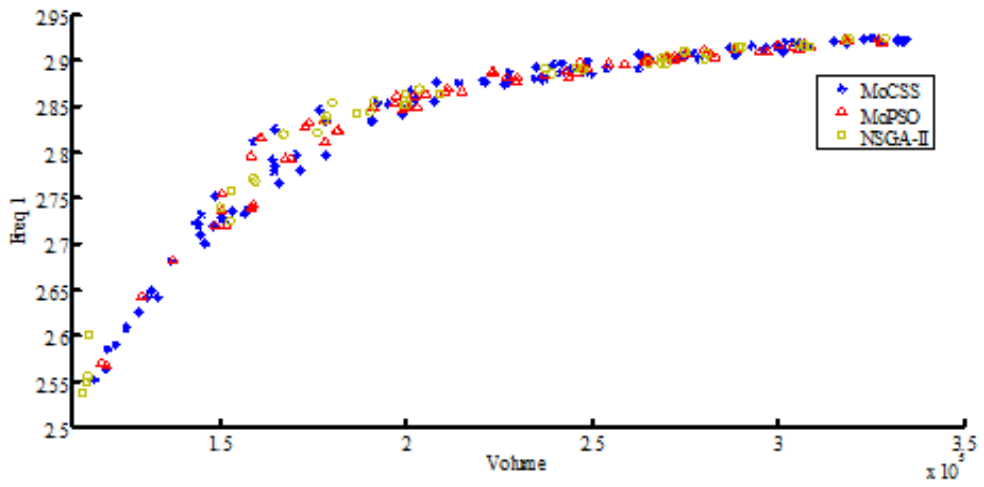

**(a)**

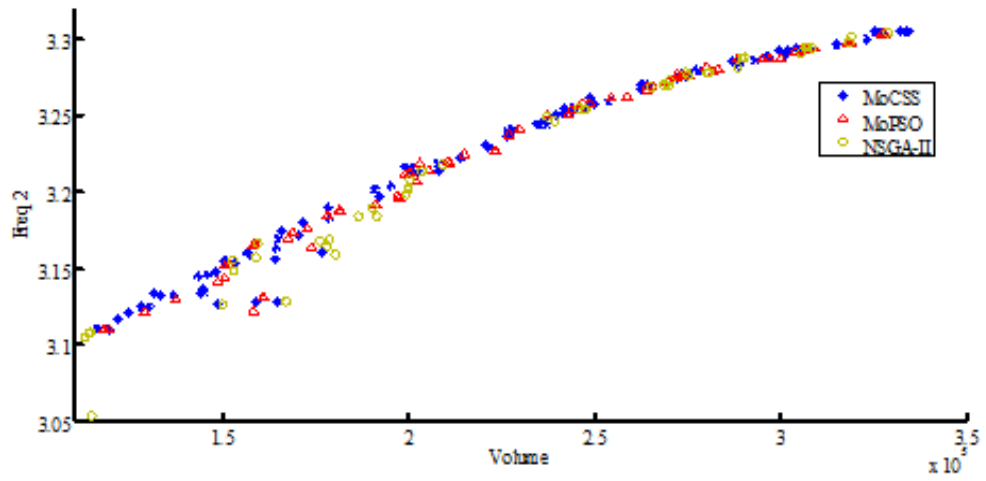

**(b)**

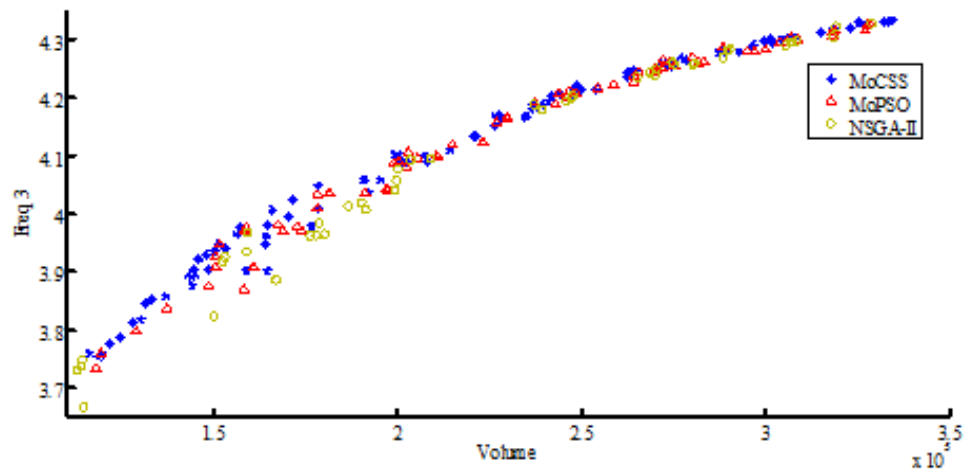
Freq 3
Volume
x 10^5
MoCSS
MoPSO
NSGA-II

**(c)**

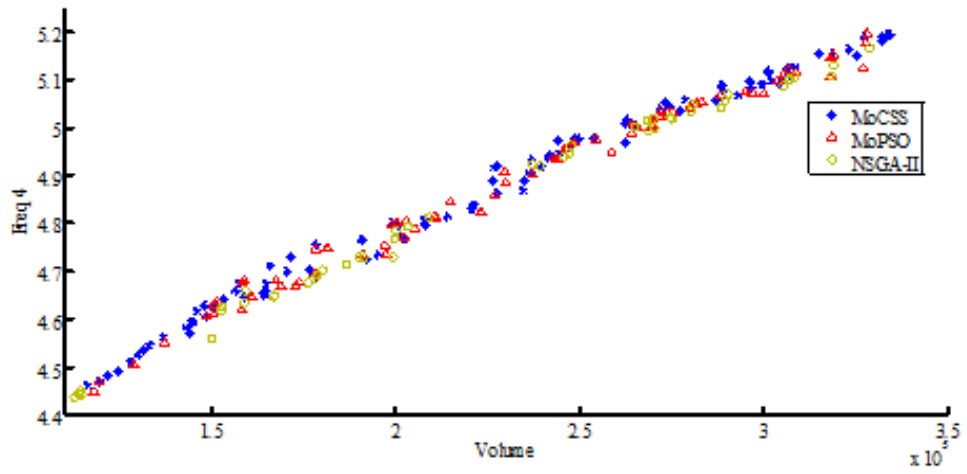
Freq 4
Volume
x 10^5
MoCSS
MoPSO
NSGA-II

**(d)**

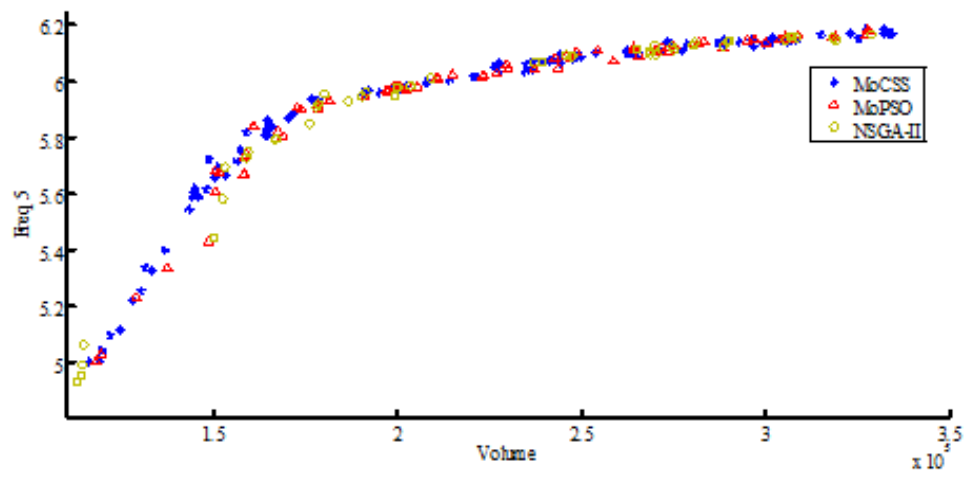
Freq 5
Volume
x 10^5
MoCSS
MoPSO
NSGA-II

**(e)**

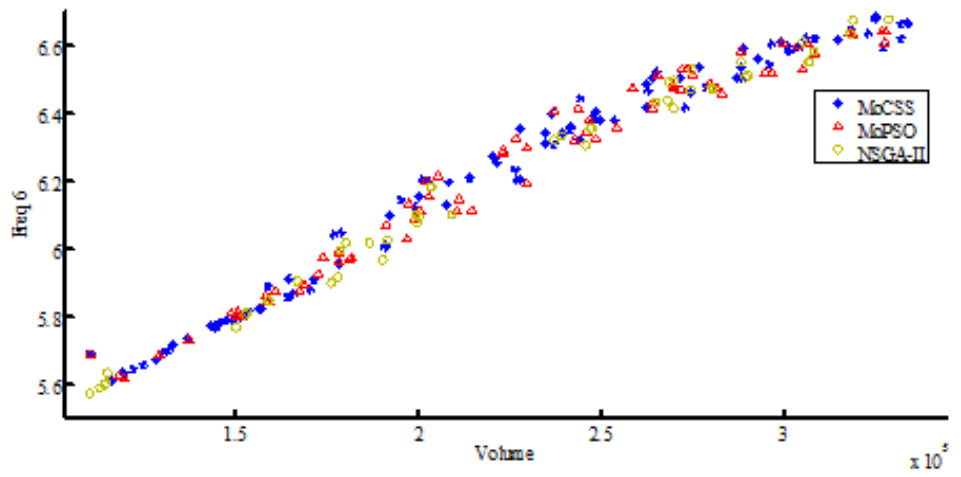
Freq 6
Volume
x 10^5
MoCSS
MoPSO
NSGA-II

**(f)**

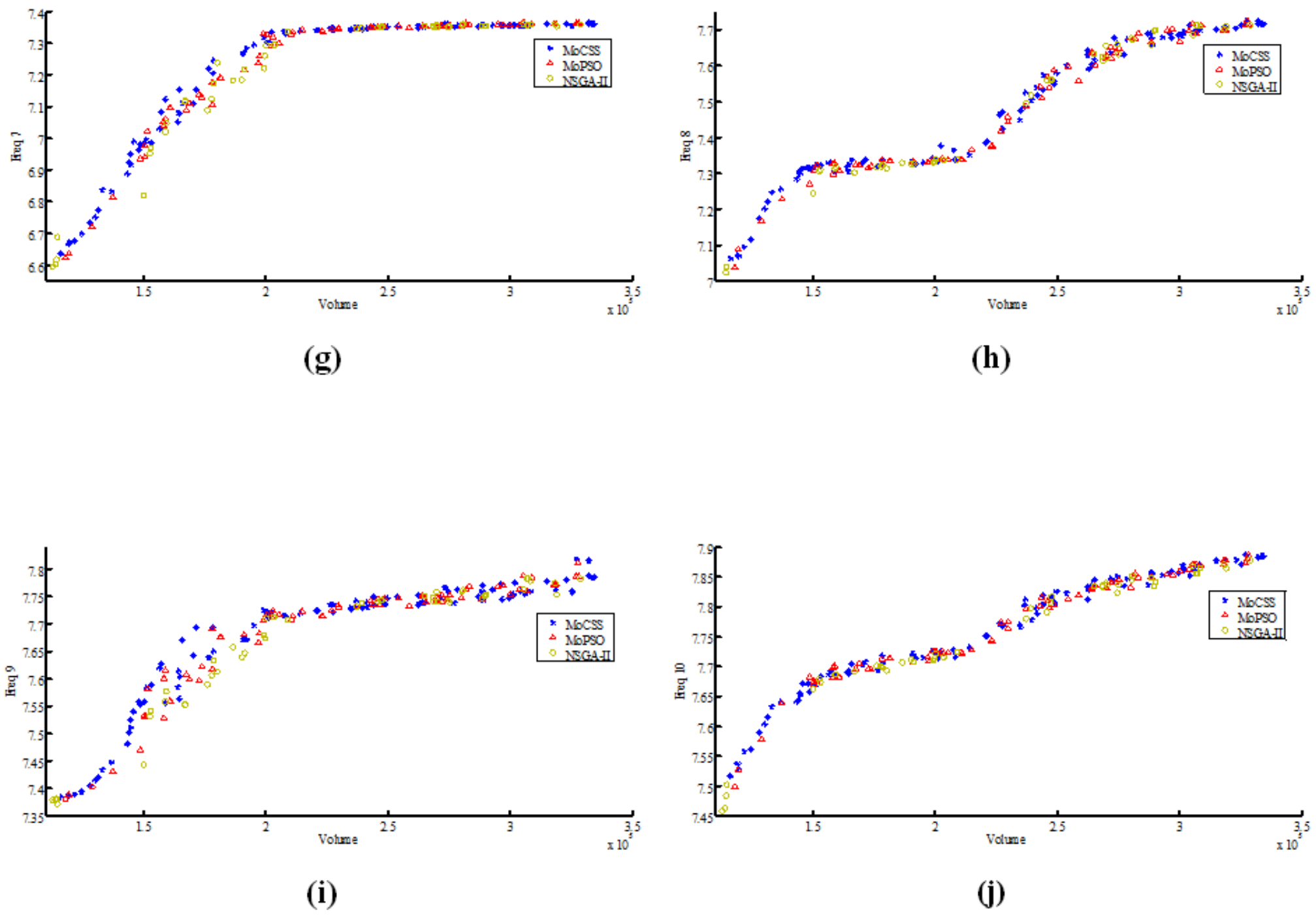

Fig. 10. The related Volume-Frequencies of three methods for P3

In order to have a proper comparison between the results and to exploit them according to decision maker criteria, five different scenarios with coefficients from 0.1 to 0.9 are considered for the all results. Tables 8, 9 and 10 illustrate the results for the P1, P2 and P3, respectively. In Table 8, the best volumes in scenarios A and C belong to the MoCSS whereas, best volume of scenario B belongs to MoPSO and best volumes of scenarios E and D belong to NSGA-II. However, if we consider the priority weights of objective functions as the comparison criterion (in the scenarios E and D, the volume in more important in the scenarios A and B, the frequencies are more important and in the scenario C, both have the equal importance),_best results of scenarios A, B, and C will belong to MoCSS and just the best result of scenario D and E will belong to NSGA-II. Same comparisons can be presented for the problems P2 and P3 as summarized in Tables 10 and 11. For the problem P2, the results of the MoCSS for scenarios A, B and C are better than those of other methods, while for problem P3, the MoCSS yields best results as shown in Table 11.

| |
|---|
| Table 8: Different possible scenarios for problems P1 corresponding solutions |
| Table 9: Different possible scenarios for problems P2 corresponding solutions |
| Table 10: Different possible scenarios for problems P3 corresponding solutions |

According to Fig. 11, if we plot the results of the MoCSS and different scenarios for the presented problems in one diagram, two striking  results will emerge. The first one is that results of P3 is better than two other problems and the second one describes that in different scenarios, the problem P3 gives less volumes with acceptable frequencies. To sum up, the most obvious finding to emerge from this study is that MoCSS outperforms its rivals and problem P3 that use multiplication of reversed frequencies as the second objective function, yields better results in comparison to two other problems.

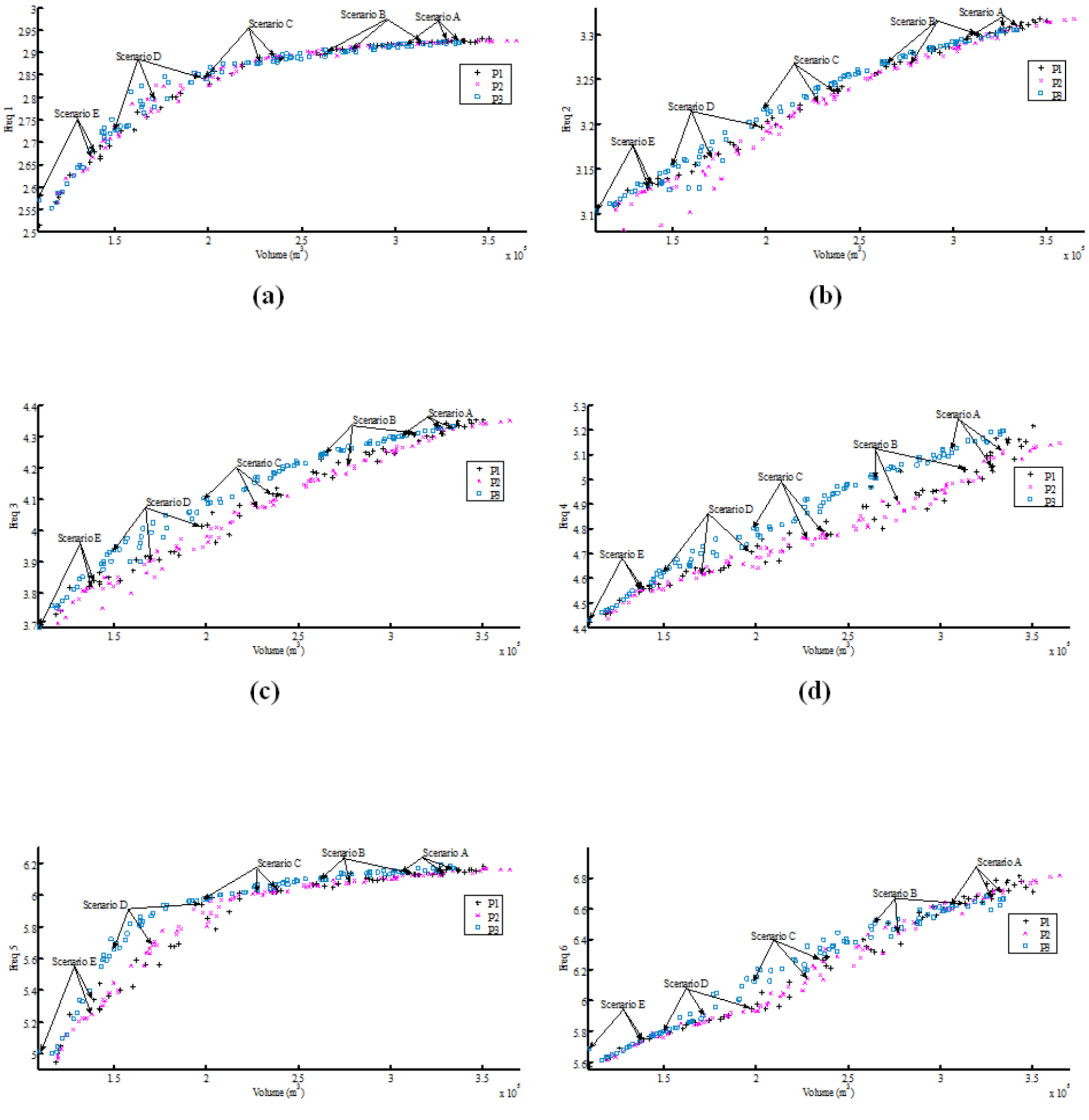

(a) (b)

(c) (d)

(e) (f)

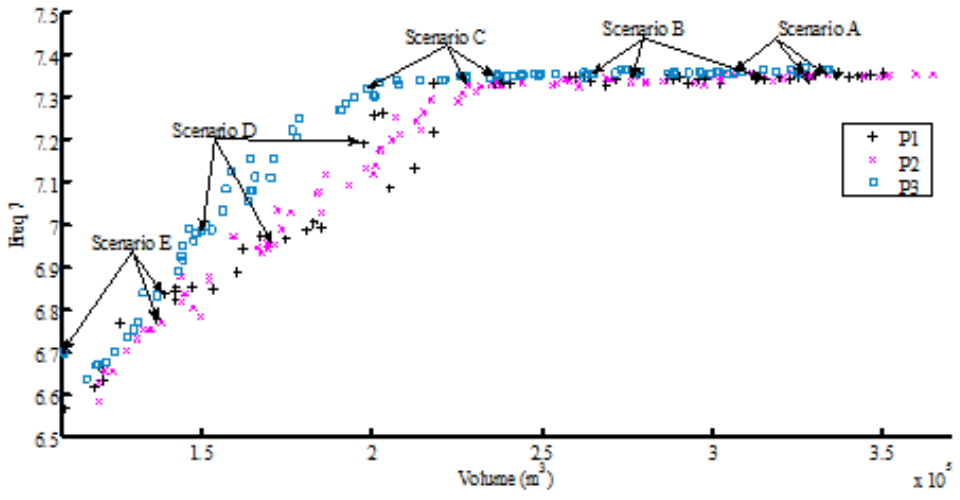

(g)

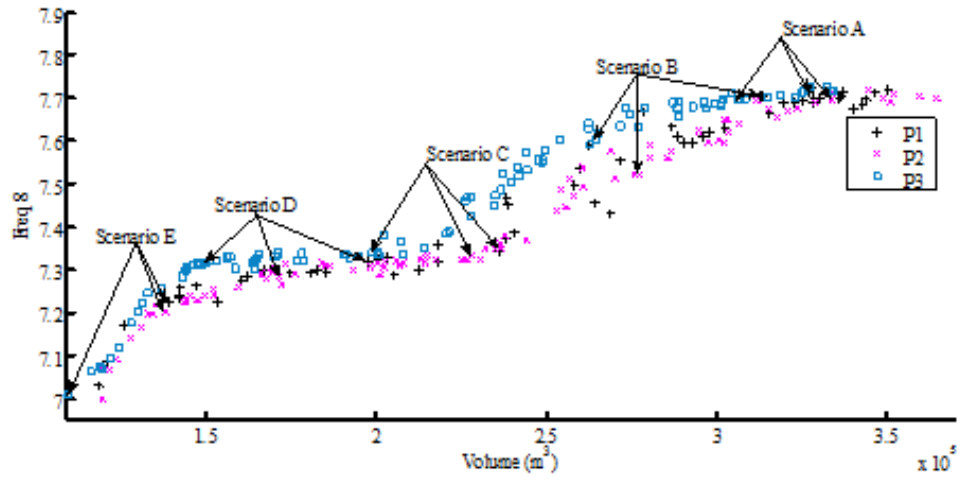

(h)

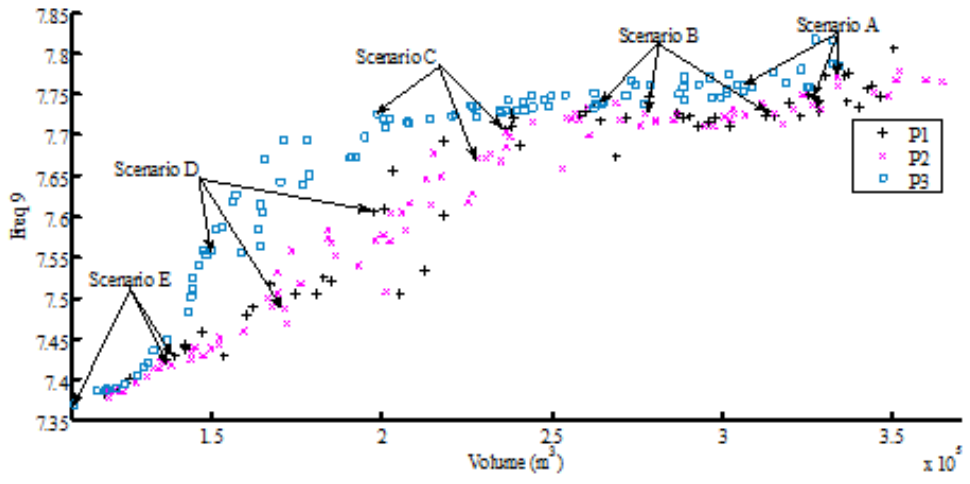

(i)

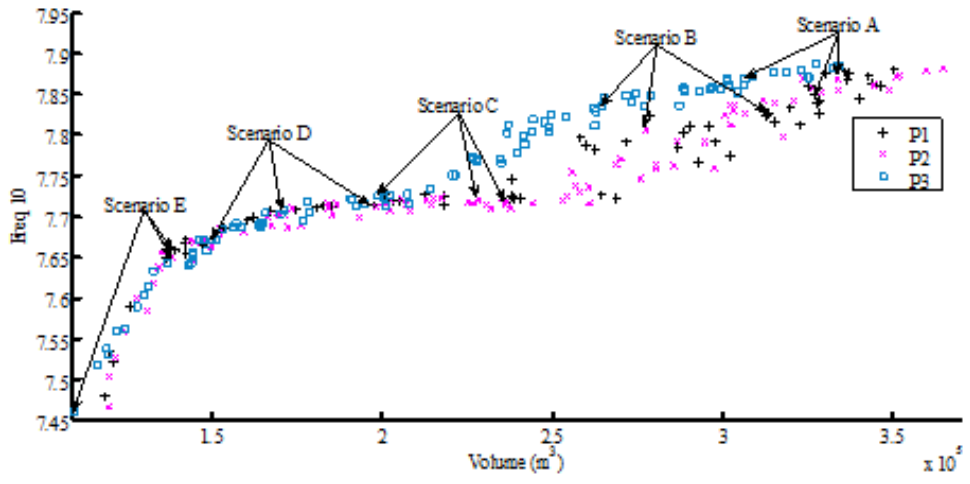

(j)

Fig. 11. Best solutions according to the five different scenarios for P1, P2 and P3 by MoCSS

## 5. Conclusion

In this paper, a new effective and robust methodology is proposed for finding the optimal shape of double-curvature arch dams. In contrast with previous researches, instead of solving a single-objective optimization, which yields only a single point as the answer, the output of proposed approach is a set of Pareto solutions in which decision maker can choose the best answer according to his preferences. The MoCSS is utilized for this purpose and developed by applying some additional steps and modifications in the standard CSS, which as a multi-agent optimization technique utilizes the governing laws from electrostatics and Newtonian mechanics to determine the movement process of the agents. A set of parallel-working APDL-MATLAB codes are developed for modeling, analyzing and obtaining fitness functions of the arch dam. These set of codes interface with the MoCSS in every step of optimization process.

In this study, the concrete volume is considered as the first objective function and a combination of reversed natural frequencies of the dam are considered as second objective function. For this purpose, three different problems were defined. Separate reversed frequencies as well as sum and multiplication of them are considered as second objective functions in P1, P2 and P3 respectively. The reservoir is assumed full and interaction of dam reservoir is considered. In order to evaluate the effectiveness of proposed methodology, the Morrow Point dam as a well-known and real-world arch dam has been investigated. Prior to implementation of the proposed algorithm, a finite element model of the dam is developed and verified with existing results from literature. The resulted Pareto solution for problems P1, P2 and P3 from the MoCSS is compared with those of the MoPSO and NSGA-II approaches. In addition, extreme values of every method for each problem are extracted and compared. Five different scenarios are defined to apply the preferences of decision maker to choose best answers. Solutions for these scenarios are compared for three method and three optimization problems. Numerical results indicate that MoCSS outperforms its rivals and P3 yields best results between three problems. Besides, results reveal that the proposed methodology is highly competitive with existing arch dam optimization approaches.

List of figures:

List of tables:

(a)

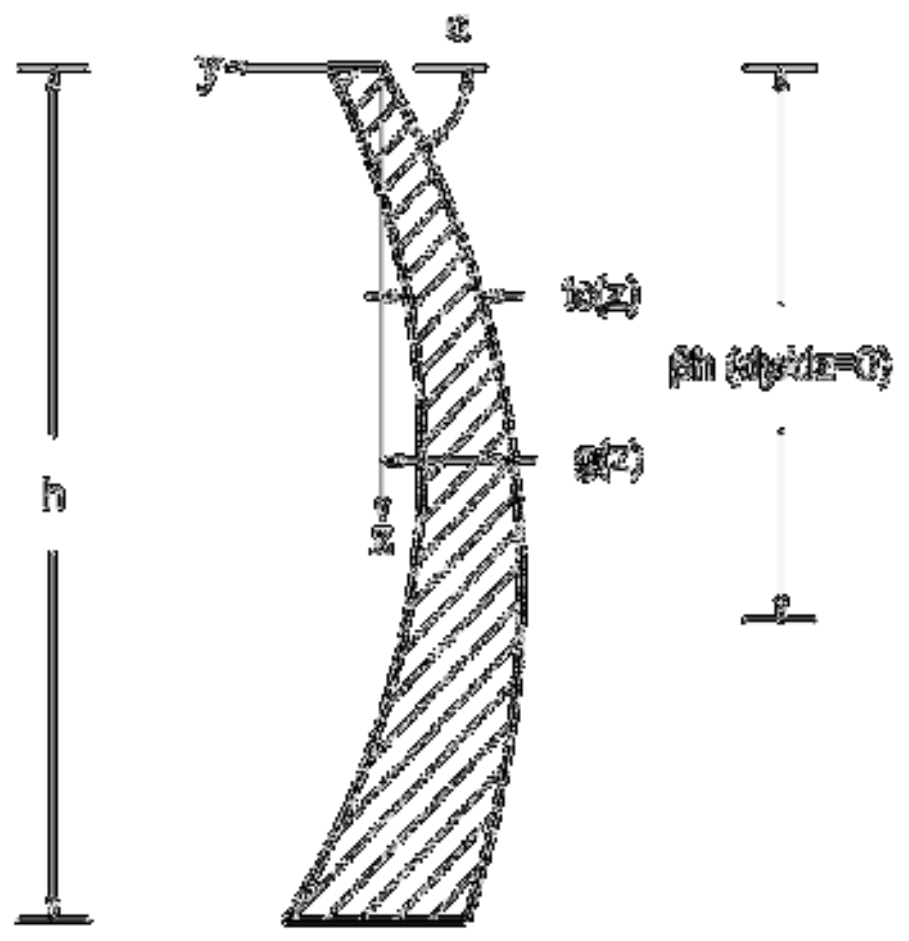

(b)

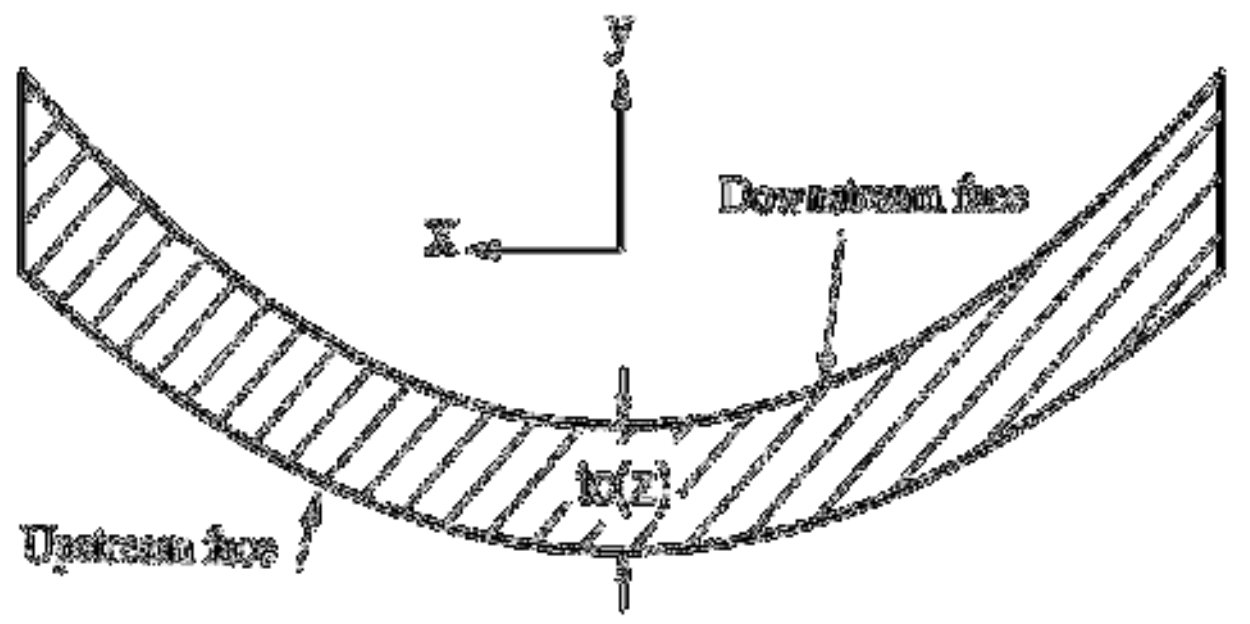

Fig. 1. Schematic of arch dam: a) Crown Cantilever profile; b) parabolic shape of an elevation

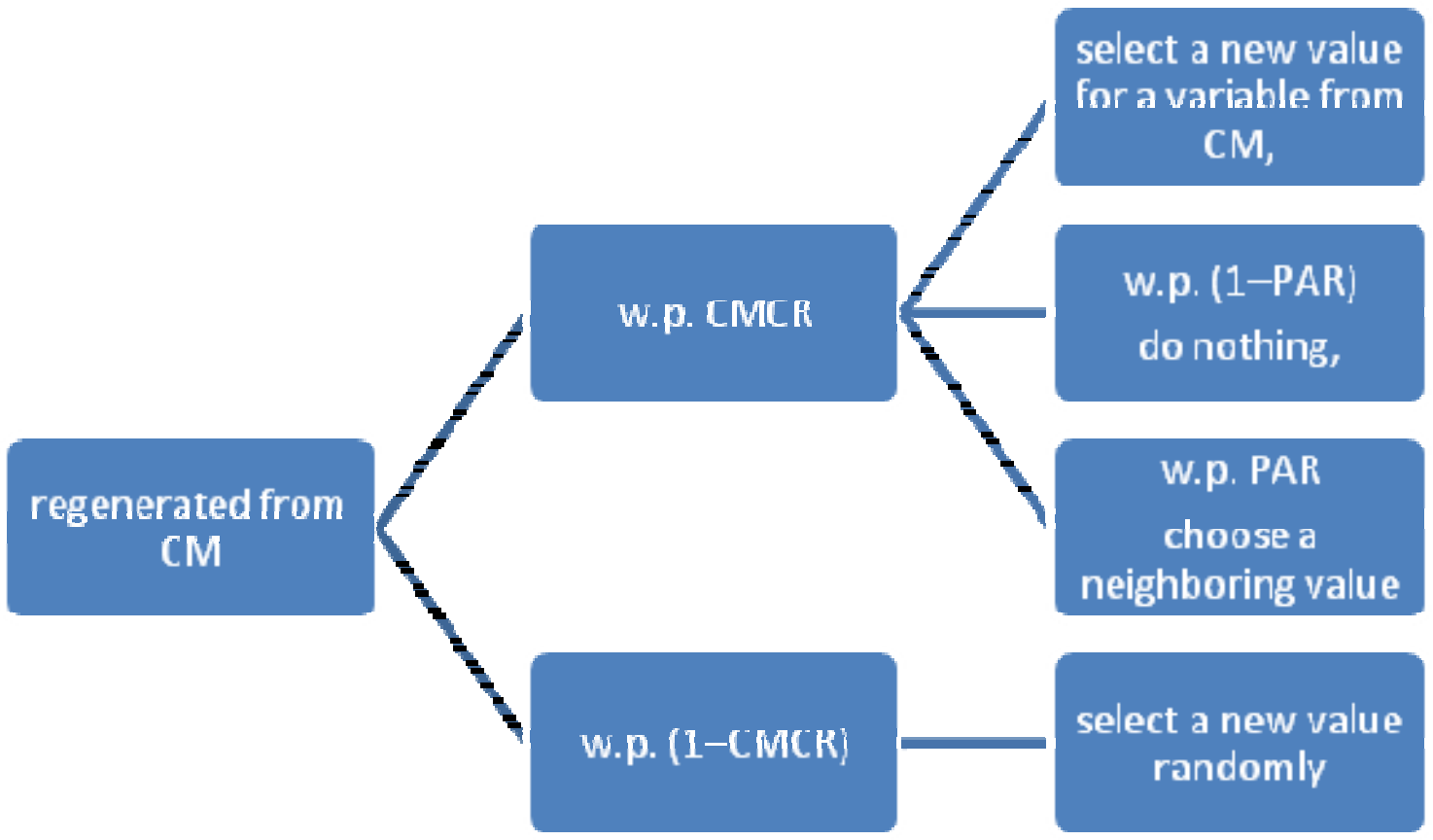

Fig. 2. Regenerating violated solutions from CM

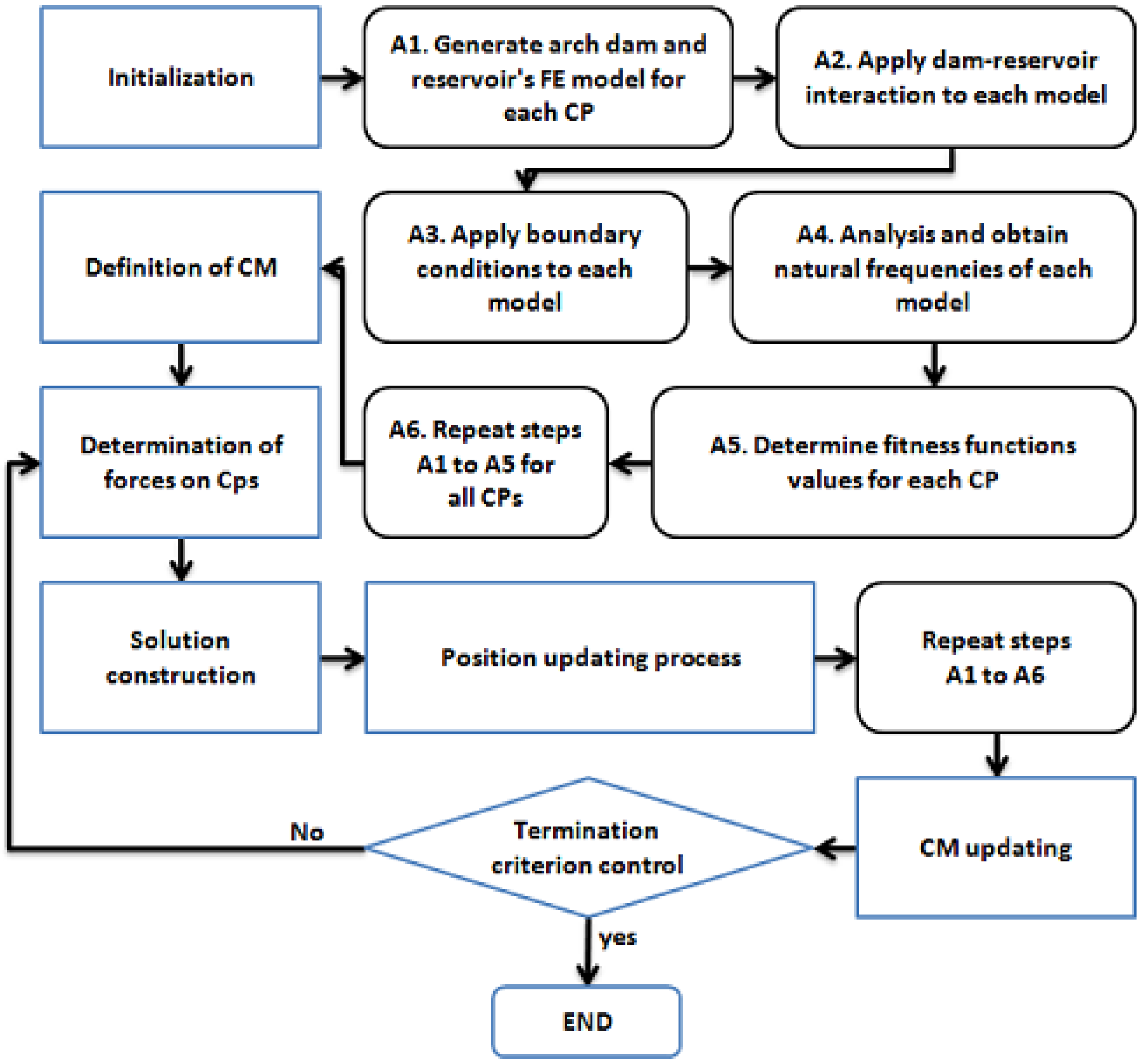

Fig. 3. Summarized flowchart of proposed method

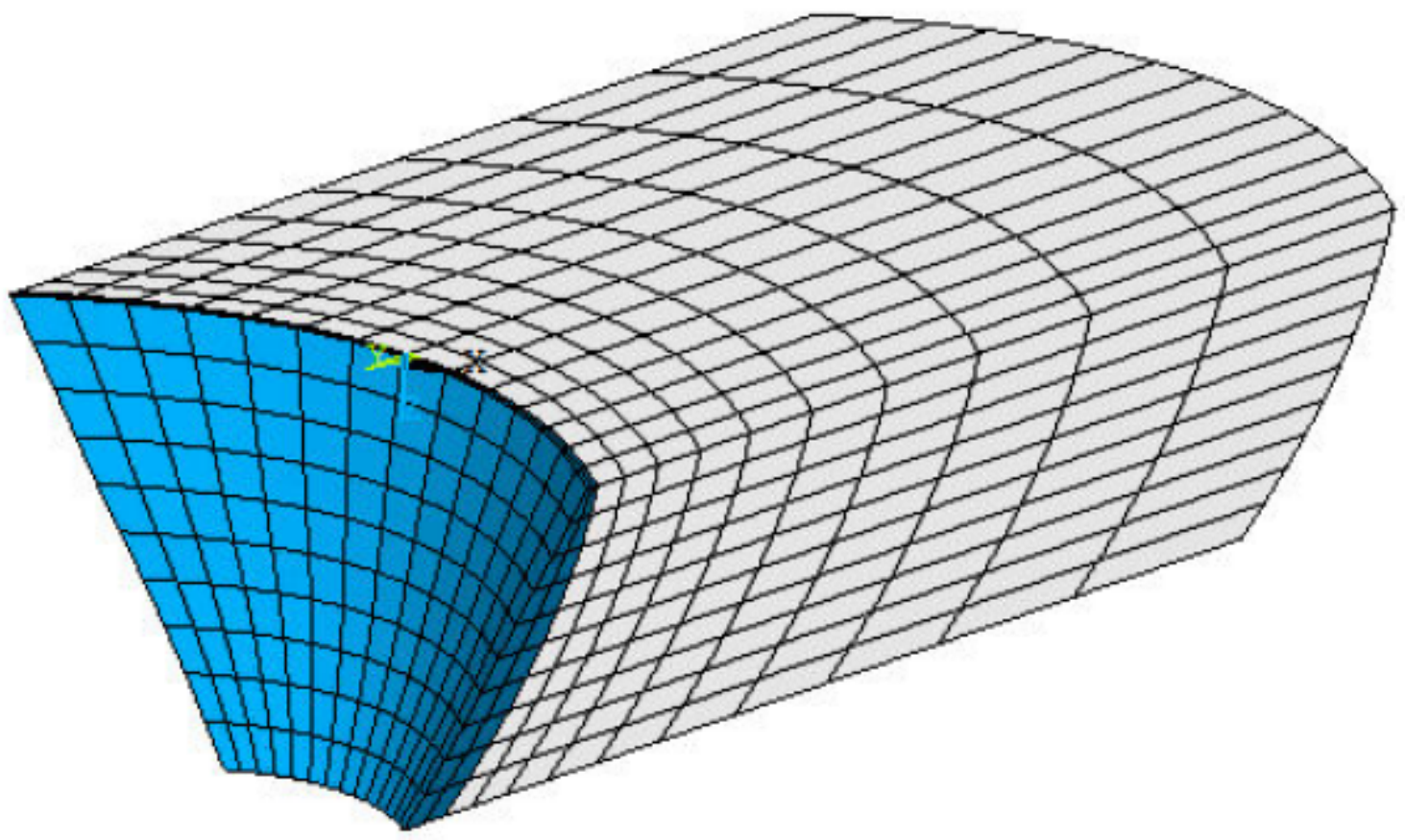

Fig. 4. Finite Element Model of Morrow Point Dam

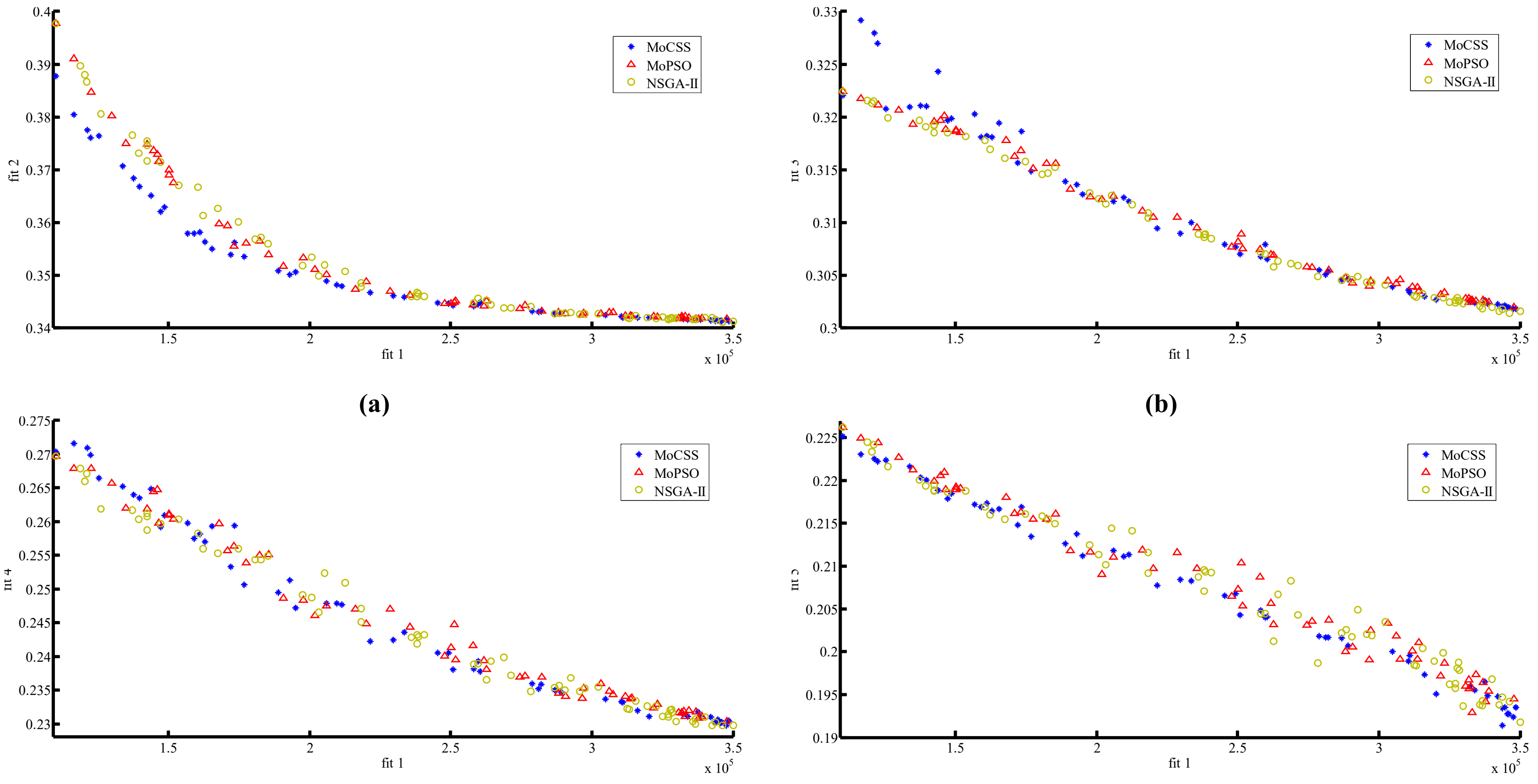
MoCSS
MoPSO
NSGA-II
fit 1
fit 2
x 10^5
1.5
2
2.5
3
3.5
0.34
0.35
0.36
0.37
0.38
0.39
0.4
0.3
0.305
0.31
0.315
0.32
0.325
0.33
0.23
0.235
0.24
0.245
0.25
0.255
0.26
0.265
0.27
0.275
0.19
0.195
0.2
0.205
0.21
0.215
0.22
0.225

**(a)** **(b)**

**(c)** **(d)**

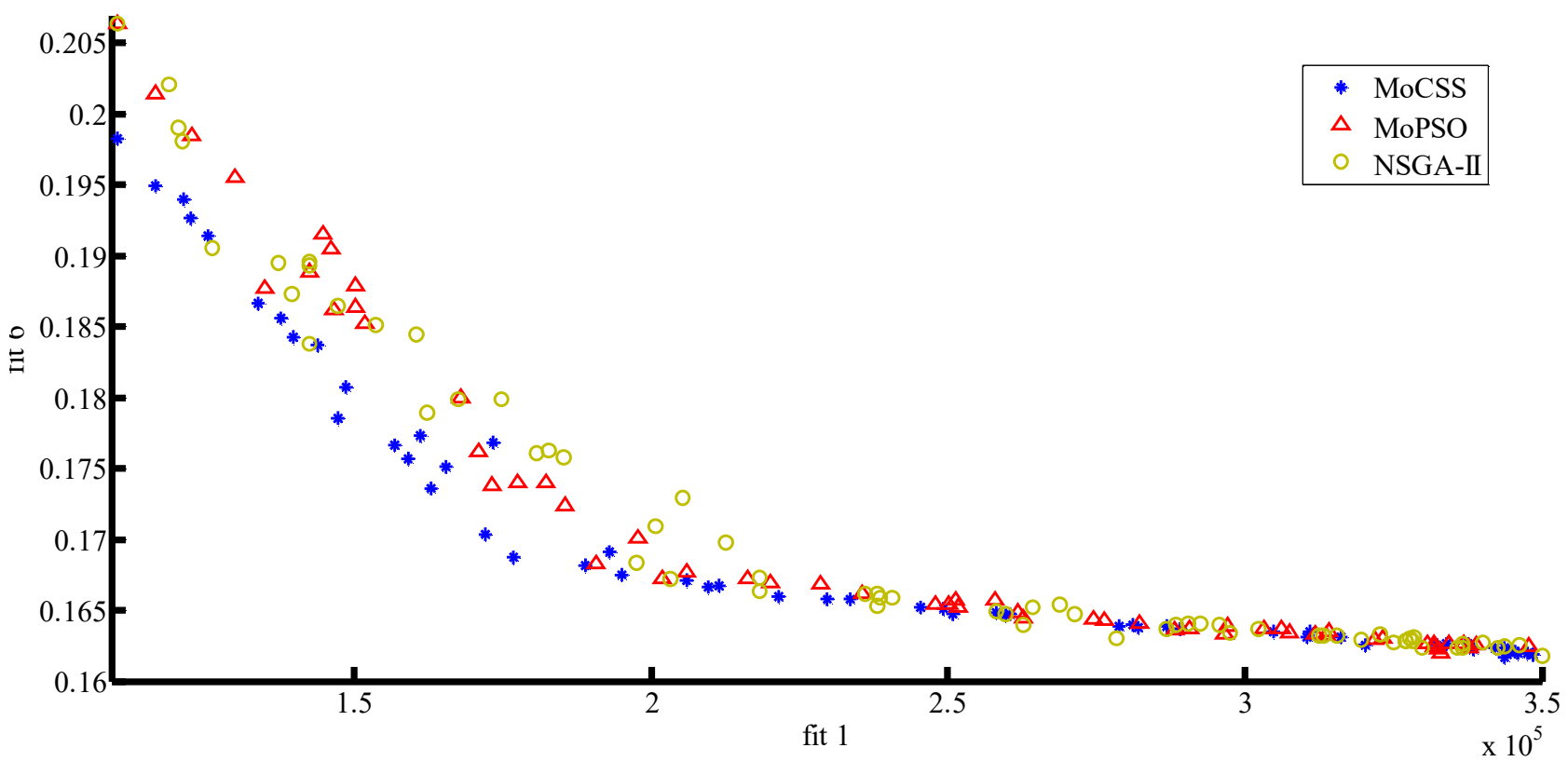
MoCSS
MoPSO
NSGA-II
fit 1
x 10^5

**(e)**

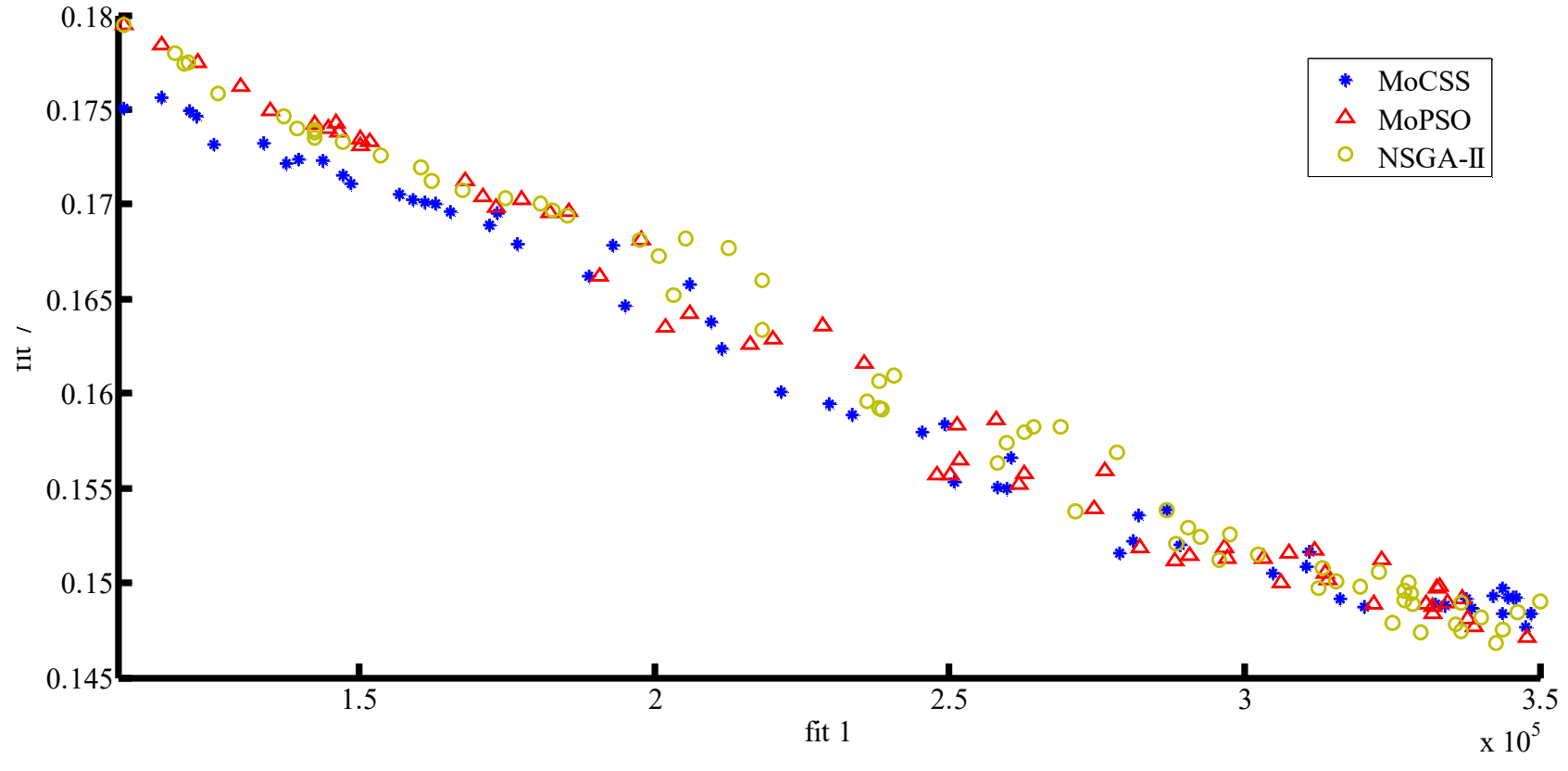
MoCSS
MoPSO
NSGA-II
fit 1
x 10^5

**(f)**

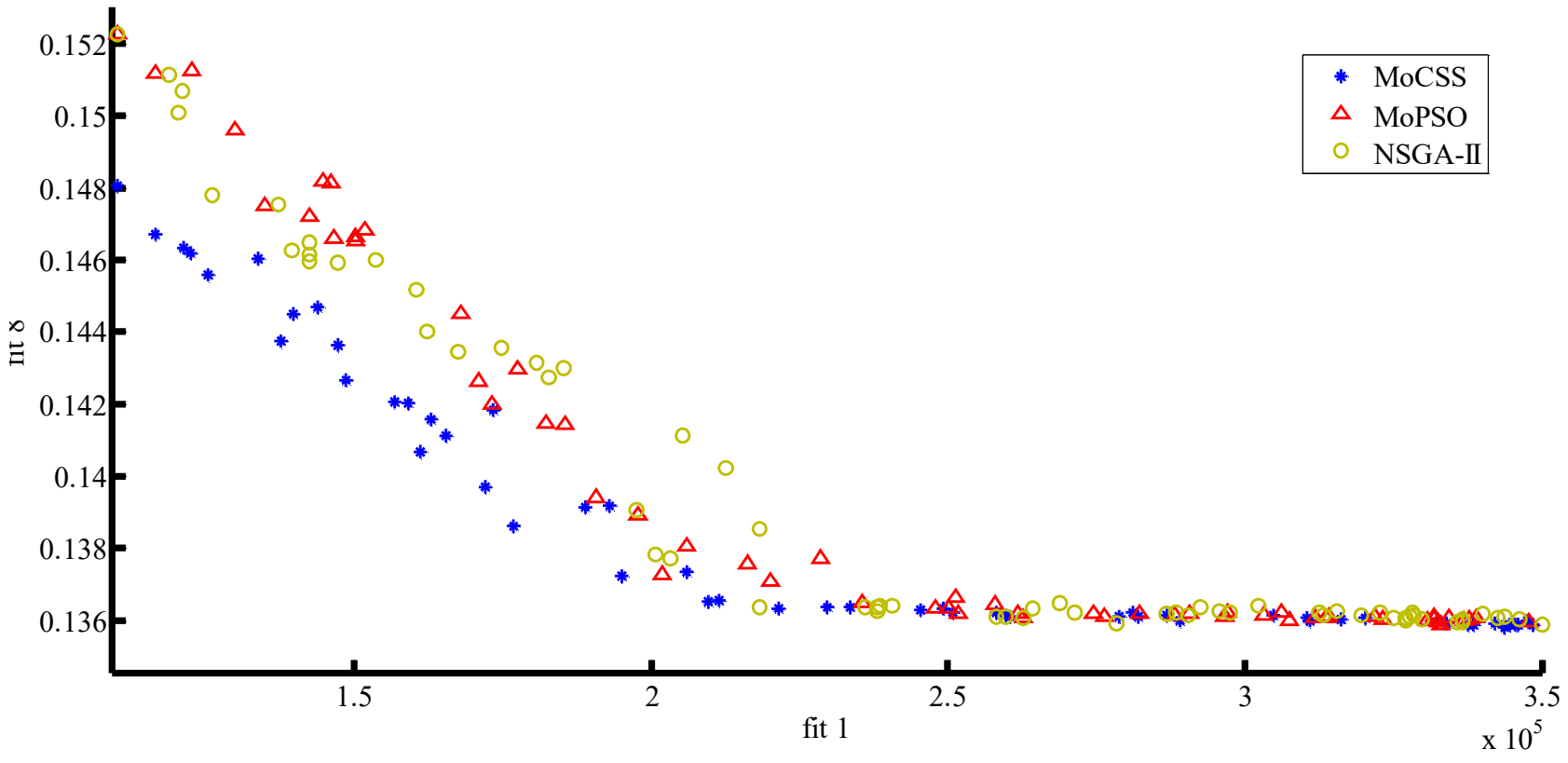
MoCSS
MoPSO
NSGA-II
fit 1
x 10^5

**(g)**

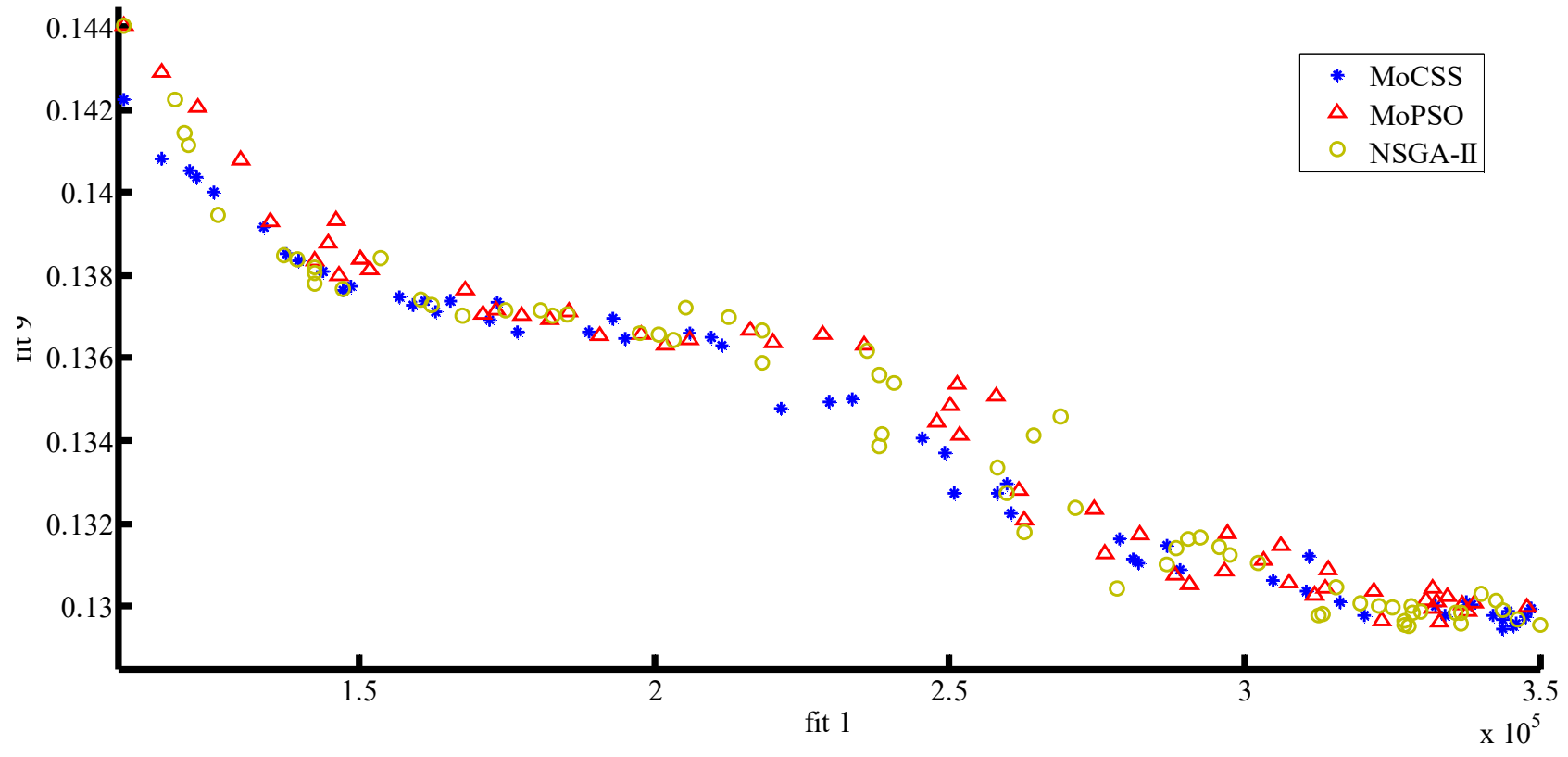
MoCSS
MoPSO
NSGA-II
fit 1
x 10^5

**(h)**

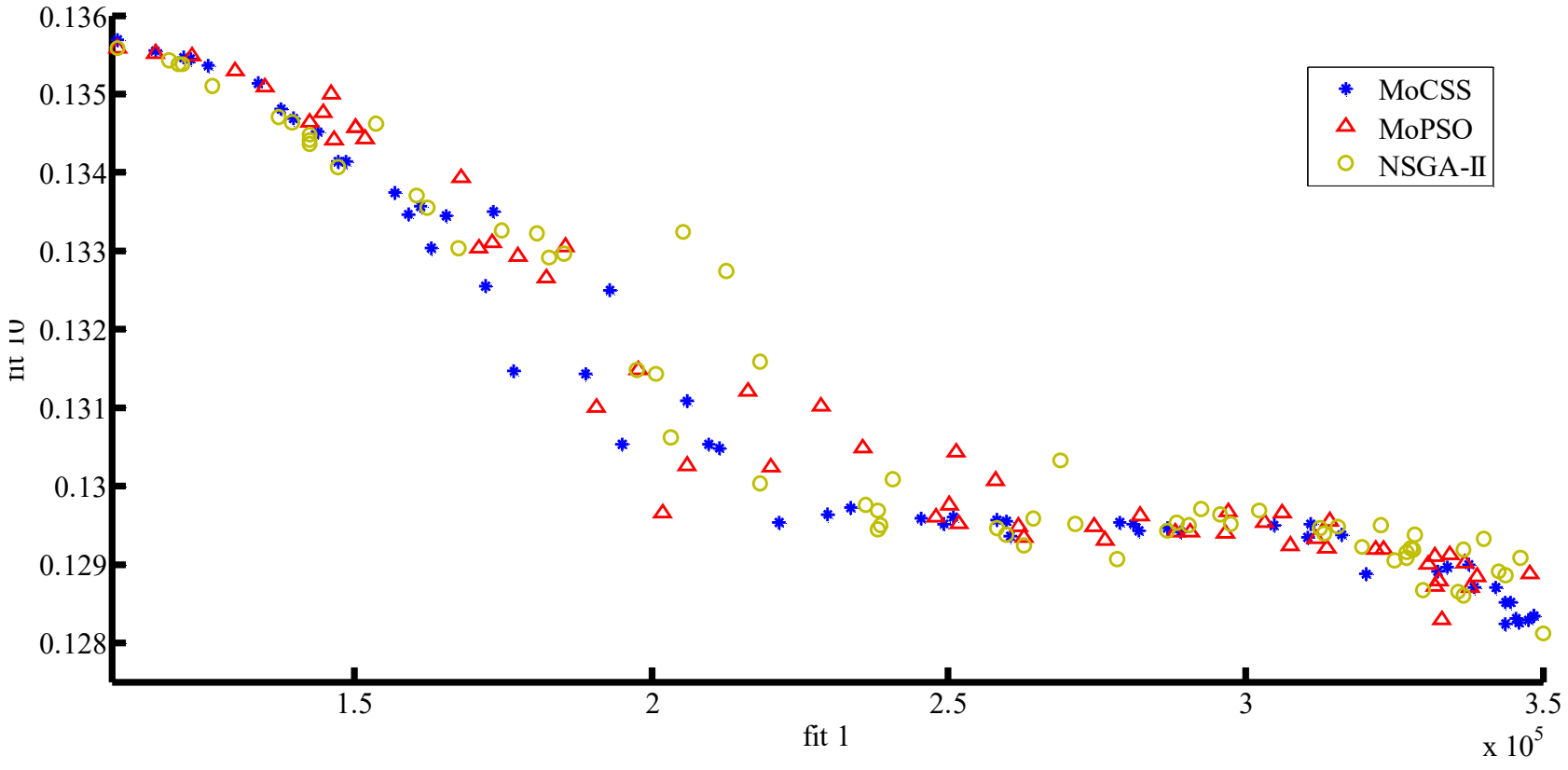

**(i)**

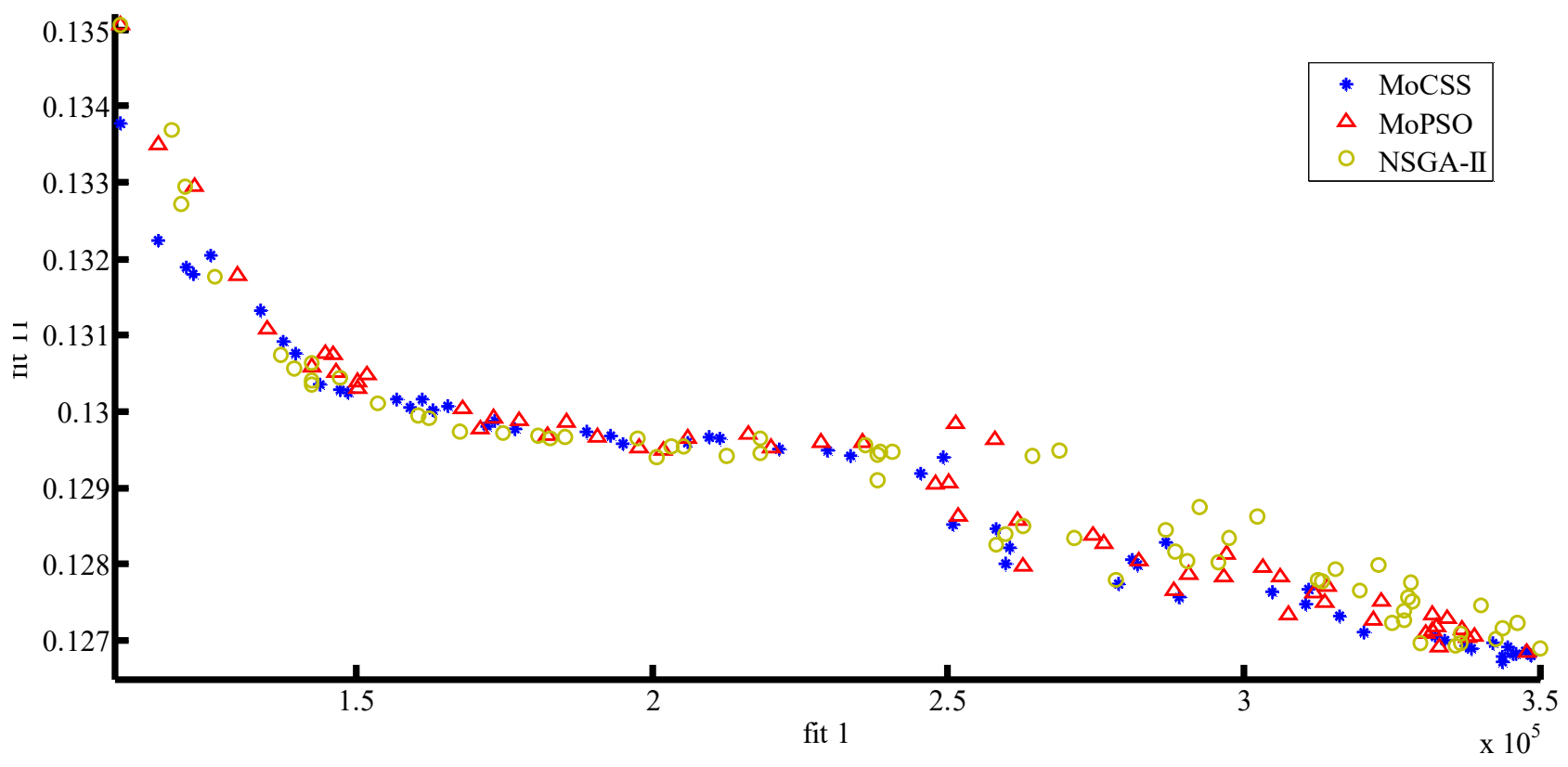

**(j)**

Fig. 5. Pareto fronts of three methods for P1; **(a)** fit1-fit2, **(b)** fit1-fit3**,** **(c)** fit1-fit4, **(d)** fit1-fit5**,** **(e)** fit1-fit6**,** **(f)** fit1-fit7**,** **(g)** fit1-fit8**,** **(h)** fit1-fit9**,** **(i)** fit1-fit10**,** **(j)** fit1-fit11**.**

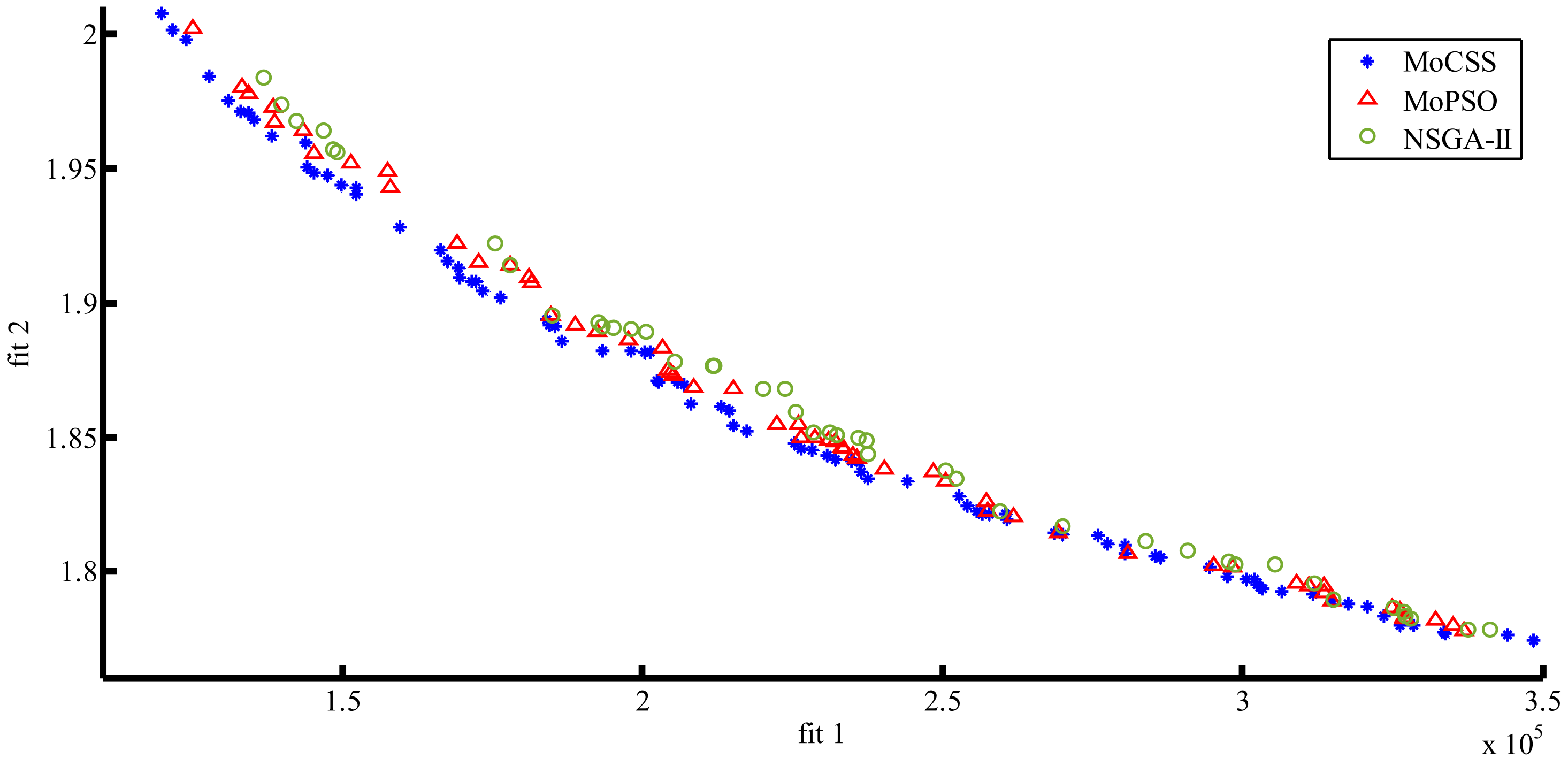

Fig. 6. Pareto front of three methods for P2

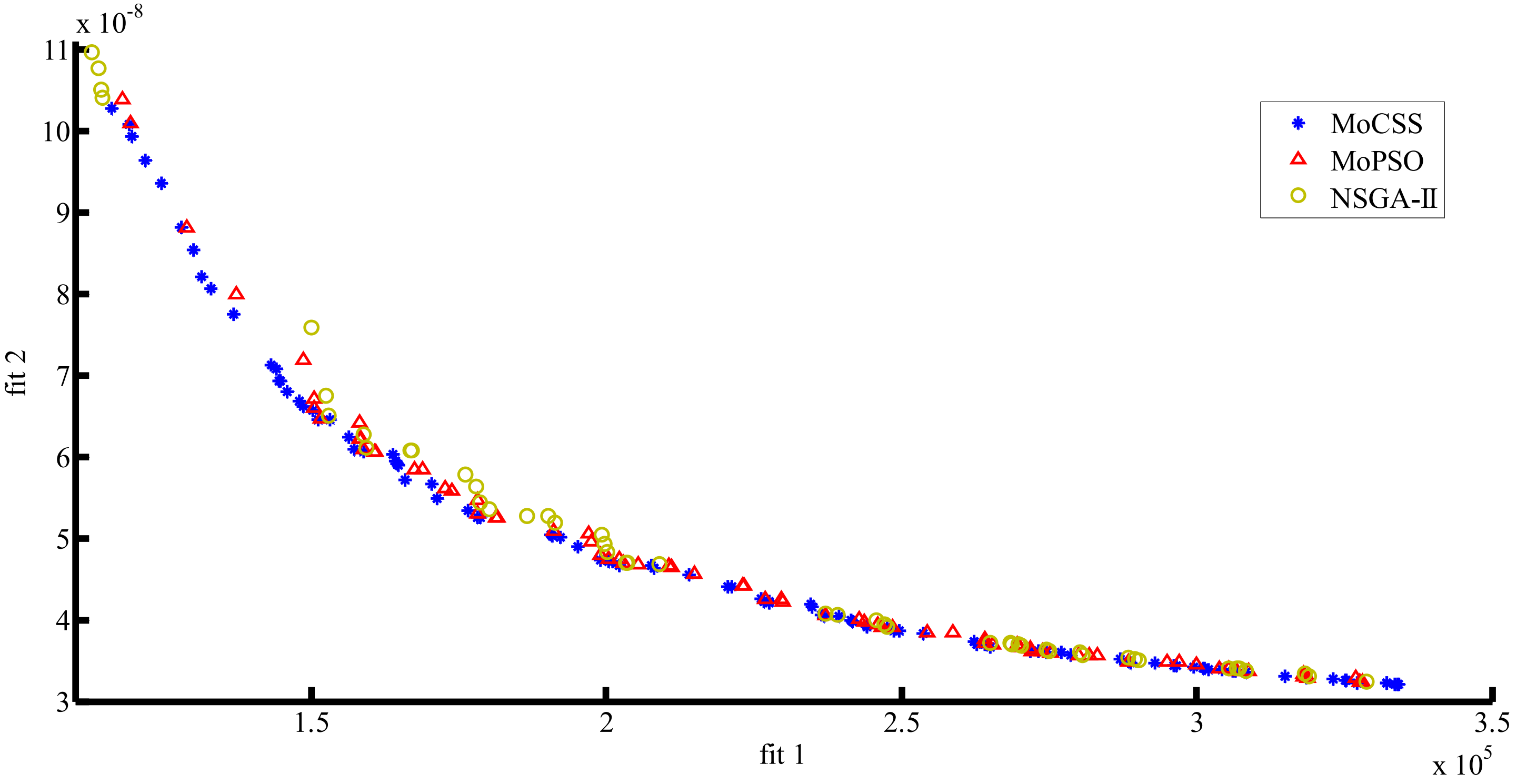

Fig. 7. Pareto front of three methods for P3

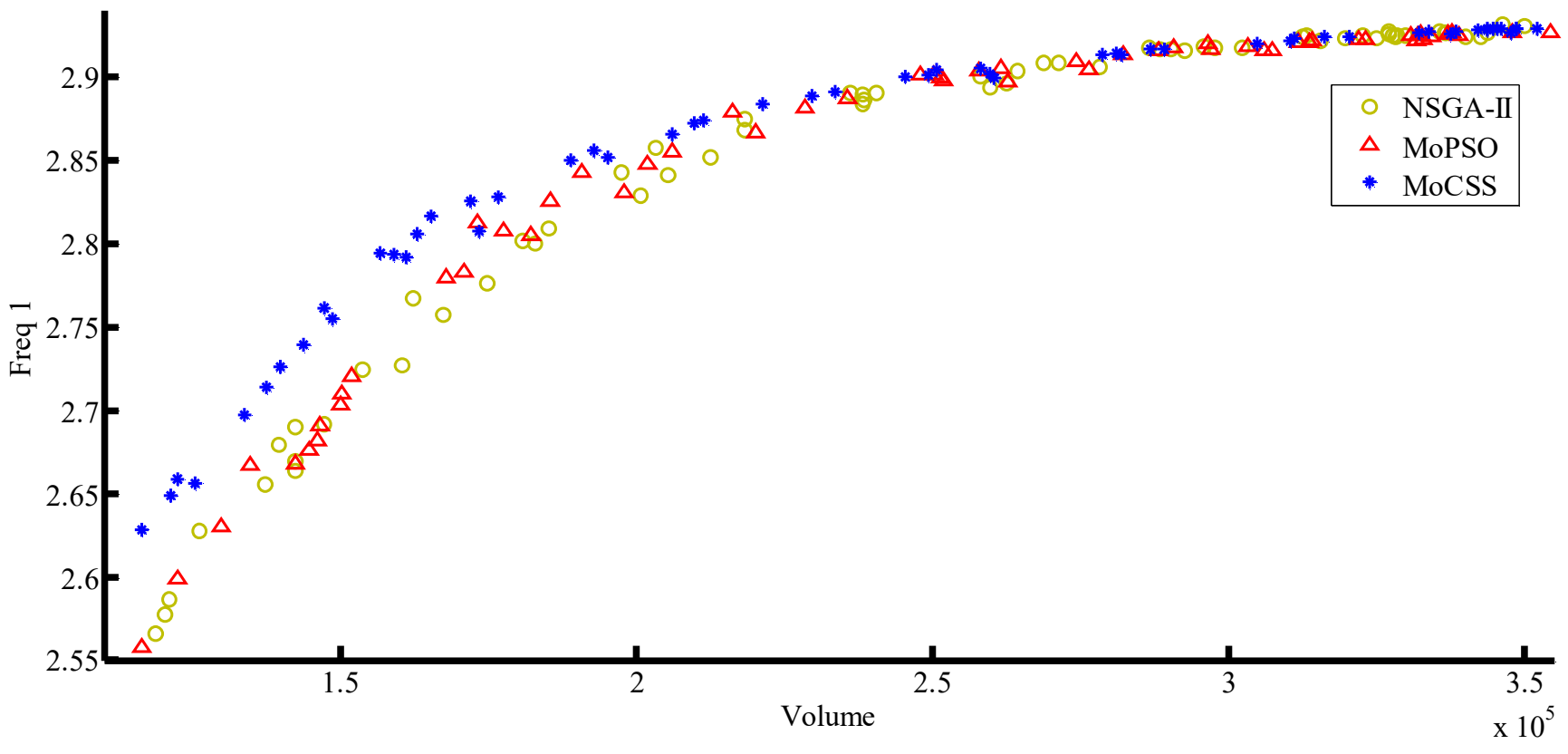
NSGA-II
MoPSO
MoCSS
Freq 1
Volume
x 10^5
2.55
2.6
2.65
2.7
2.75
2.8
2.85
2.9
1.5
2
2.5
3
3.5

**(a)**

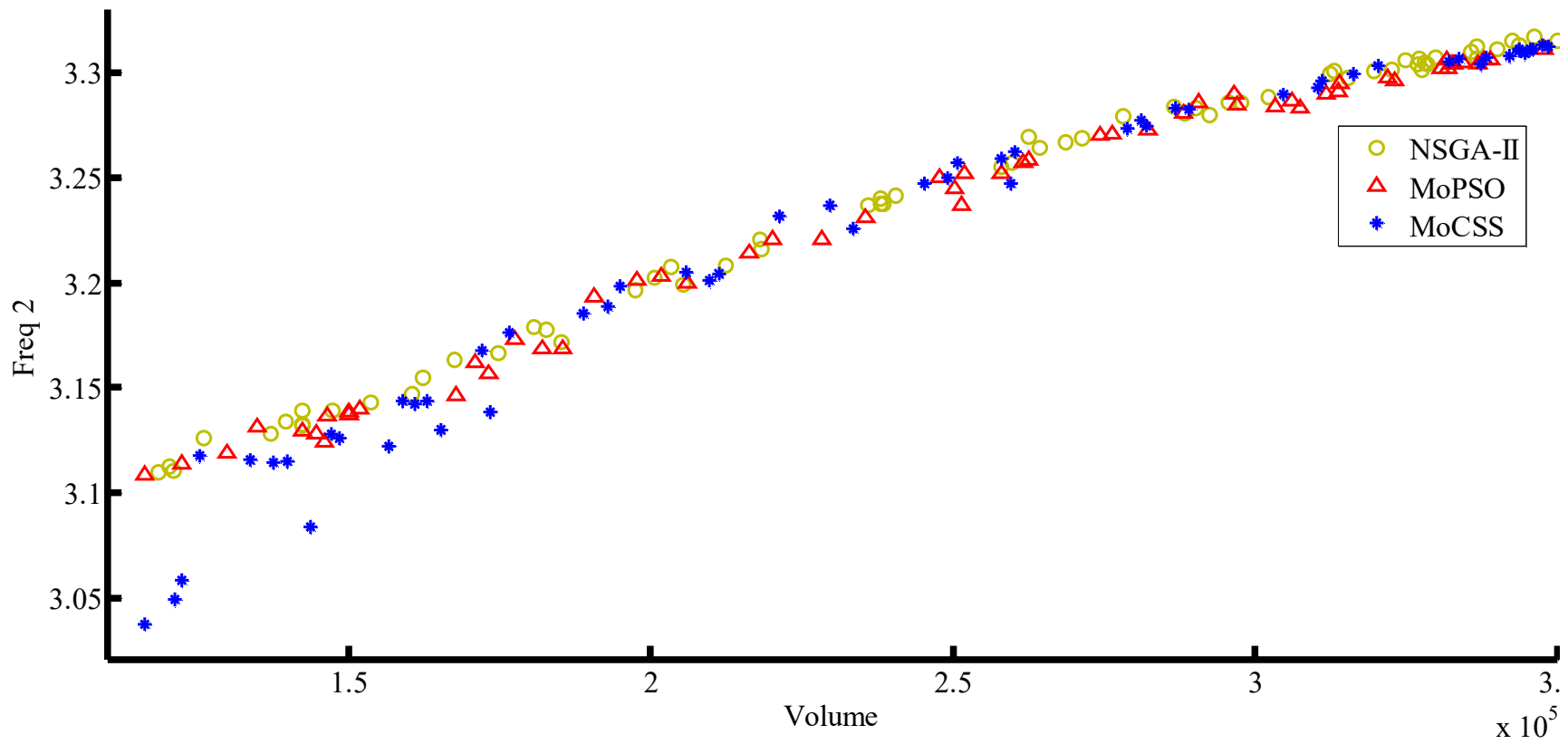
NSGA-II
MoPSO
MoCSS
Freq 2
Volume
x 10^5
3.05
3.1
3.15
3.2
3.25
3.3
1.5
2
2.5
3
3.

**(b)**

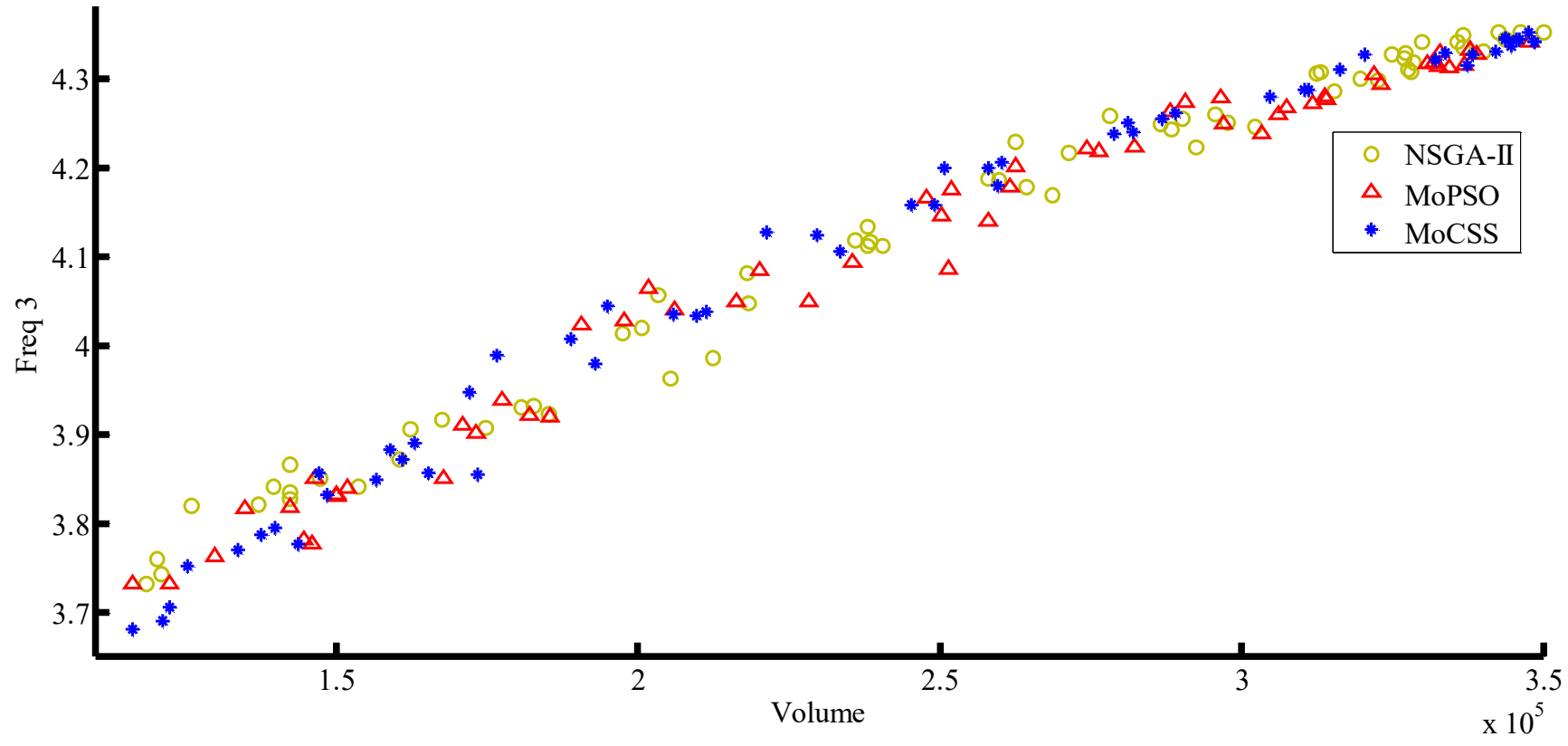
NSGA-II
MoPSO
MoCSS
Freq 3
Volume
x 10^5
3.7
3.8
3.9
4
4.1
4.2
4.3
1.5
2
2.5
3
3.5

**(c)**

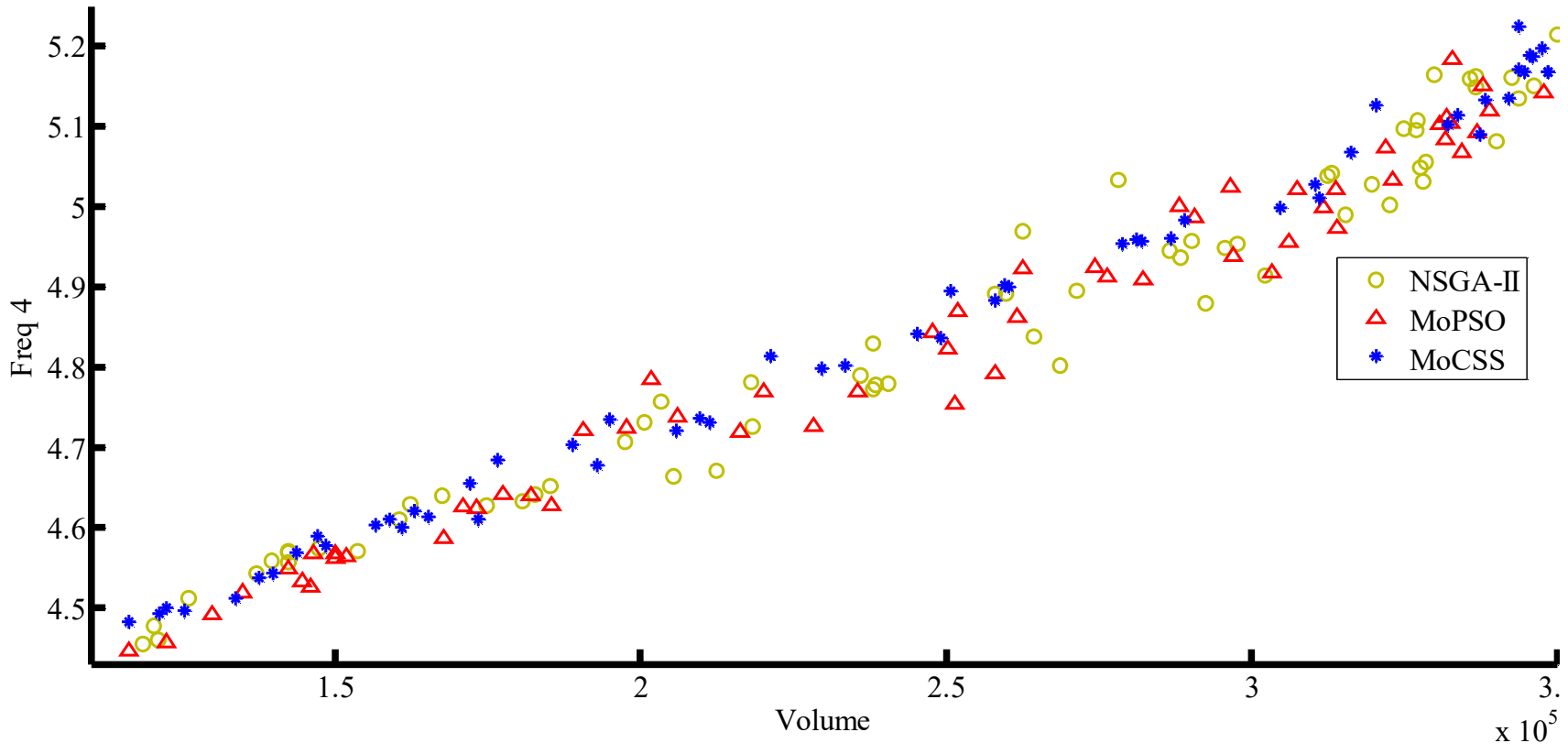
NSGA-II
MoPSO
MoCSS
Freq 4
Volume
x 10^5
4.5
4.6
4.7
4.8
4.9
5
5.1
5.2
1.5
2
2.5
3
3.

**(d)**

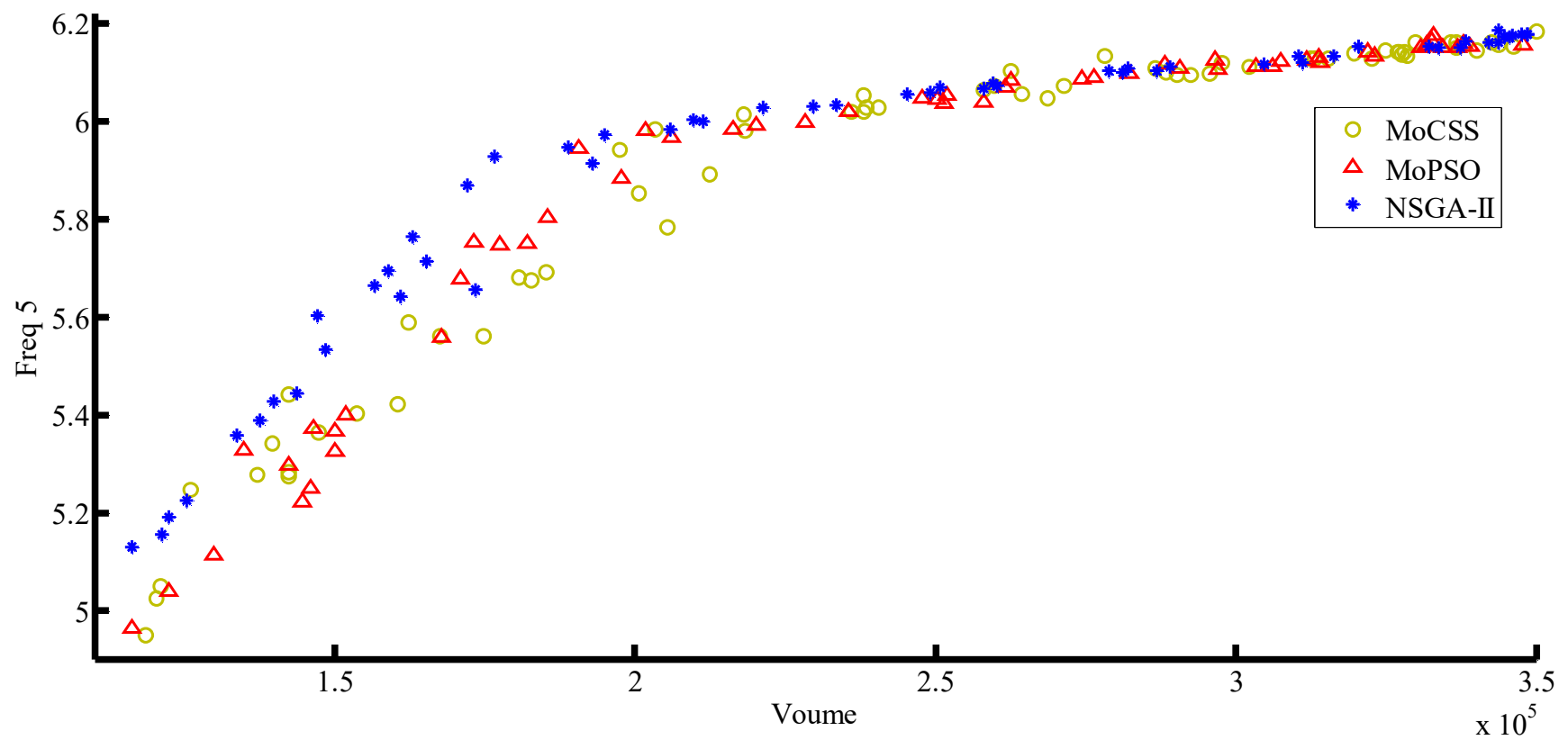

6.2
6
5.8
5.6
5.4
5.2
5
Freq 5
1.5
2
2.5
3
3.5
Voume
x 10^5
MoCSS
MoPSO
NSGA-II

**(e)**

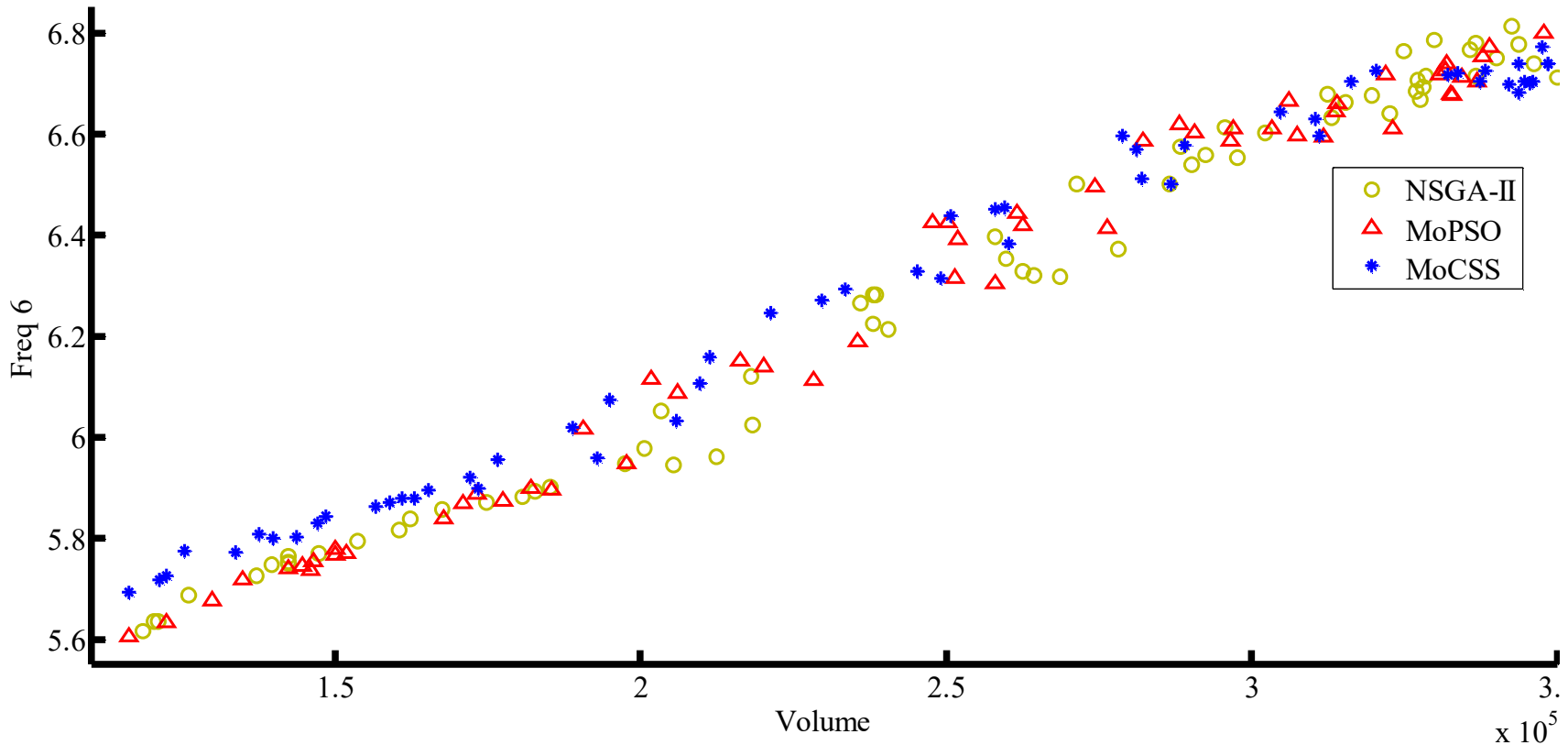

6.8
6.6
6.4
6.2
6
5.8
5.6
Freq 6
1.5
2
2.5
3
3.
Volume
x 10^5
NSGA-II
MoPSO
MoCSS

**(f)**

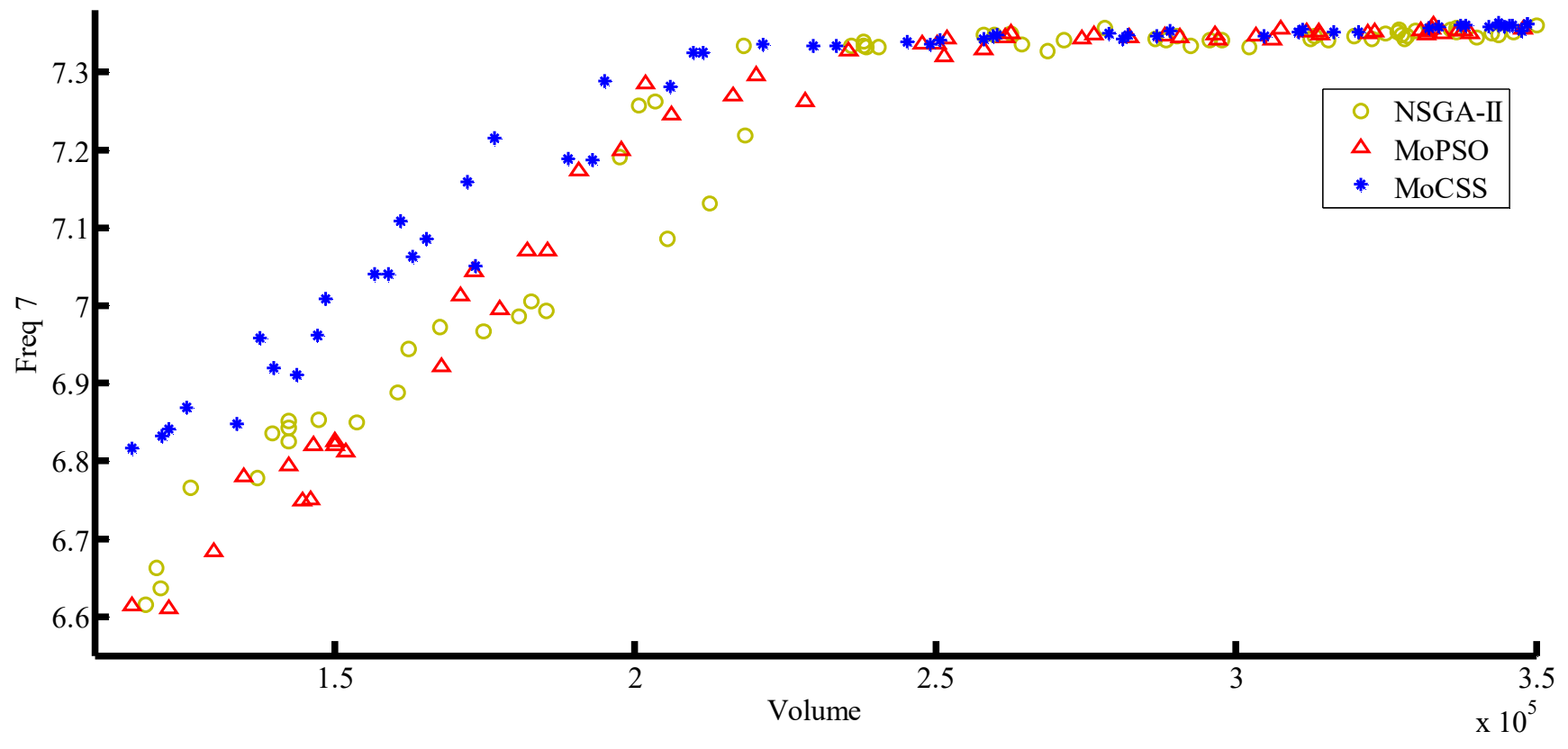

7.3
7.2
7.1
7
6.9
6.8
6.7
6.6
Freq 7
1.5
2
2.5
3
3.5
Volume
x 10^5
NSGA-II
MoPSO
MoCSS

**(g)**

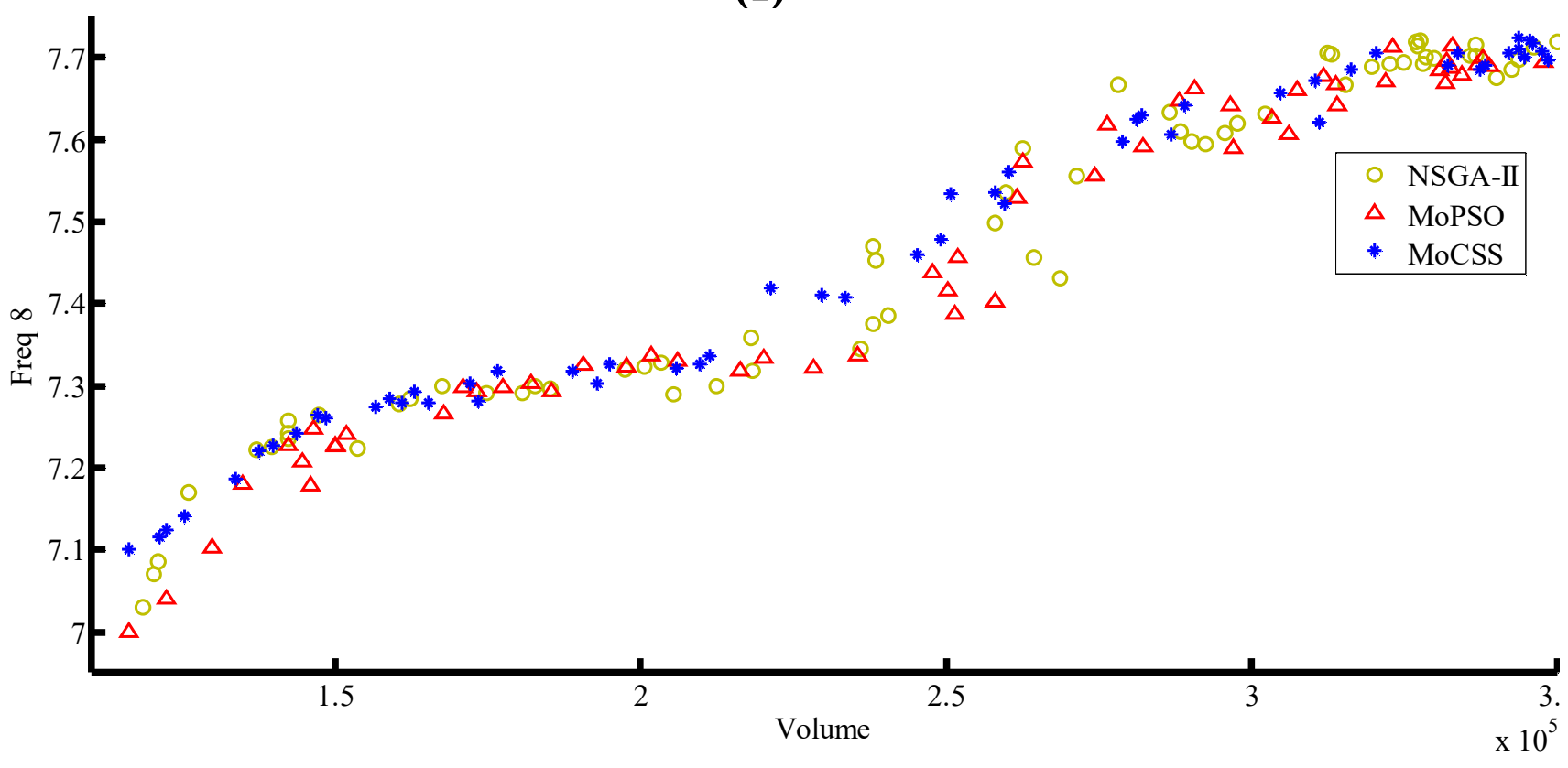

7.7
7.6
7.5
7.4
7.3
7.2
7.1
7
Freq 8
1.5
2
2.5
3
3.
Volume
x 10^5
NSGA-II
MoPSO
MoCSS

**(h)**

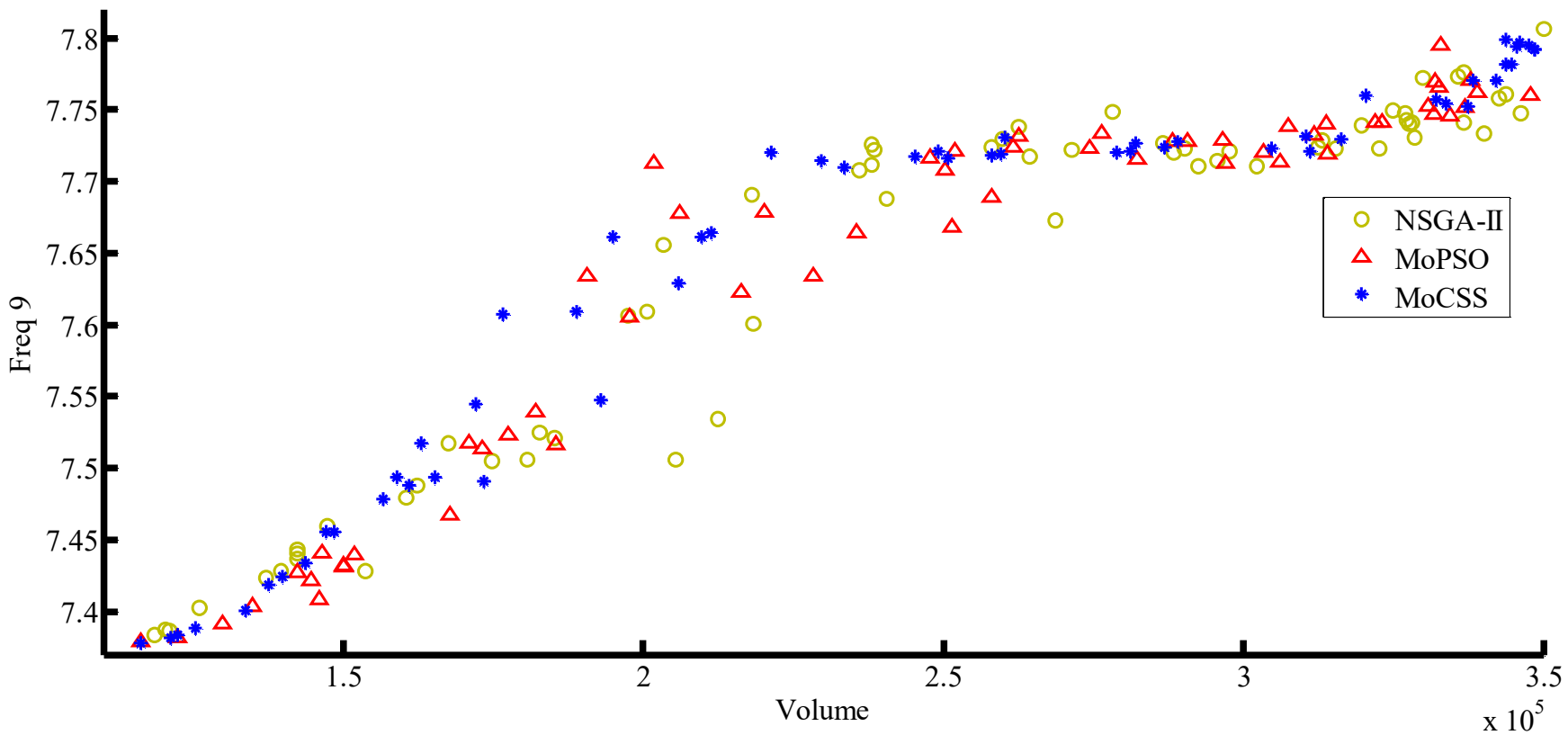

**(i)**

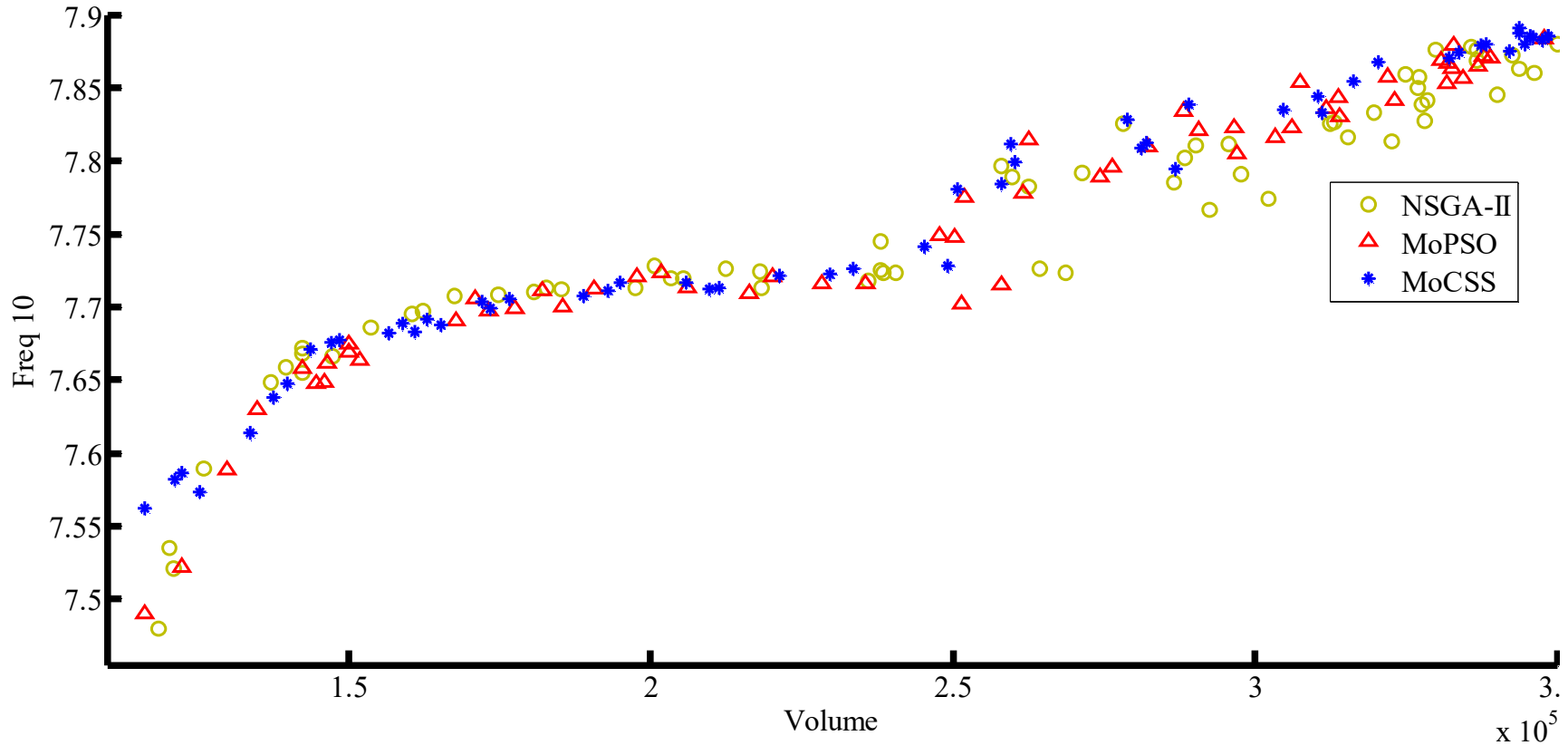

**(j)**

Fig. 8. The related Volume-Frequencies of three methods for P1; **(a)** Vol-Fr1, **(b)** Vol-Fr2**, (c)** Vol-Fr3, **(d)** Vol-Fr4**, (e)** Vol-Fr5**, (f)** Vol-Fr6**, (g)** Vol-Fr7**, (h)** Vol-Fr8**, (i)** Vol-Fr9**, (j)** Vol-Fr10**.**

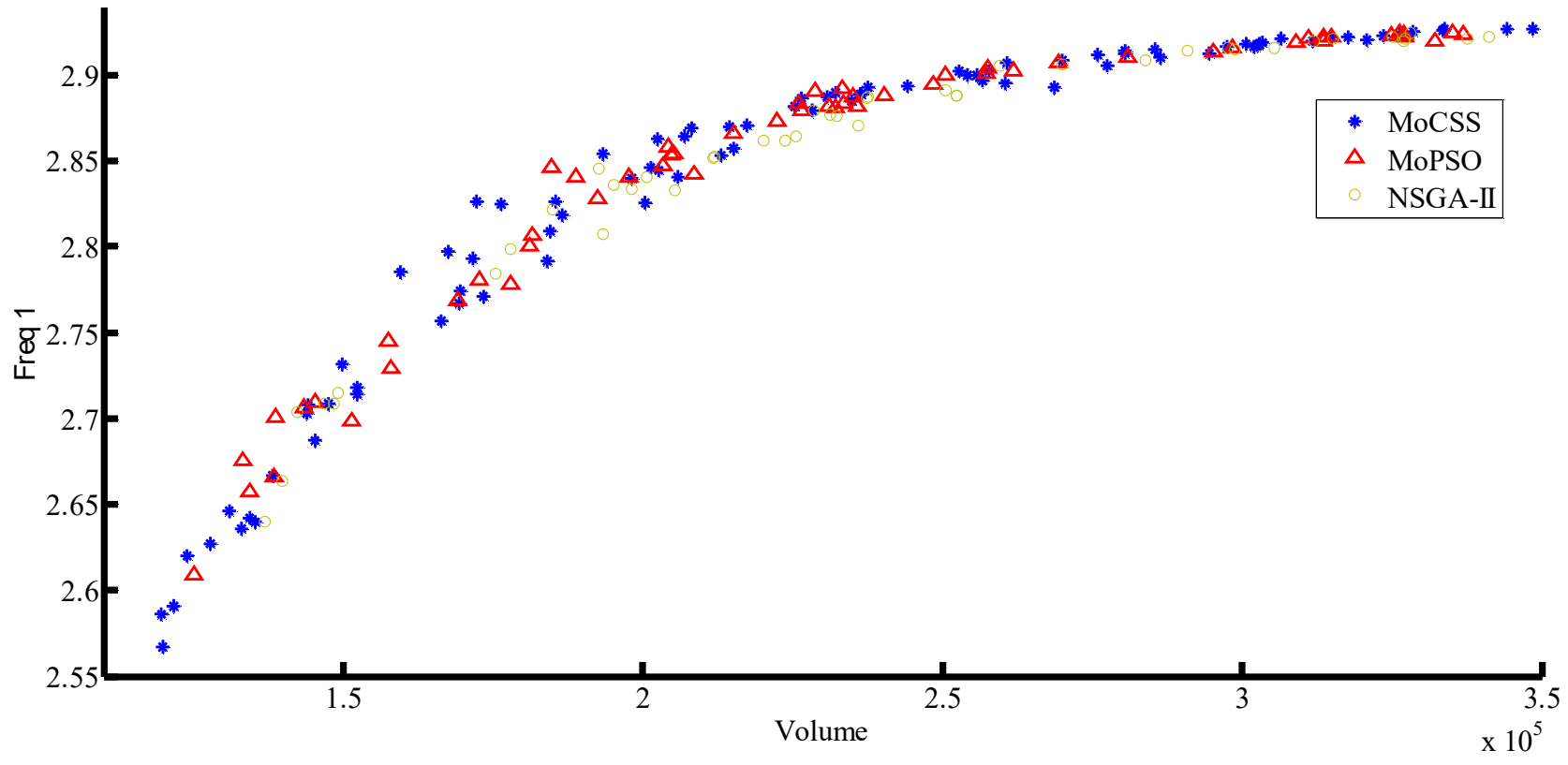
Freq 1
2.55
2.6
2.65
2.7
2.75
2.8
2.85
2.9
1.5
2
2.5
3
3.5
Volume
x 10^5
MoCSS
MoPSO
NSGA-II

**(a)**

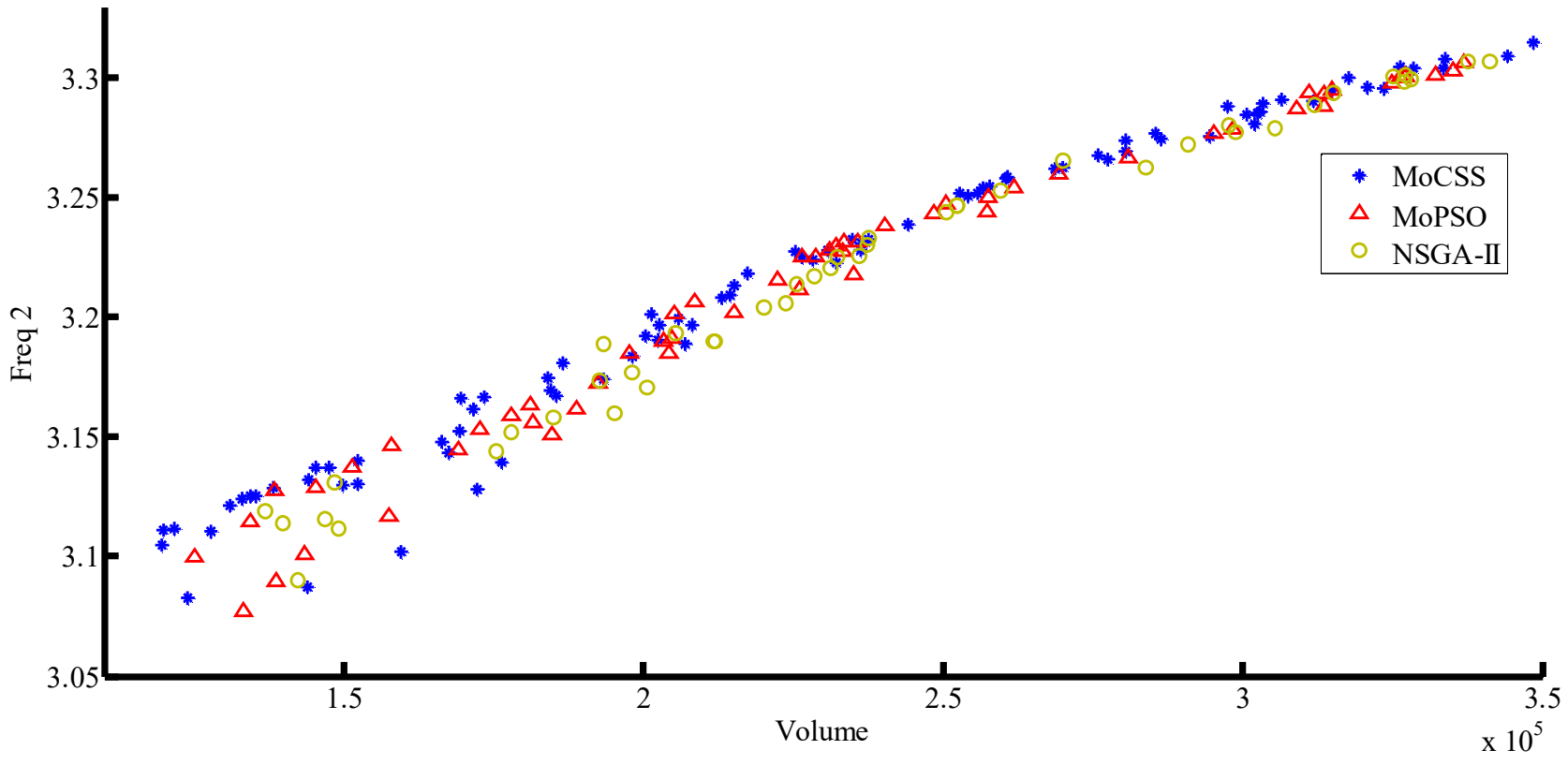
Freq 2
3.05
3.1
3.15
3.2
3.25
3.3
1.5
2
2.5
3
3.5
Volume
x 10^5
MoCSS
MoPSO
NSGA-II

**(b)**

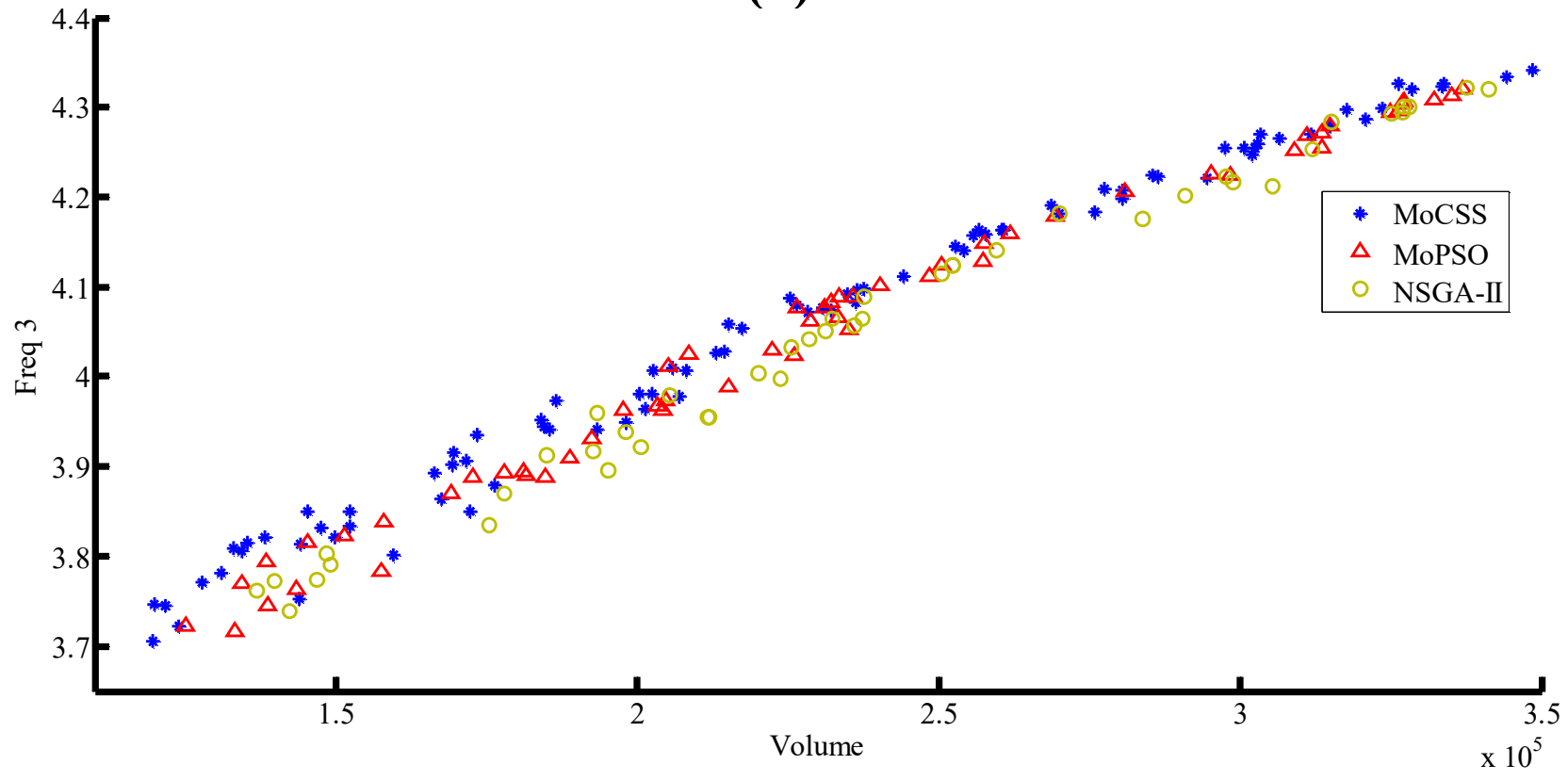
Freq 3
3.7
3.8
3.9
4
4.1
4.2
4.3
4.4
1.5
2
2.5
3
3.5
Volume
x 10^5
MoCSS
MoPSO
NSGA-II

**(c)**

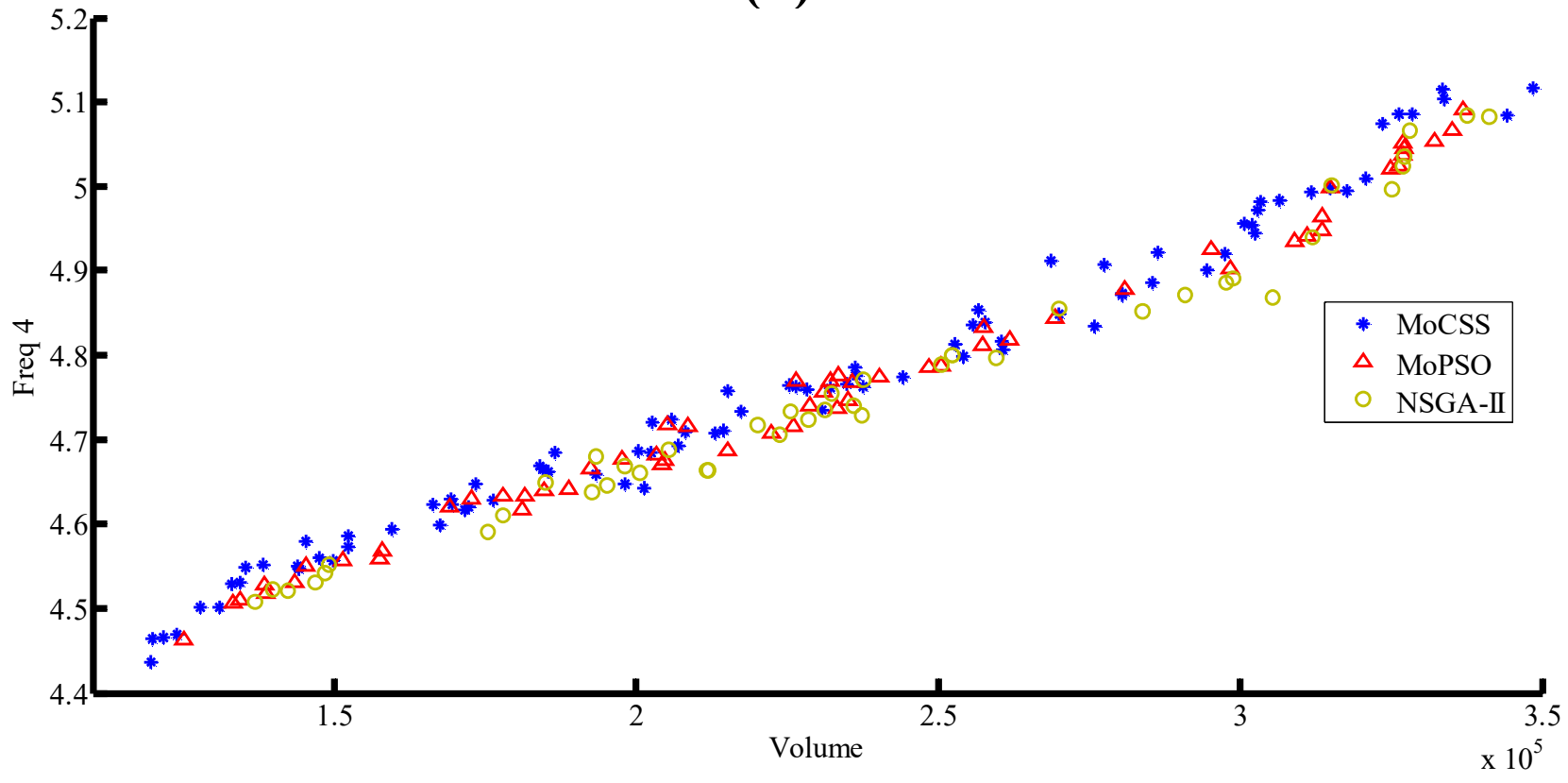
Freq 4
4.4
4.5
4.6
4.7
4.8
4.9
5
5.1
5.2
1.5
2
2.5
3
3.5
Volume
x 10^5
MoCSS
MoPSO
NSGA-II

**(d)**

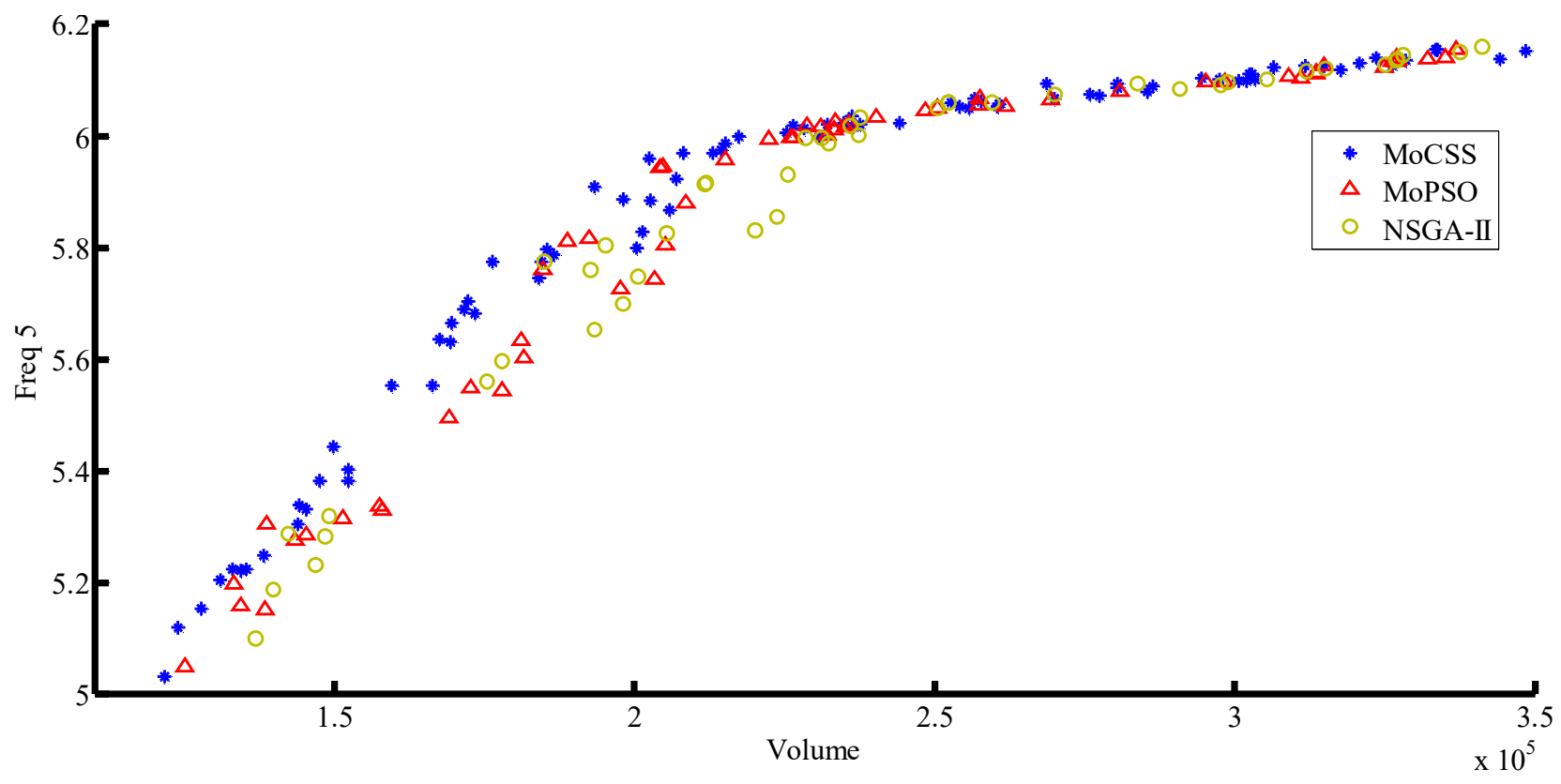

Freq 5
Volume
x 10^5
MoCSS
MoPSO
NSGA-II

**(e)**

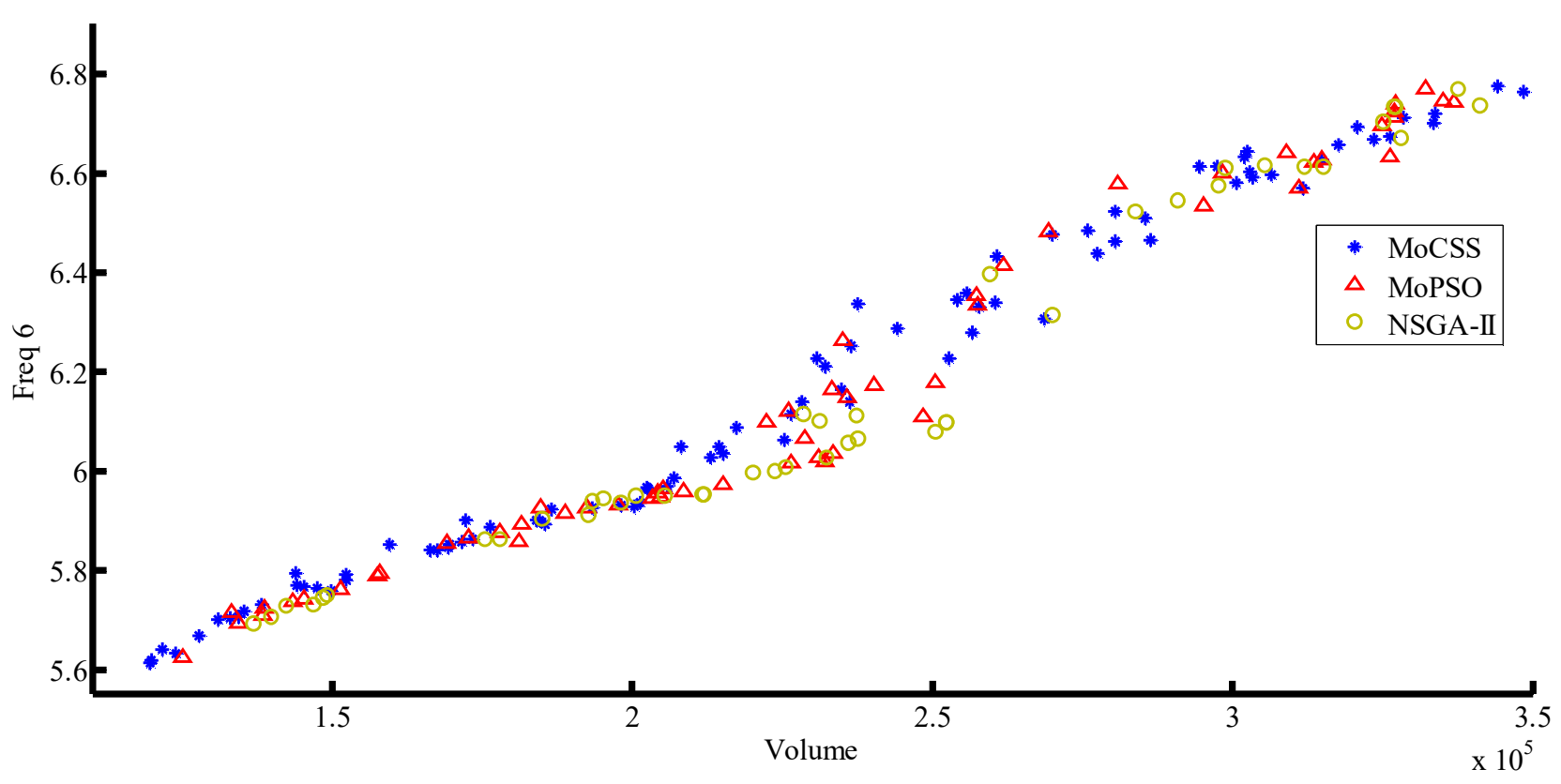

Freq 6
Volume
x 10^5
MoCSS
MoPSO
NSGA-II

**(f)**

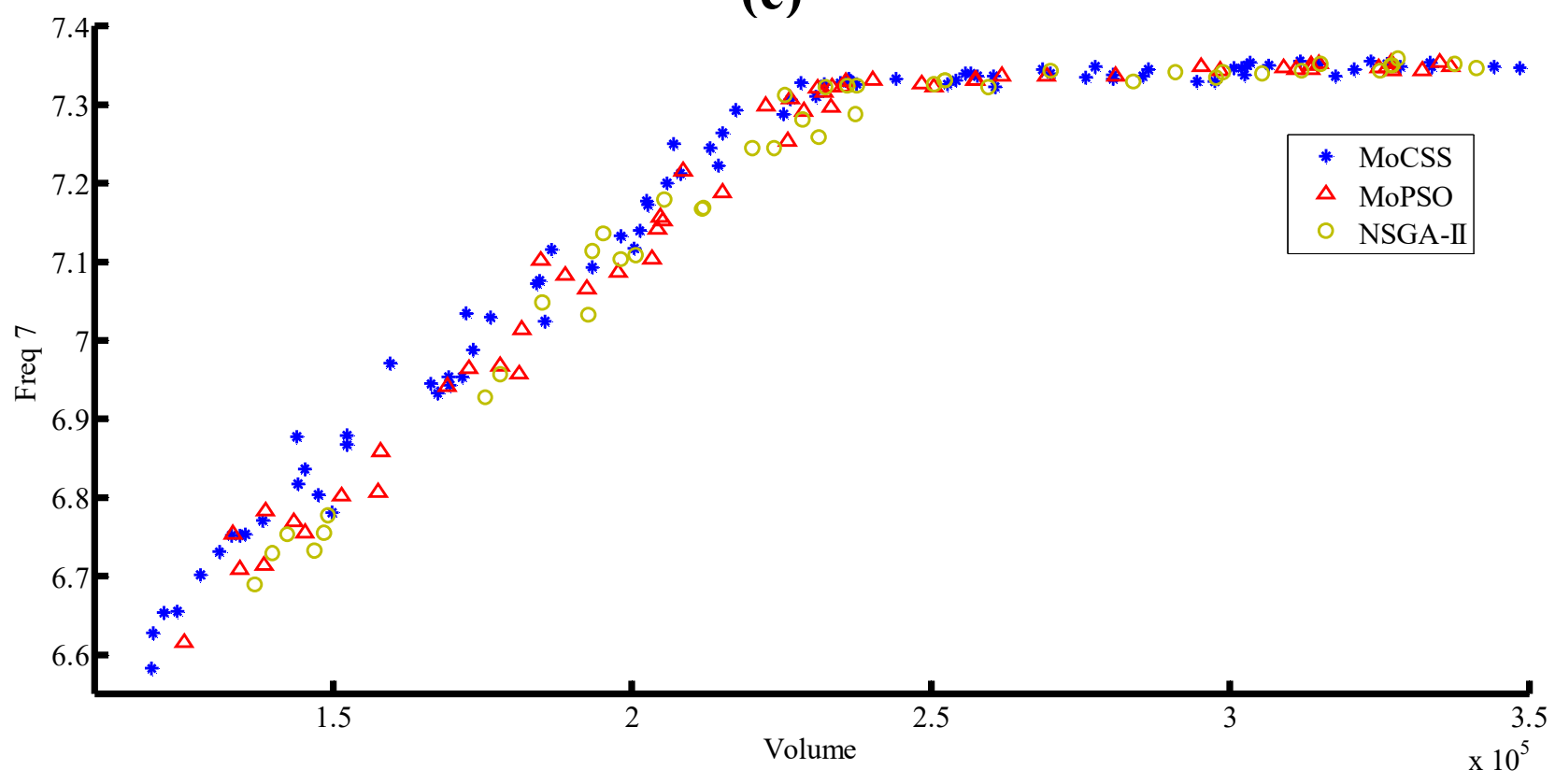

Freq 7
Volume
x 10^5
MoCSS
MoPSO
NSGA-II

**(g)**

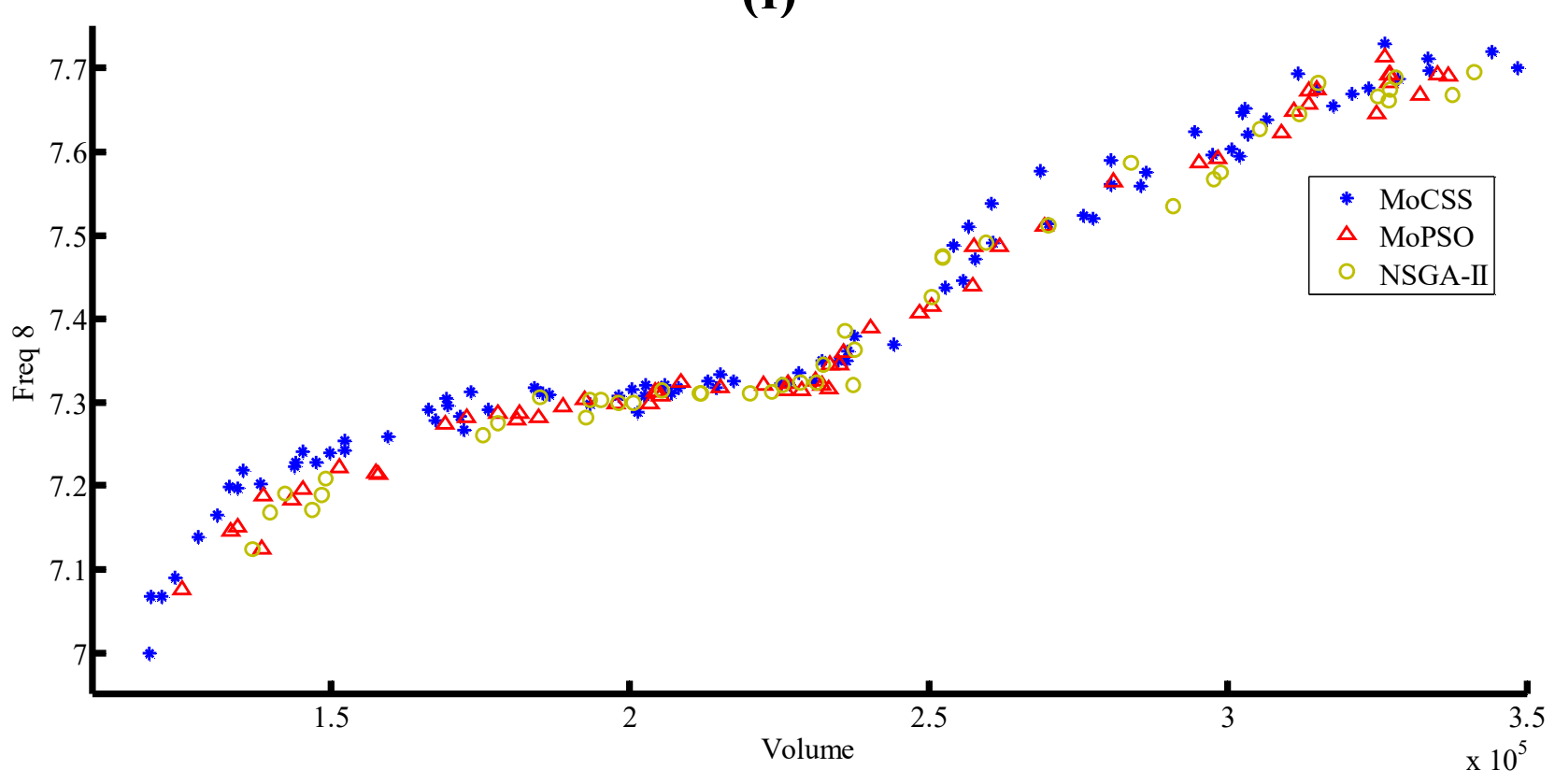

Freq 8
Volume
x 10^5
MoCSS
MoPSO
NSGA-II

**(h)**

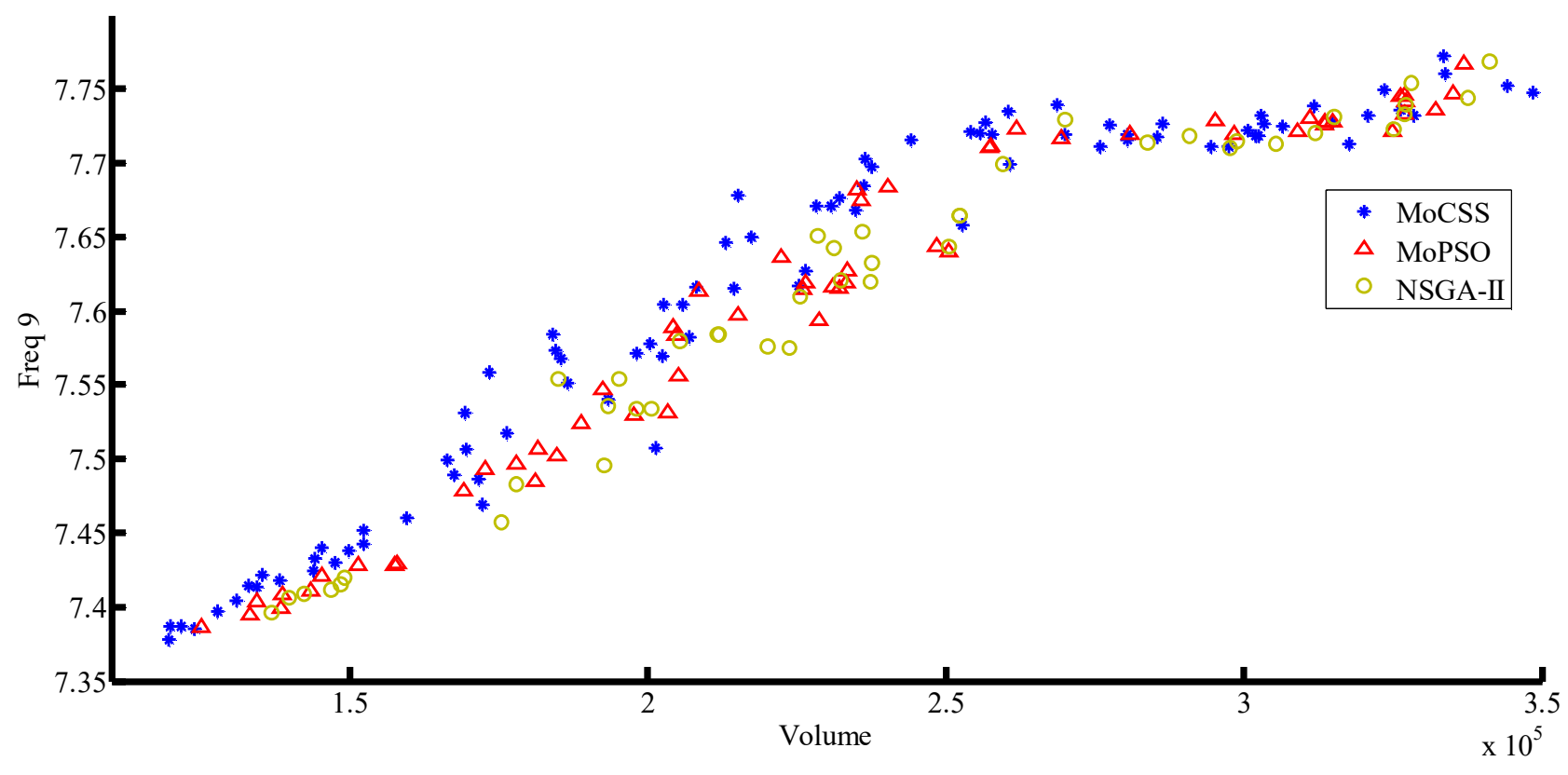

**(i)**

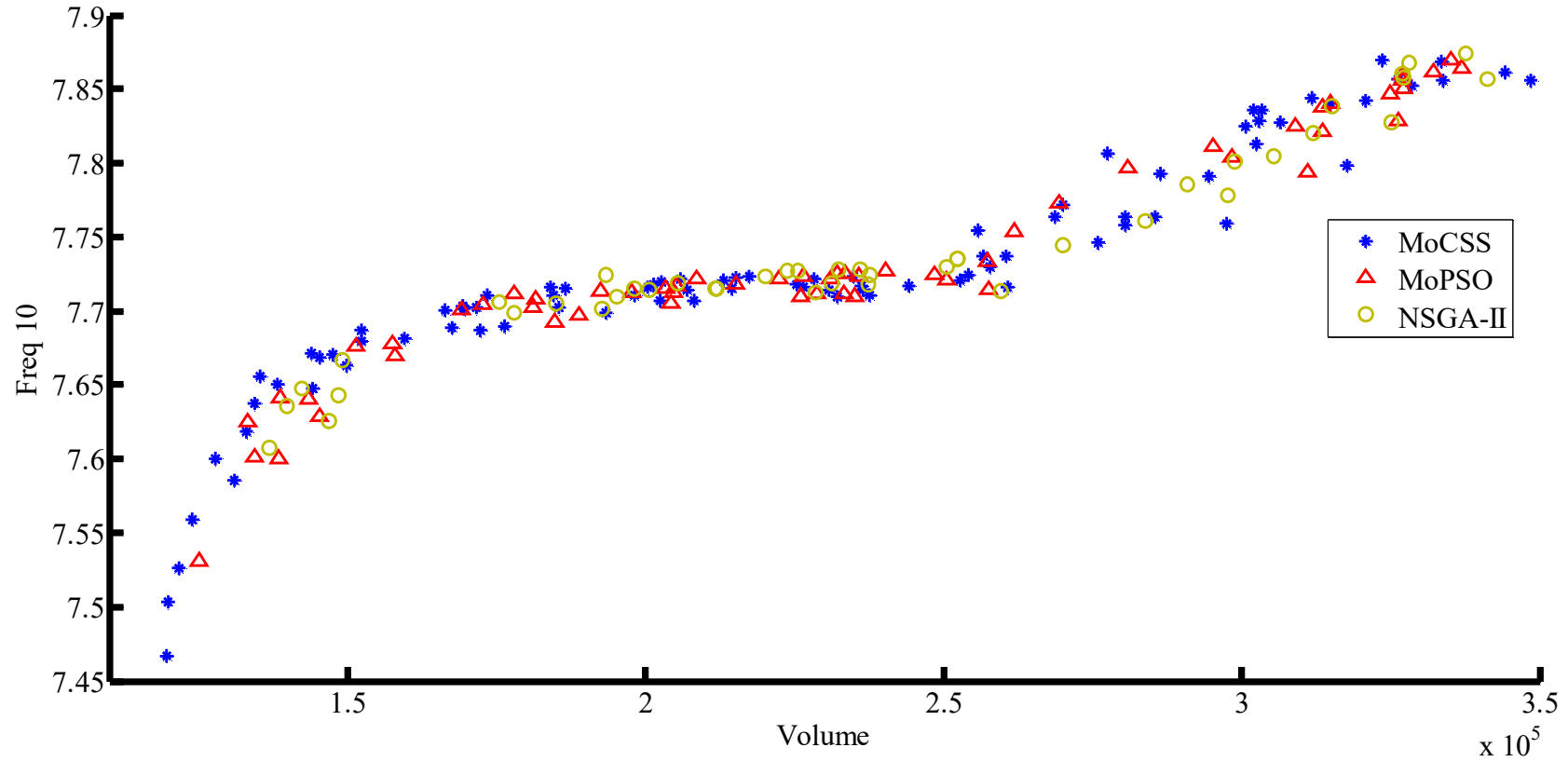

**(j)**

Fig. 9. The related Volume-Frequencies of three methods for P2; **(a)** Vol-Fr1, **(b)** Vol-Fr2**, (c)** Vol-Fr3, **(d)** Vol-Fr4**, (e)** Vol-Fr5**, (f)** Vol-Fr6**, (g)** Vol-Fr7**, (h)** Vol-Fr8**, (i)** Vol-Fr9**, (j)** Vol-Fr10**.**

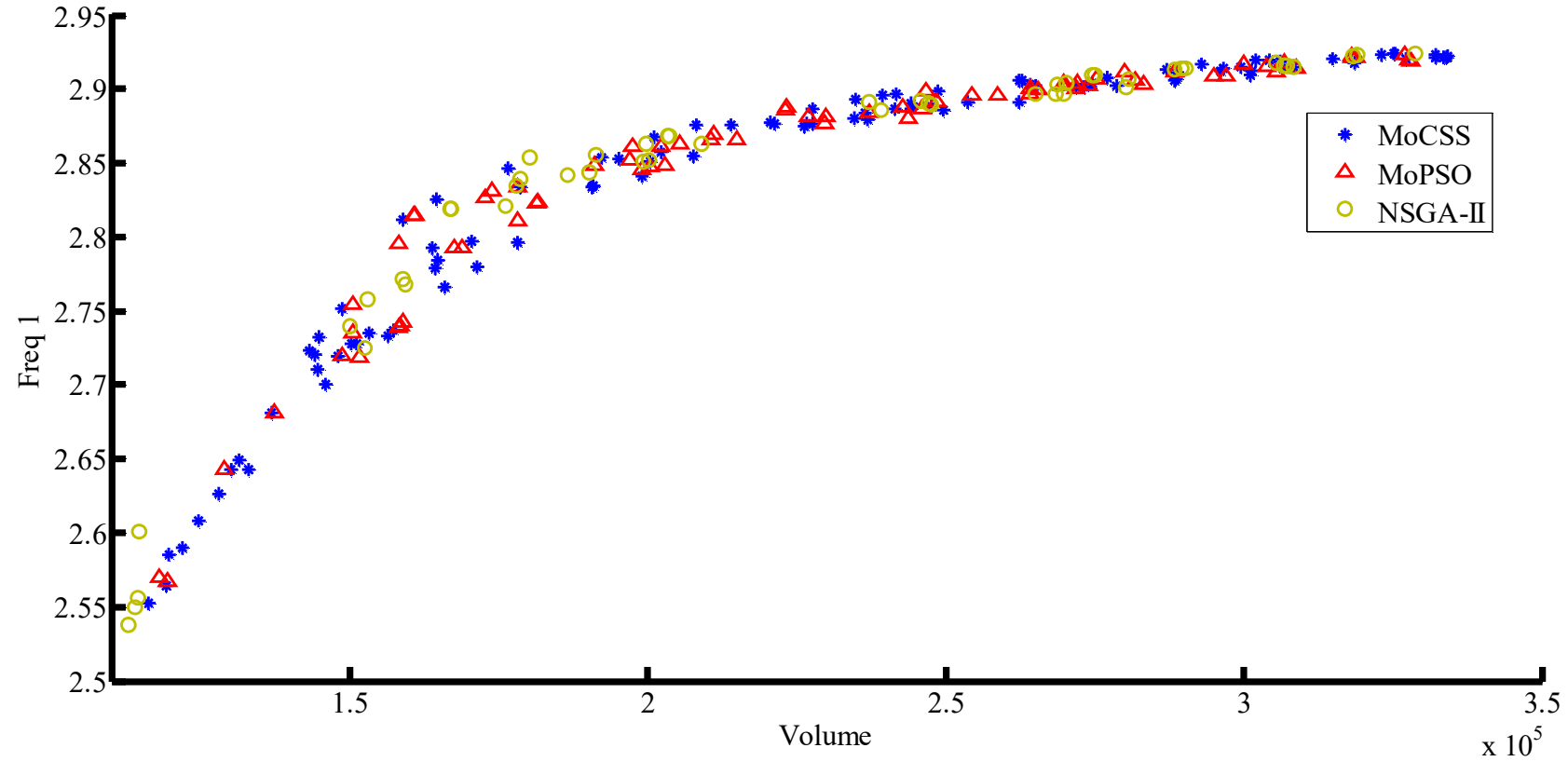

**(a)**

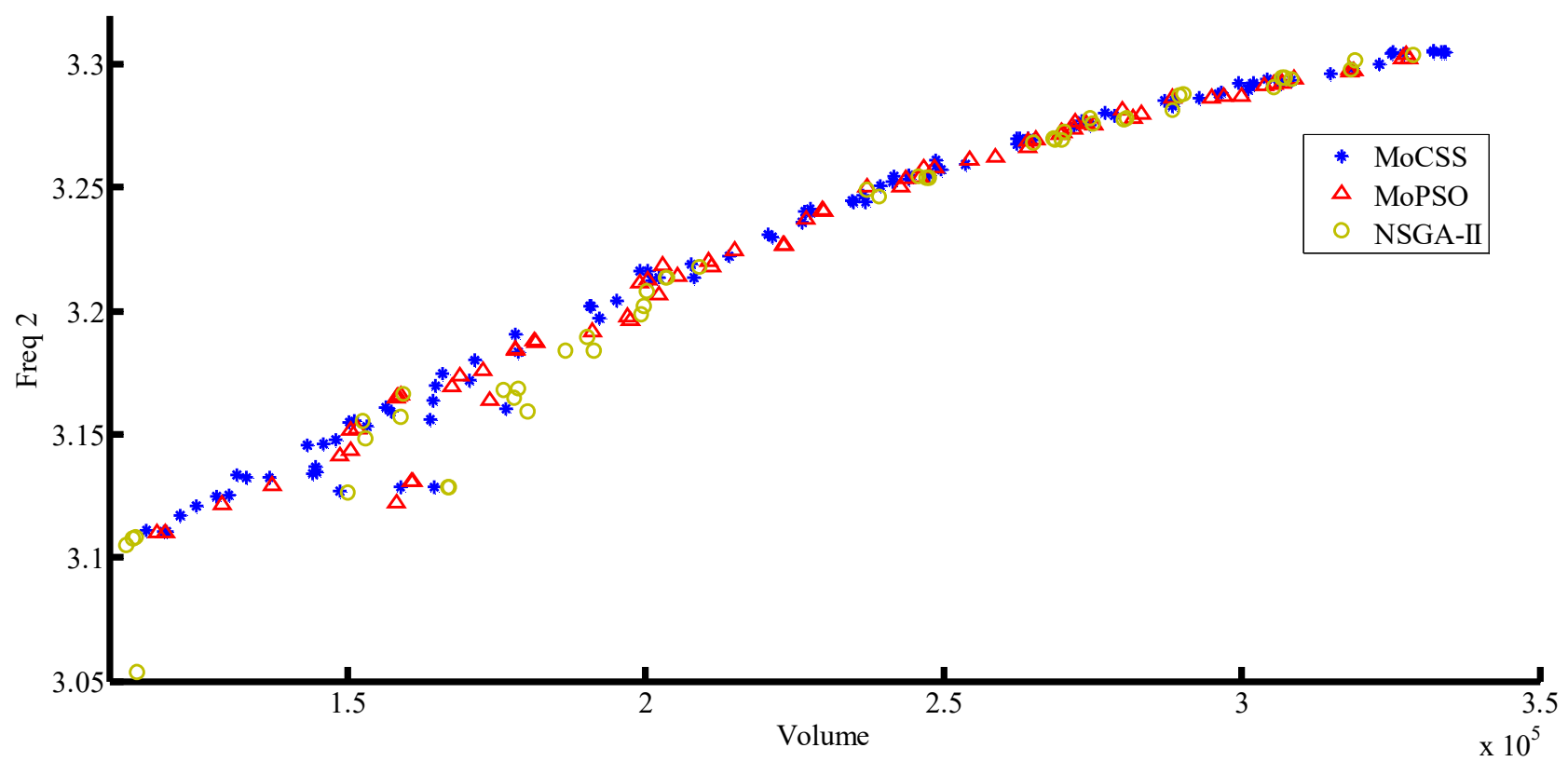

**(b)**

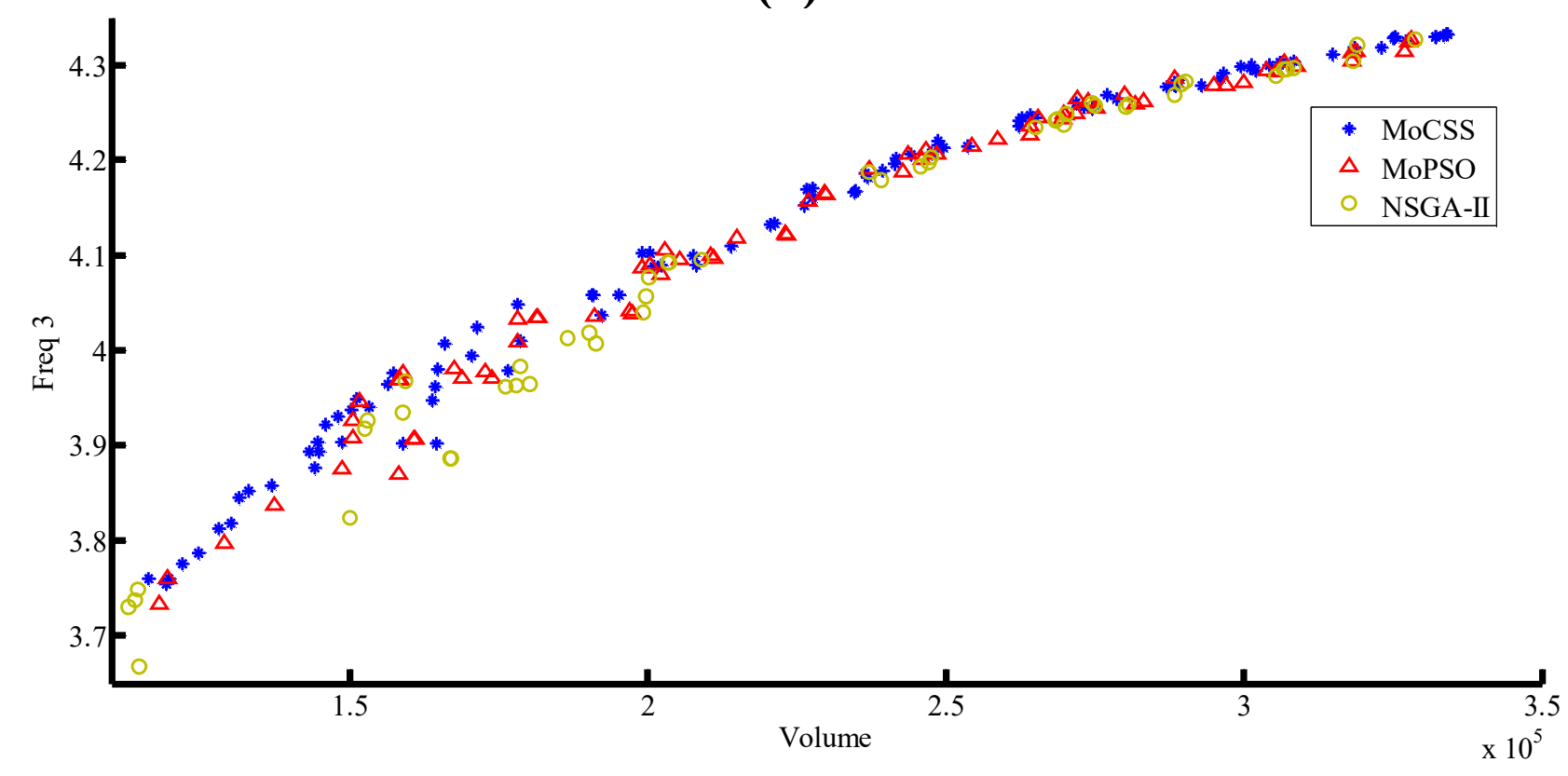

**(c)**

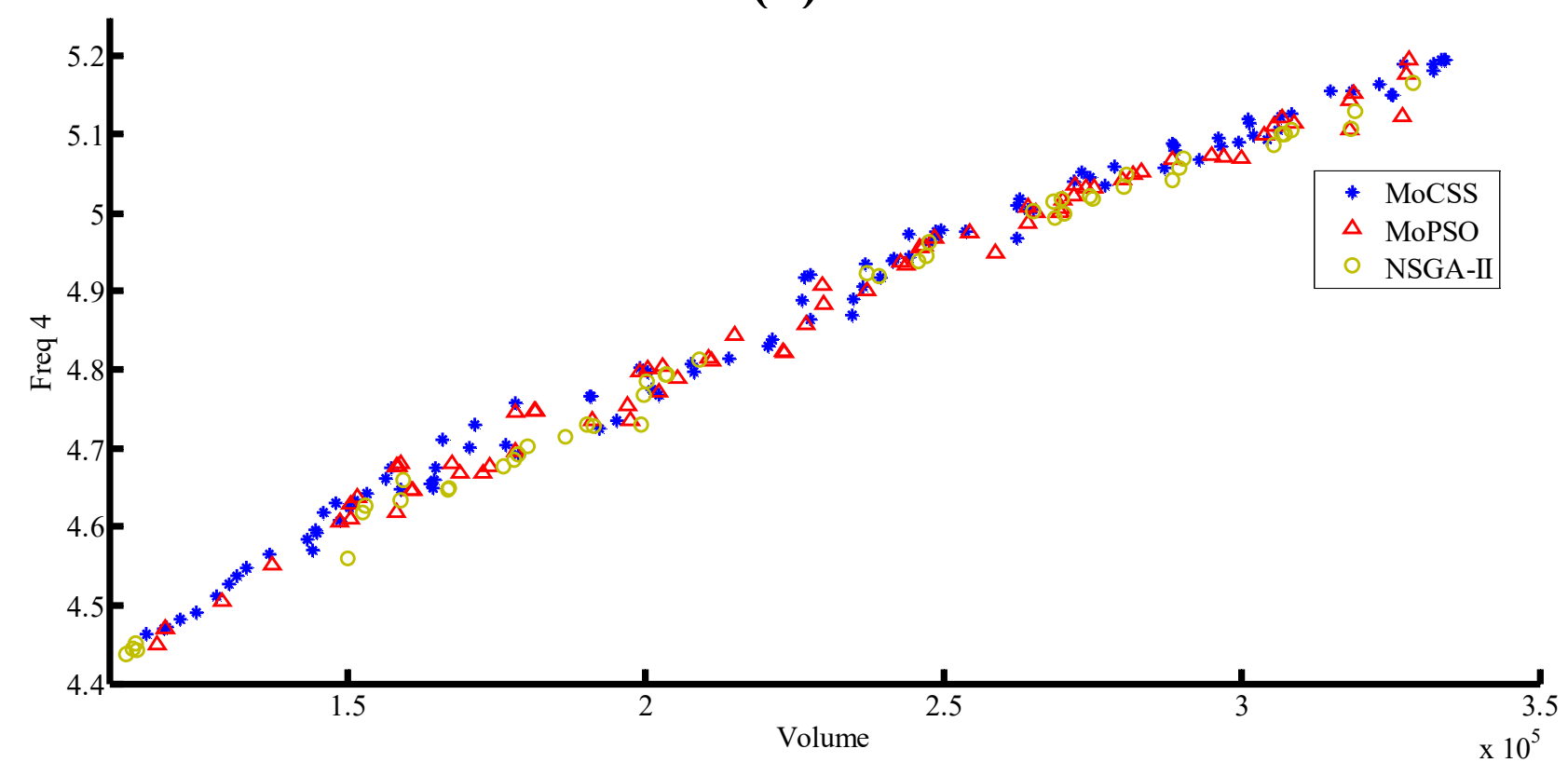

**(d)**

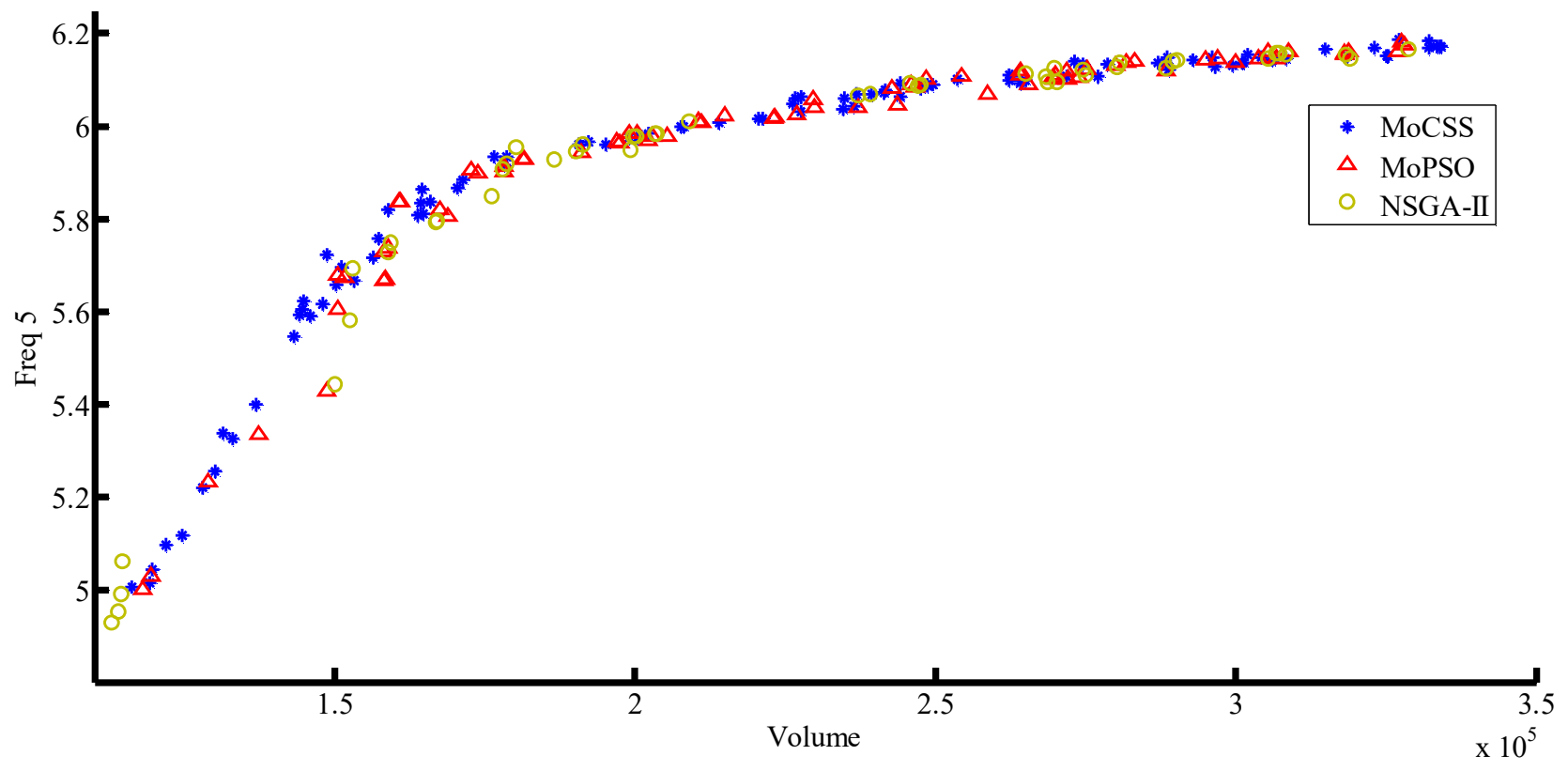

**(e)**

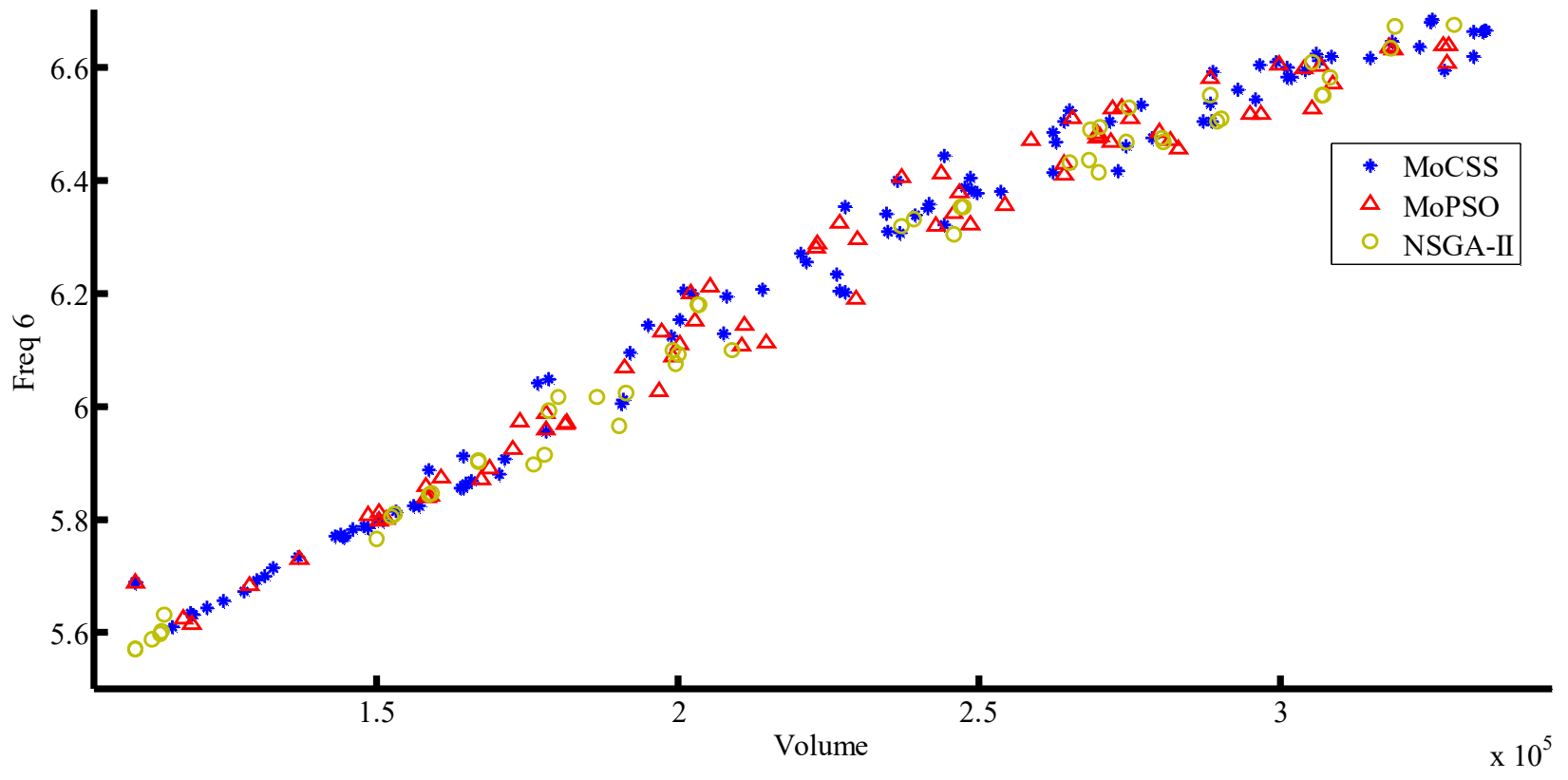

**(f)**

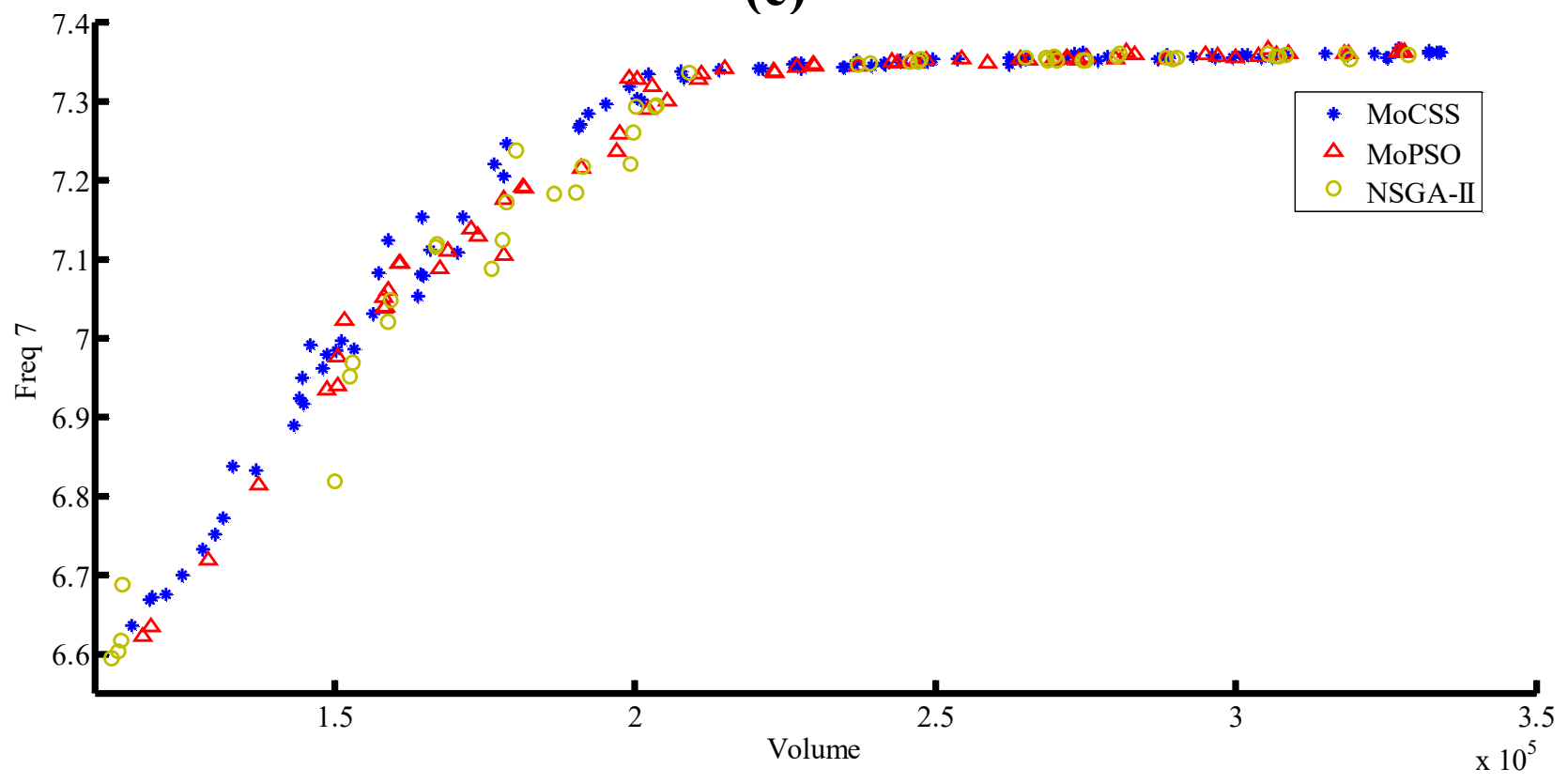

**(g)**

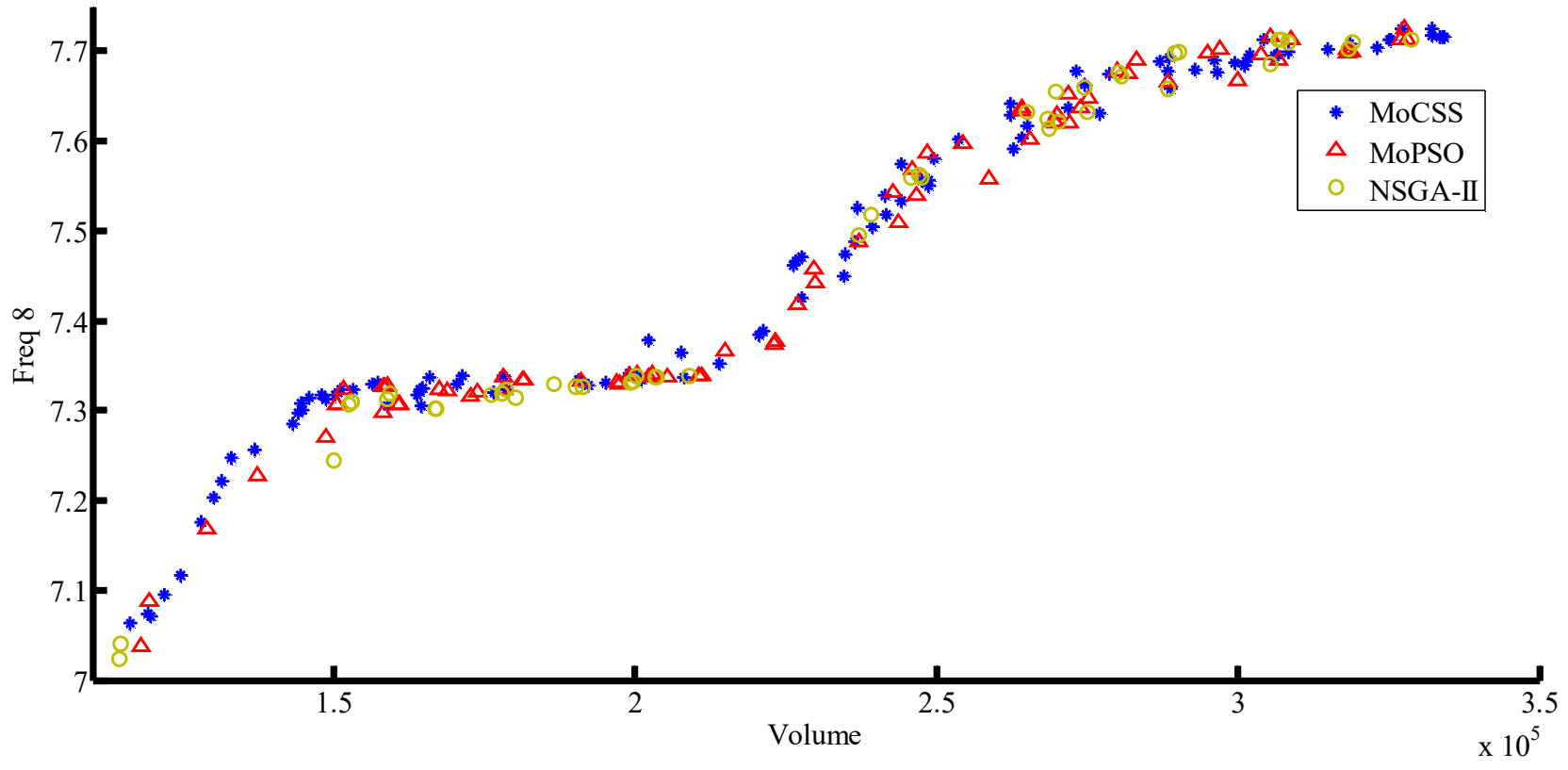

**(h)**

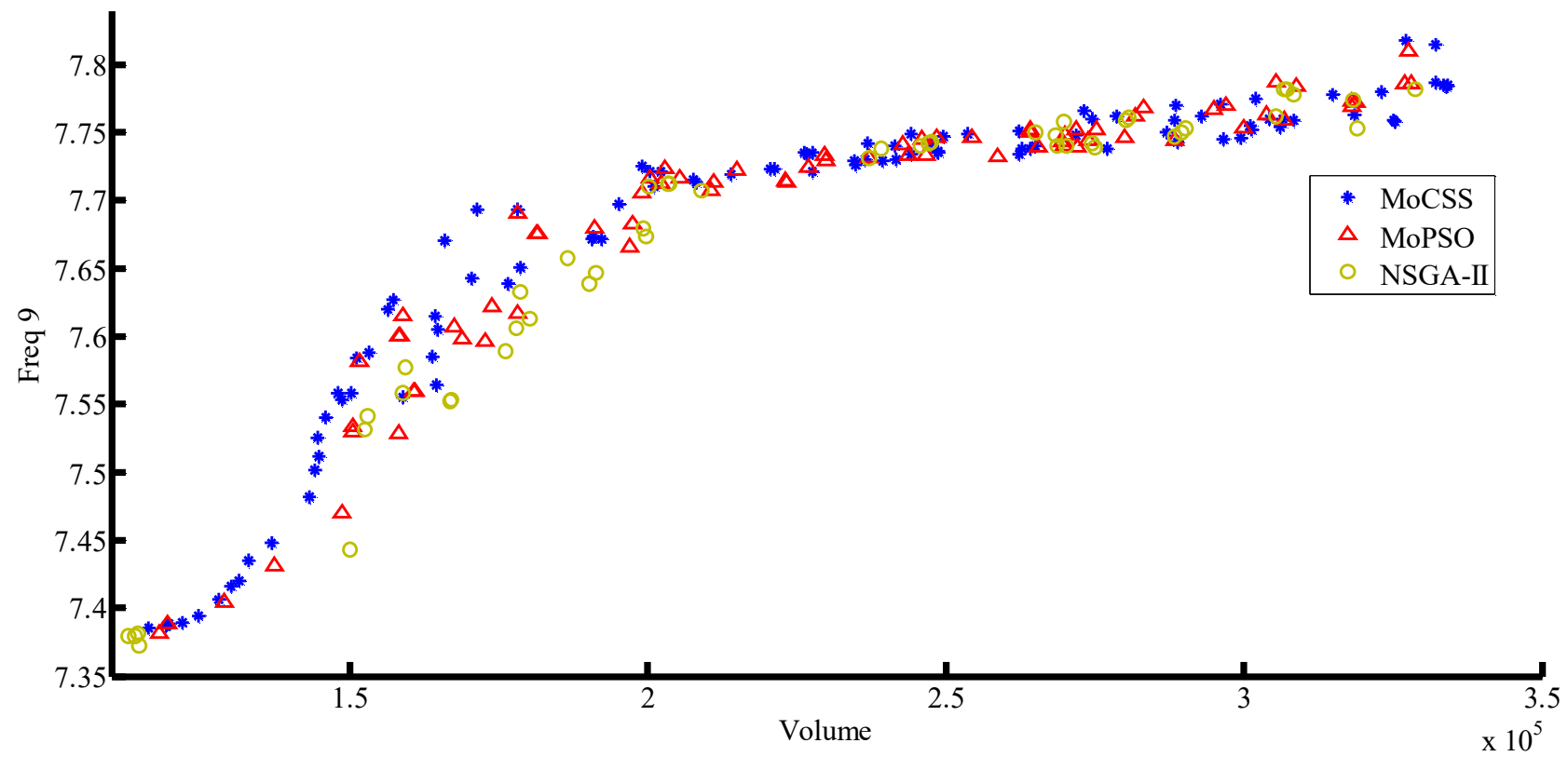

**(i)**

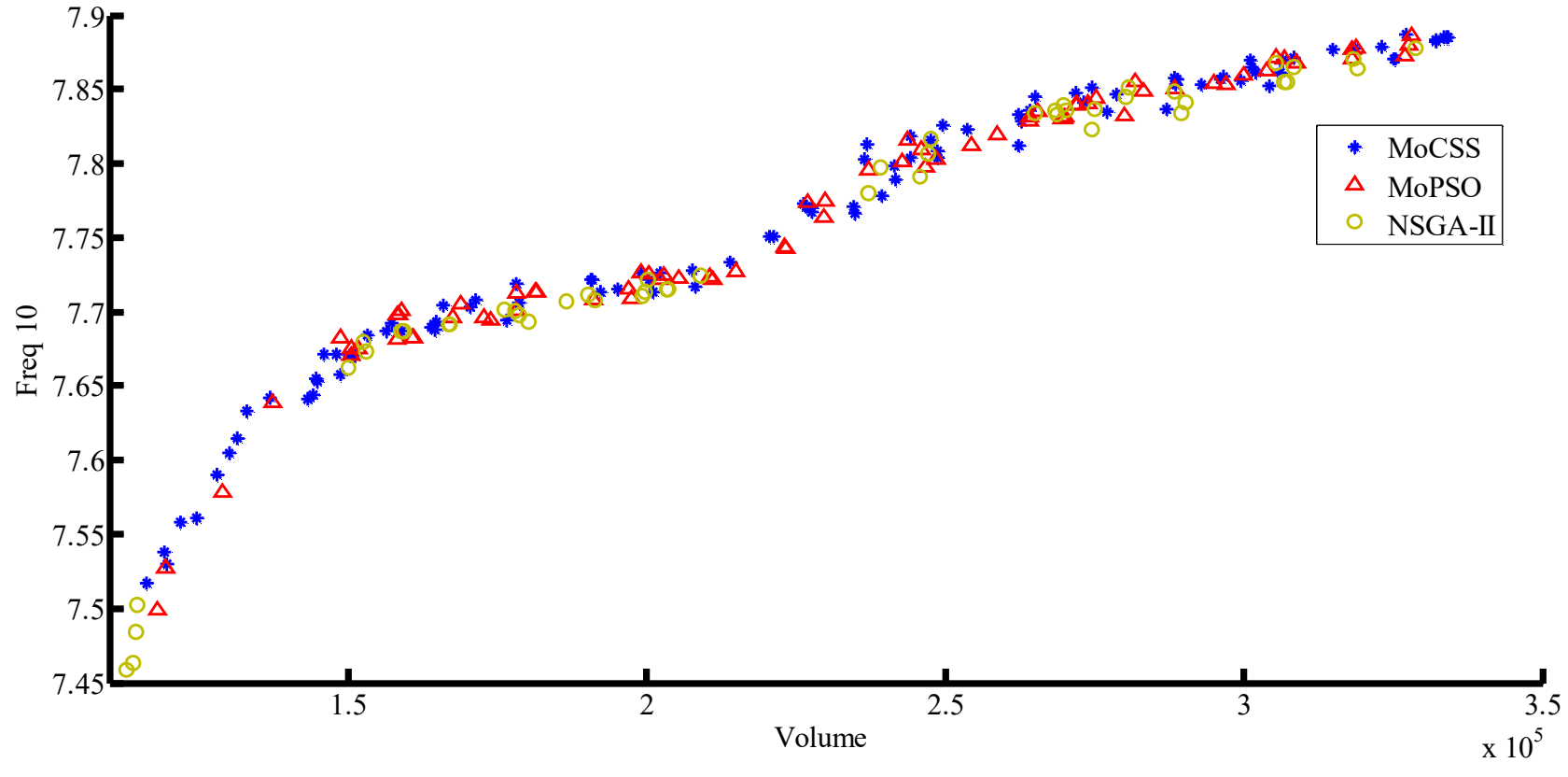

**(j)**

Fig. 10. The related Volume-Frequencies of three methods for P3; **(a)** Vol-Fr1, **(b)** Vol-Fr2**, (c)** Vol-Fr3, **(d)** Vol-Fr4**, (e)** Vol-Fr5**, (f)** Vol-Fr6**, (g)** Vol-Fr7**, (h)** Vol-Fr8**, (i)** Vol-Fr9**, (j)** Vol-Fr10**.**

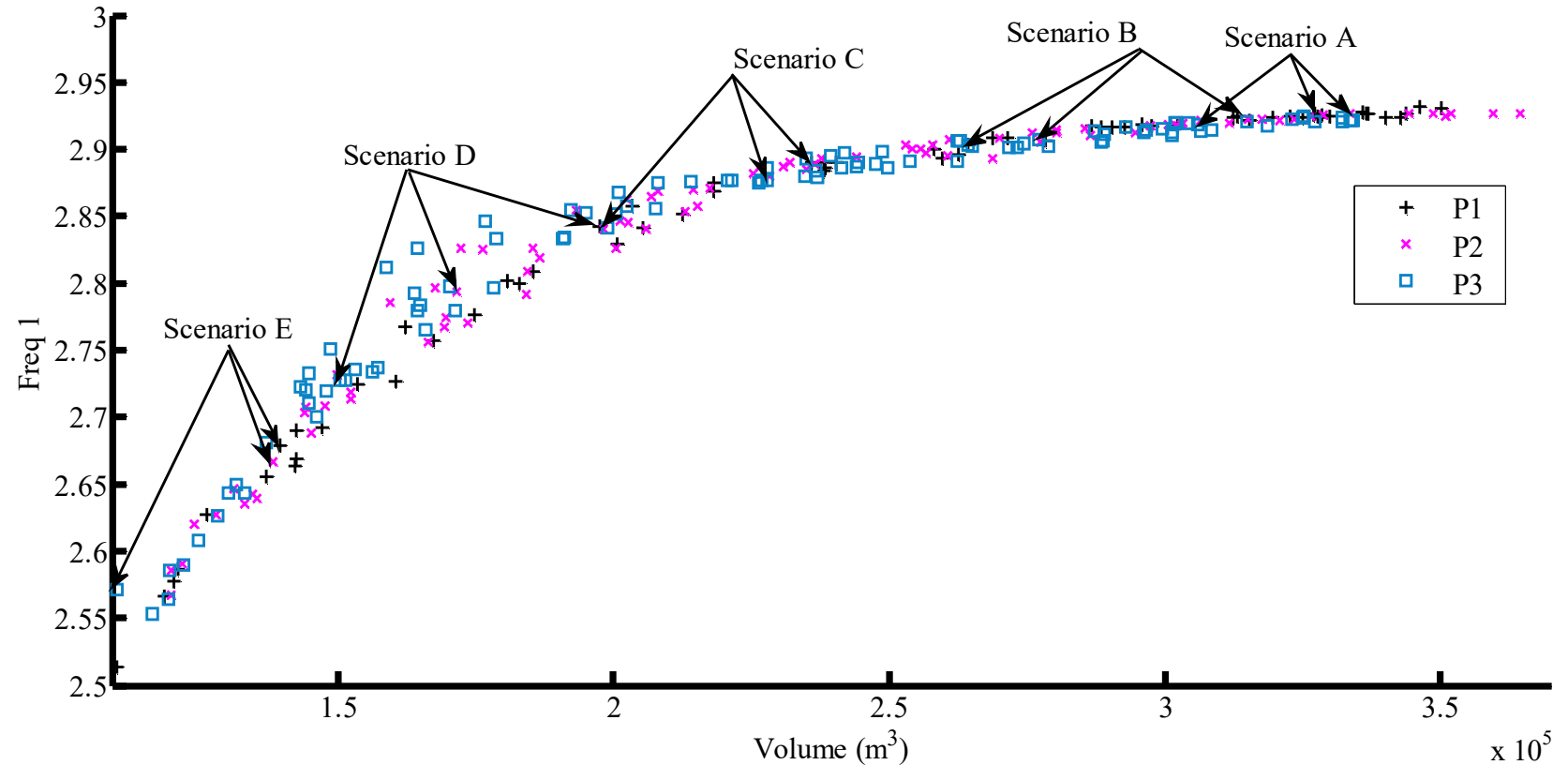

**(a)**

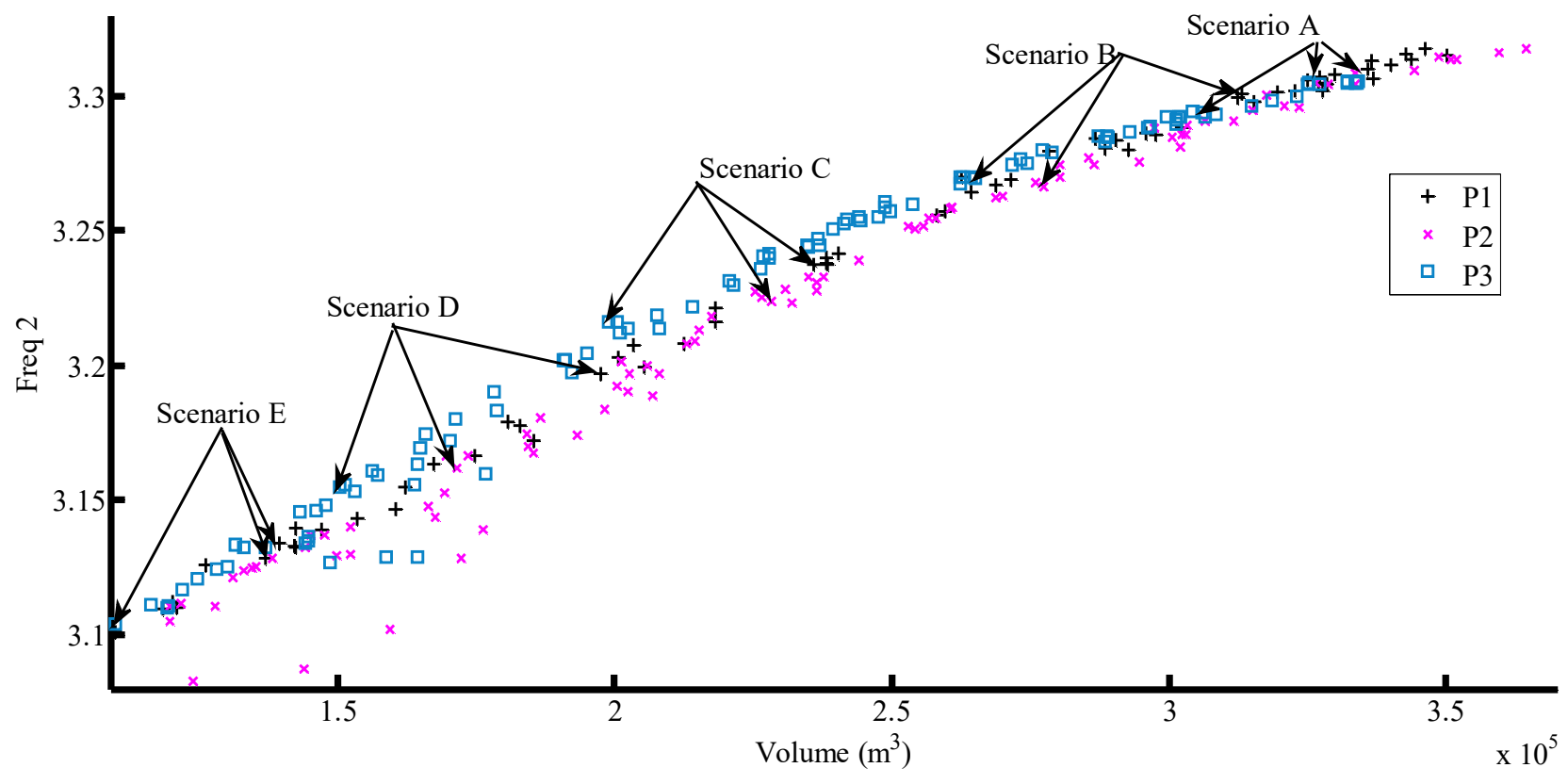

**(b)**

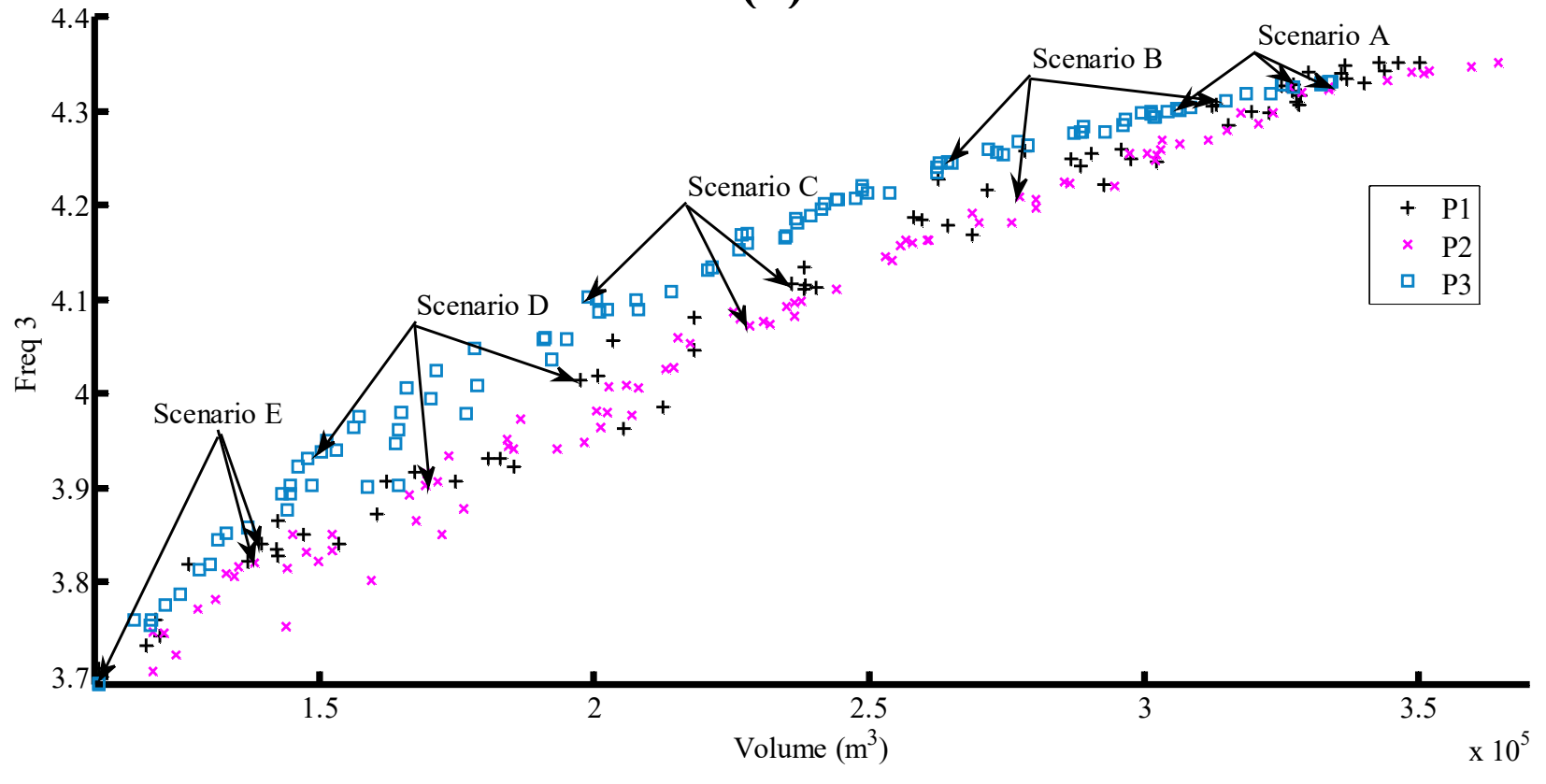

**(c)**

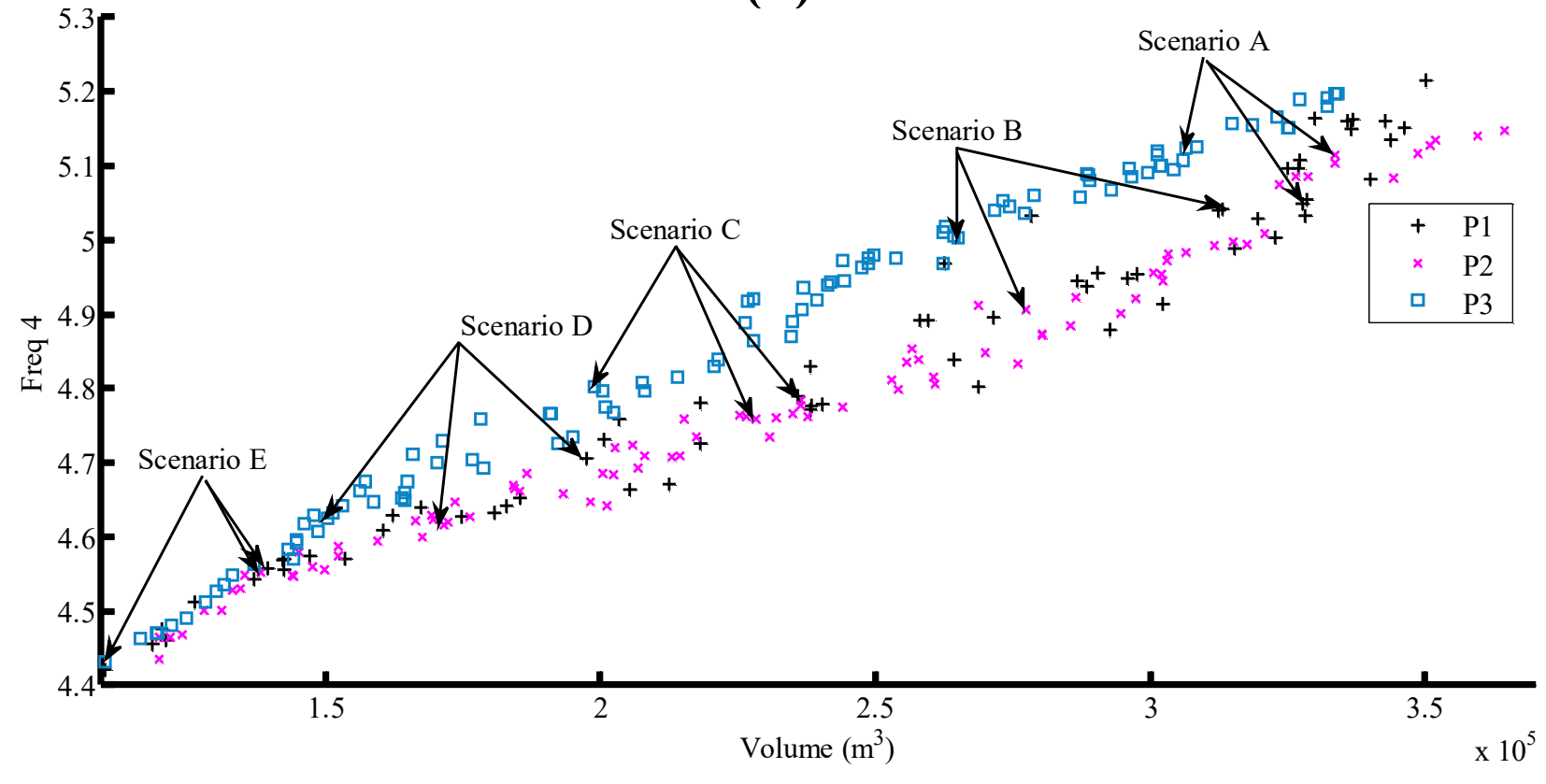

**(d)**

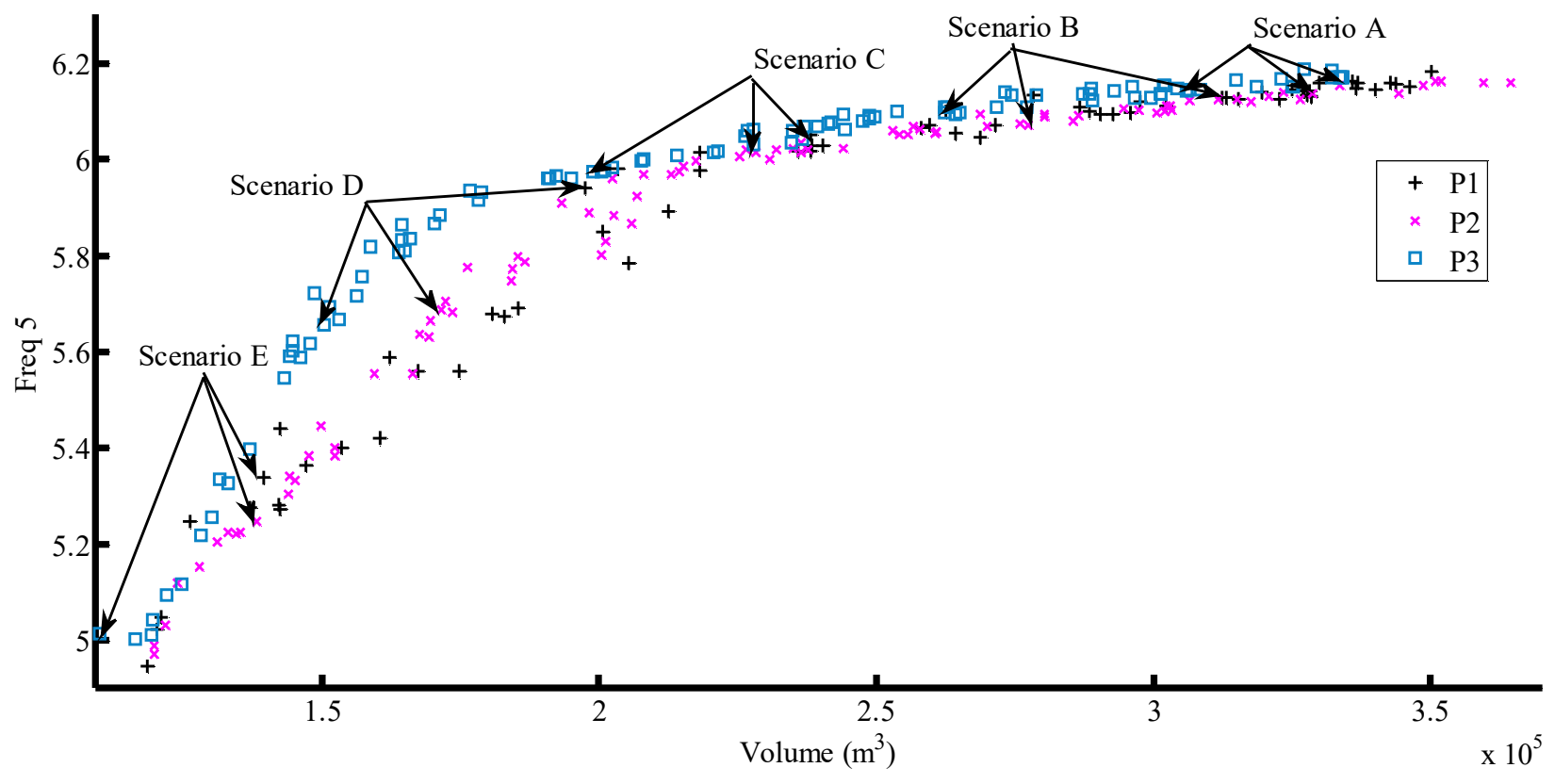

**(e)**

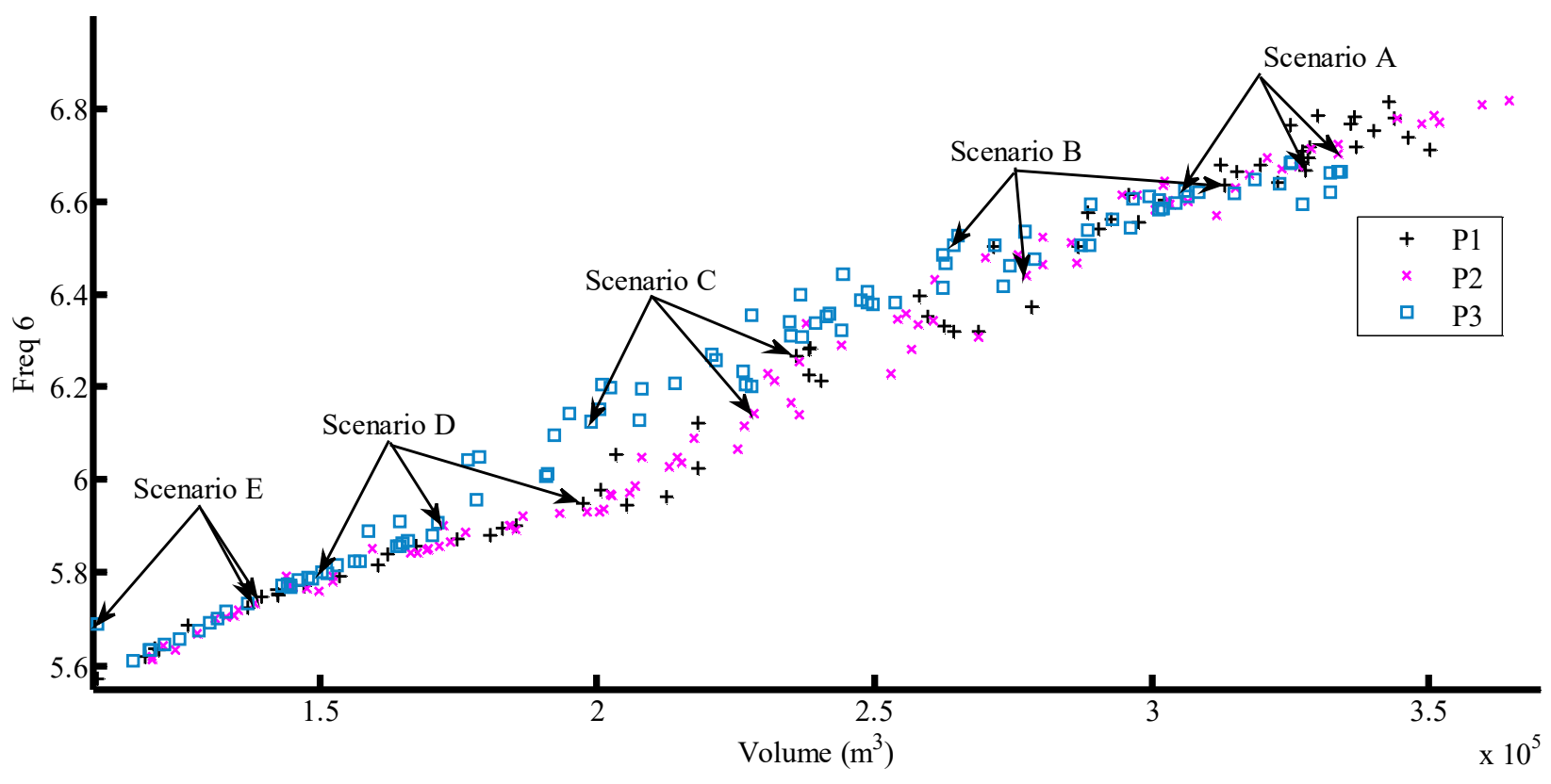

**(f)**

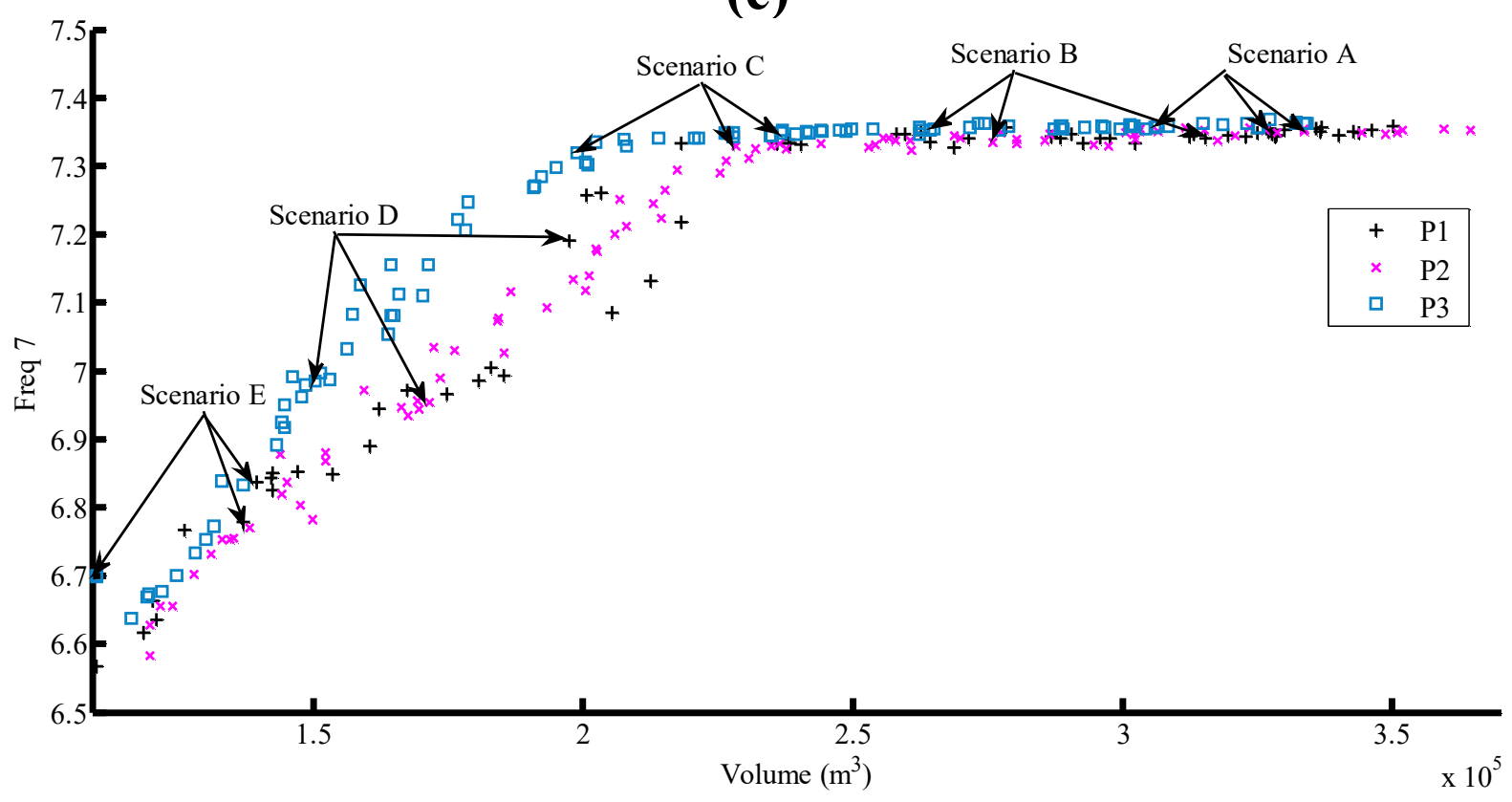

**(g)**

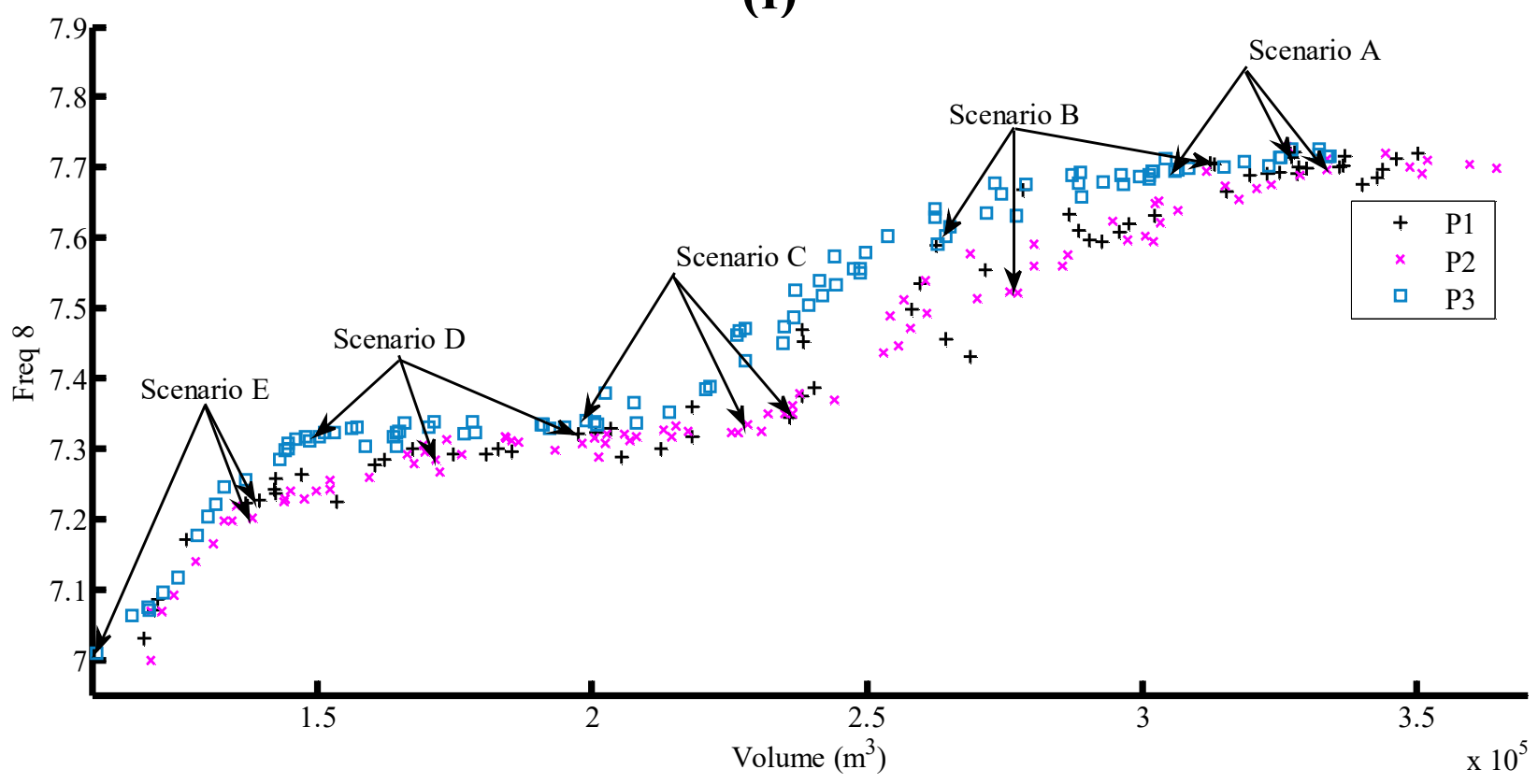

**(h)**

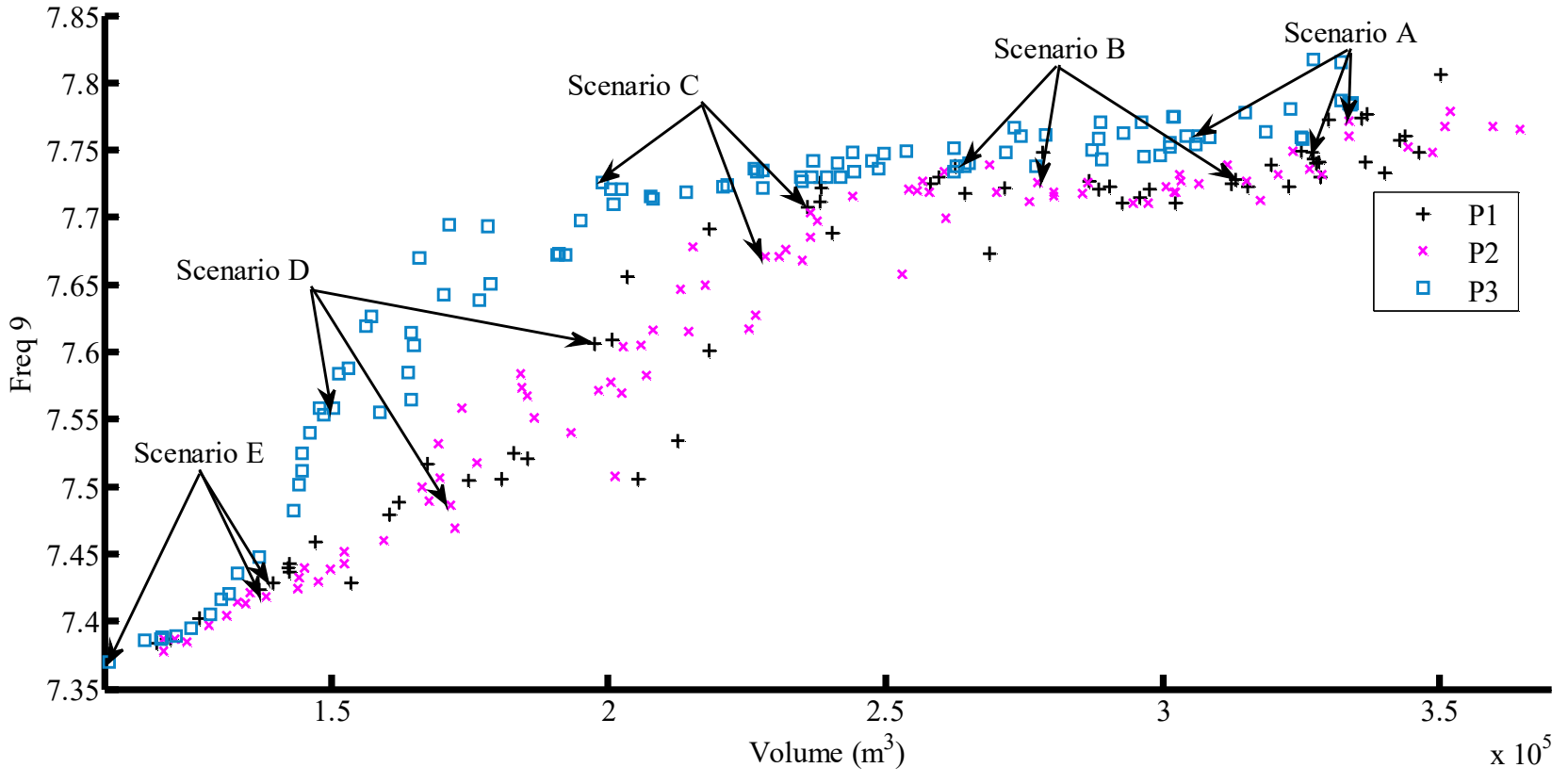

(i)

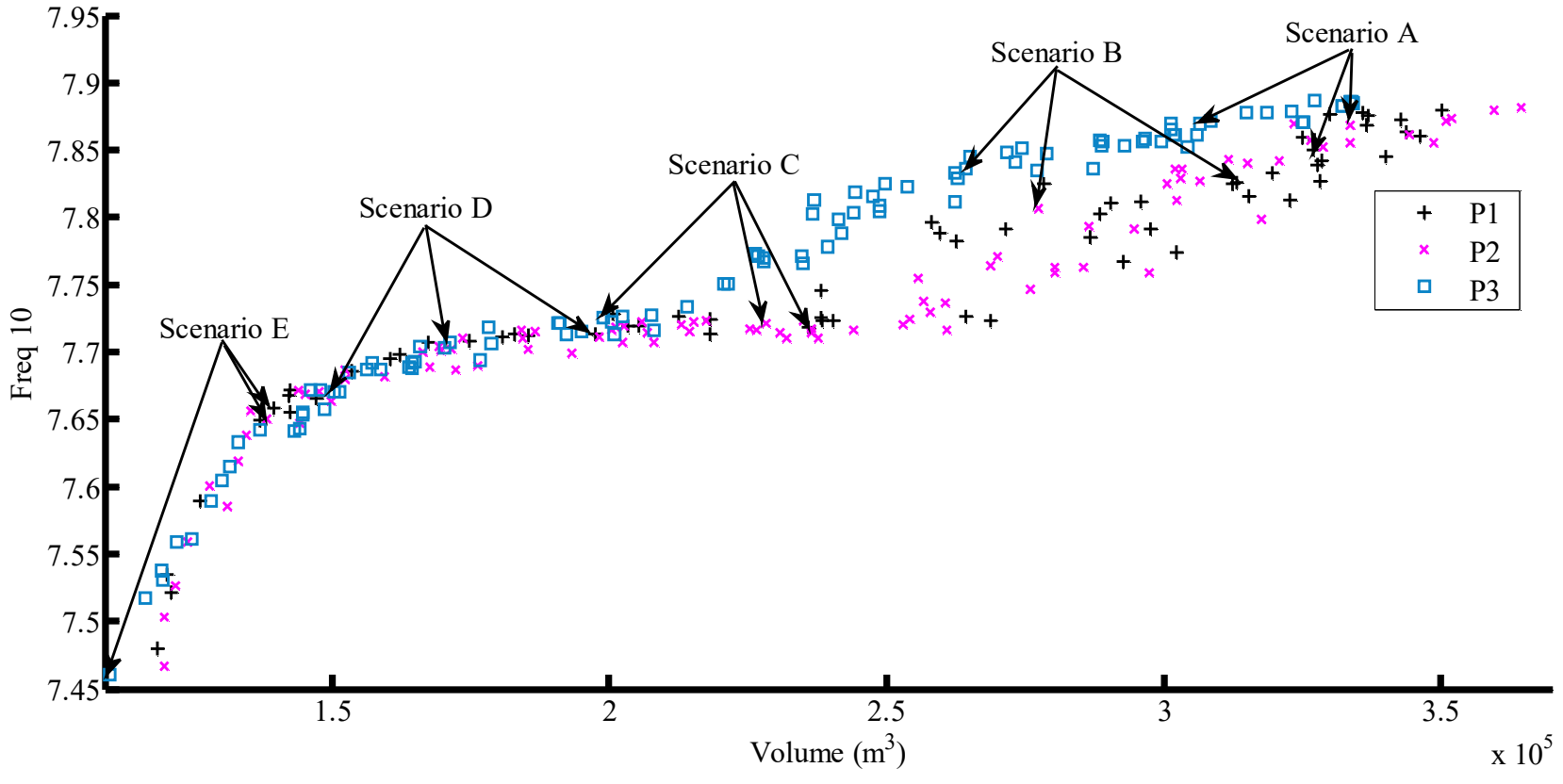

(j)

Fig. 11. Best solutions according to the five different scenarios for P1, P2 and P3 by MoCSS

Table 2: The lower and upper bounds of design variables

| | | |
|---|---|---|
| $0 \le \gamma \le 0.3$ | $0.5 \le \beta \le 1$ | |
| $3m \le t_{c1} \le 10m$ | $104m \le r_{u1} \le 135m$ | $104m \le r_{d1} \le 135m$ |
| $5m \le t_{c2} \le 14m$ | $91m \le r_{u2} \le 118m$ | $91m \le r_{d2} \le 118m$ |
| $7m \le t_{c3} \le 19m$ | $78m \le r_{u3} \le 101m$ | $78m \le r_{d3} \le 101m$ |
| $9m \le t_{c4} \le 23m$ | $65m \le r_{u4} \le 85m$ | $65m \le r_{d4} \le 85m$ |
| $11m \le t_{c5} \le 26m$ | $52m \le r_{u5} \le 68m$ | $52m \le r_{d5} \le 68m$ |
| $12m \le t_{c6} \le 31m$ | $39m \le r_{u6} \le 51m$ | $39m \le r_{d6} \le 51m$ |

Table 3: Material Property of Mass Concrete and Water

| Material | Property | Value | unit |
|---|---|---|---|
| **Mass Concrete** | Modulus of Elasticity | 27.579 | GP |
| | Poisson's Ratio | 0.2 | - |
| | Mass Density | 2483 | $Kg/m^3$ |
| **Water** | Bulk Modulus | 2.15 | GP |
| | Mass Density | 1000 | Kg/m3 |
| | Velocity of Pressure Waves | 1438.66 | m/s |
| | Wave reflection coefficient | 1 | - |

Table 4: Comparison the natural frequencies from the literature with the FEM

| mode | Tan and Chopra [24] | | Duron and Hall [25] | | Present work | |
|---|---|---|---|---|---|---|
| | empty | Full | Full | | empty | full |
| | FEA | | FEA | experimental | FEA | |
| **1** | 4.27 | 2.82 | 3.05 | 2.95 | 4.2897 | 2.9607 |
| **2** | - | - | - | 3.3 | - | 3.3046 |
| **3** | - | - | 4.21 | 3.95 | - | 4.2425 |
| **4** | - | - | 5.96 | 5.4 | - | 4.9522 |
| **5** | - | - | - | 6.21 | - | 6.1156 |

Table 5: Parameters of the MoCSS method

| | | |
|---|---|---|
| **MoCSS** | Number of agents | 100 |
| | Maximum number of iterations | 200 |
| | $k_a$ | 2 |
| | $k_v$ | 2 |
| | $\alpha$ | 1 |
| | $k_t$ | 0.75 |
| **MoPSO** | Swarm size | 100 |
| | Archive size | 100 |
| | Mutation rate | 0.5 |
| | Number of divisions for the adaptive grid | 30 |
| **NSGA-II** | Population size | 100 |
| | Crossover probability ( $p_c$ ) | 0.9 |
| | Mutation rate | $1/N_{\text{var}}$ |
| | Crossover distribution rate ( $\eta_c$ ) | 20 |
| | Mutation distribution rate ( $\eta_m$ ) | 20 |

Table 6: Extreme values for the problem P1

| Obtained extreme values | MoCSS | | | | | | | | | | | NSGA-II | | | | | | | | | | | MoPSO | | | | | | | | | | |
|---|---|---|---|---|---|---|---|---|---|---|---|---|---|---|---|---|---|---|---|---|---|---|---|---|---|---|---|---|---|---|---|---|---|
| | ***fit1*** | ***fit2*** | ***fit3*** | ***fit4*** | ***fit5*** | ***fit6*** | ***fit7*** | ***fit8*** | ***fit9*** | ***fit10*** | ***fit11*** | ***fit1*** | ***fit2*** | ***fit3*** | ***fit4*** | ***fit5*** | ***fit6*** | ***fit7*** | ***fit8*** | ***fit9*** | ***fit10*** | ***fit11*** | ***fit1*** | ***fit2*** | ***fit3*** | ***fit4*** | ***fit5*** | ***fit6*** | ***fit7*** | ***fit8*** | ***fit9*** | ***fit10*** | ***fit11*** |
| ***Vol (1e5)*** | **1.099** | 3.521 | 3.521 | 3.521 | 3.521 | 3.521 | 3.476 | 3.521 | 3.521 | 3.521 | 3.521 | 1.099 | 3.501 | 3.463 | 3.427 | 3.501 | 3.501 | 3.427 | 3.501 | 3.278 | 3.501 | 3.501 | 1.099 | 3.545 | 3.545 | 3.545 | 3.329 | 3.329 | 3.478 | 3.478 | 3.478 | 3.329 | 3.478 |
| ***Fr1*** | 2.579 | 2.929 | 2.929 | 2.929 | 2.929 | 2.929 | 2.927 | 2.929 | 2.929 | 2.929 | 2.929 | 2.514 | **2.931** | 2.931 | 2.924 | 2.930 | 2.930 | 2.924 | 2.930 | 2.930 | 2.930 | 2.930 | 2.514 | 2.926 | 2.926 | 2.926 | 2.926 | 2.926 | 2.926 | 2.926 | 2.926 | 2.922 | 2.926 |
| ***Fr2*** | 3.105 | 3.318 | **3.318** | 3.318 | 3.318 | 3.318 | 3.313 | 3.318 | 3.318 | 3.318 | 3.318 | 3.101 | 3.315 | 3.318 | 3.315 | 3.315 | 3.315 | 3.315 | 3.315 | 3.315 | 3.315 | 3.315 | 3.101 | 3.312 | 3.312 | 3.312 | 3.311 | 3.312 | 3.311 | 3.311 | 3.311 | 3.304 | 3.311 |
| ***Fr3*** | 3.699 | 4.366 | 4.366 | **4.366** | 4.366 | 4.366 | 4.352 | 4.366 | 4.366 | 4.366 | 4.366 | 3.708 | 4.351 | 4.352 | 4.352 | 4.351 | 4.351 | 4.352 | 4.351 | 4.351 | 4.351 | 4.351 | 3.708 | 4.342 | 4.342 | 4.342 | 4.341 | 4.342 | 4.341 | 4.341 | 4.341 | 4.328 | 4.341 |
| ***Fr4*** | 4.441 | 5.312 | 5.312 | 5.312 | **5.315** | 5.315 | 5.198 | 5.312 | 5.312 | 5.312 | 5.312 | 4.421 | 5.215 | 5.151 | 5.160 | 5.215 | 5.215 | 5.160 | 5.215 | 5.215 | 5.215 | 5.215 | 4.421 | 5.136 | 5.136 | 5.136 | 5.183 | 5.136 | 5.142 | 5.142 | 5.142 | 5.183 | 5.142 |
| ***Fr5*** | 5.044 | 6.211 | 6.211 | 6.211 | 6.211 | **6.211** | 6.176 | 6.211 | 6.211 | 6.211 | 6.211 | 4.845 | 6.181 | 6.152 | 6.160 | 6.181 | 6.181 | 6.160 | 6.181 | 6.181 | 6.181 | 6.181 | 4.845 | 6.157 | 6.157 | 6.157 | 6.173 | 6.173 | 6.155 | 6.155 | 6.155 | 6.173 | 6.155 |
| ***Fr6*** | 5.711 | 6.682 | 6.682 | 6.682 | 6.682 | 6.682 | 6.772 | 6.682 | 6.682 | 6.682 | 6.682 | 5.572 | 6.712 | 6.738 | 6.813 | 6.712 | 6.712 | **6.813** | 6.712 | 6.712 | 6.712 | 6.712 | 5.572 | 6.793 | 6.793 | 6.793 | 6.677 | 6.793 | 6.800 | 6.800 | 6.800 | 6.677 | 6.800 |
| ***Fr7*** | 6.753 | 7.367 | 7.367 | 7.367 | 7.367 | 7.367 | 7.354 | **7.367** | 7.367 | 7.367 | 7.367 | 6.568 | 7.359 | 7.352 | 7.350 | 7.359 | 7.359 | 7.350 | 7.359 | 7.359 | 7.359 | 7.359 | 6.568 | 7.351 | 7.351 | 7.351 | 7.360 | 7.360 | 7.355 | 7.355 | 7.355 | 7.360 | 7.355 |
| ***Fr8*** | 7.030 | 7.725 | 7.725 | 7.725 | 7.725 | 7.725 | 7.707 | 7.725 | 7.725 | 7.725 | 7.725 | 6.943 | 7.719 | 7.711 | 7.685 | 7.719 | 7.719 | 7.685 | 7.719 | 7.721 | 7.719 | 7.719 | 6.943 | 7.702 | 7.702 | 7.702 | 7.715 | 7.715 | 7.793 | 7.793 | **7.793** | 7.715 | 7.693 |
| ***Fr9*** | 7.370 | 7.821 | 7.821 | 7.821 | 7.821 | 7.821 | 7.795 | 7.821 | 7.821 | 7.821 | 7.821 | 7.376 | 7.806 | 7.748 | 7.758 | 7.806 | 7.806 | 7.758 | 7.806 | 7.740 | 7.806 | 7.806 | 7.376 | 7.765 | 7.765 | 7.765 | 7.795 | 7.795 | 7.760 | 7.760 | 7.760 | **7.795** | 7.760 |
| ***Fr10*** | 7.475 | 7.900 | 7.900 | 7.900 | 7.900 | 7.900 | 7.882 | 7.900 | 7.900 | 7.900 | **7.900** | 7.404 | 7.880 | 7.860 | 7.873 | 7.880 | 7.880 | 7.873 | 7.880 | 7.839 | 7.880 | 7.880 | 7.404 | 7.876 | 7.876 | 7.876 | 7.879 | 7.879 | 7.884 | 7.884 | 7.884 | 7.879 | 7.884 |

Table 7: Comparison of the results for problem P2

| **Optimization methods** | **MoCSS** | | **NSGA-II** | | **MoPSO** | |
|---|---|---|---|---|---|---|
| **Obtained extreme values** | *fit1* | *fit2* | *fit1* | *fit2* | *fit1* | *fit2* |
| ***Vol (1e5)*** | **1.197** | 3.643 | 1.369 | 3.560 | 1.250 | 3.522 |
| ***Fr1*** | 2.586 | 2.927 | 2.640 | 2.927 | 2.609 | **2.930** |
| ***Fr2*** | 3.105 | **3.318** | 3.119 | 3.314 | 3.099 | 3.318 |
| ***Fr3*** | 3.705 | **4.352** | 3.761 | 4.340 | 3.723 | 4.351 |
| ***Fr4*** | 4.436 | **5.147** | 4.508 | 5.118 | 4.463 | 5.118 |
| ***Fr5*** | 4.974 | **6.160** | 5.101 | 6.157 | 5.050 | 6.147 |
| ***Fr6*** | 5.614 | **6.817** | 5.695 | 6.785 | 5.624 | 6.755 |
| ***Fr7*** | 6.583 | **7.353** | 6.691 | 7.348 | 6.617 | 7.346 |
| ***Fr8*** | 7.000 | 7.699 | 7.125 | 7.702 | 7.706 | **7.716** |
| ***Fr9*** | 7.378 | 7.765 | 7.396 | **7.768** | 7.386 | 7.744 |
| ***Fr10*** | 7.467 | **7.881** | 7.608 | 7.866 | 7.531 | 7.846 |

Table 8: Comparison of the results for problem P3

| Optimization methods | MoCSS | | NSGA-II | | MoPSO | |
|---|---|---|---|---|---|---|
| Obtained extreme values | *fit1* | *fit2* | *fit1* | *fit2* | *fit1* | *fit2* |
| ***Vol (1e5)*** | **1.099** | 3.340 | **1.099** | 3.288 | **1.099** | 3.280 |
| ***Fr1*** | 2.572 | 2.922 | 2.514 | **2.924** | 2.572 | 2.918 |
| ***Fr2*** | 3.104 | **3.305** | 3.101 | 3.304 | 3.104 | 3.302 |
| ***Fr3*** | 3.691 | **4.332** | 3.708 | 4.326 | 3.691 | 4.326 |
| ***Fr4*** | 4.433 | **5.196** | 4.421 | 5.167 | 4.433 | 5.196 |
| ***Fr5*** | 5.015 | 6.171 | 4.845 | 6.160 | 5.015 | **6.174** |
| ***Fr6*** | 5.689 | 6.664 | 5.572 | **6.675** | 5.689 | 6.638 |
| ***Fr7*** | 6.699 | 7.362 | 6.568 | 7.360 | 6.699 | **7.364** |
| ***Fr8*** | 7.008 | **7.716** | 6.943 | 7.712 | 7.008 | 7.713 |
| ***Fr9*** | 7.370 | 7.784 | 7.376 | 7.782 | 7.370 | **7.786** |
| ***Fr10*** | 7.460 | 7.885 | 7.404 | 7.878 | 7.460 | **7.886** |

Table 9: Different possible scenarios for problems P1 corresponding solutions

| Scenario | Importance of criteria | Possible priority weights | Selected solution by MTDM | | | | | | | | | | |
|---|---|---|---|---|---|---|---|---|---|---|---|---|---|
| | | | **MoCSS** | | | | | | | | | | |
| | | | *Vol* | *Fr 1* | *Fr 2* | *Fr 3* | *Fr 4* | *Fr 5* | *Fr 6* | *Fr 6* | *Fr 7* | *Fr 8* | *Fr 10* |
| A | $C1>>C2$ | [0.9,0.1] | 327156.241 | 2.927 | 3.307 | 4.329 | 5.108 | 6.142 | 6.708 | 7.354 | 7.713 | 7.743 | 7.857 |
| B | $C1>C2$ | [0.7,0.3] | 313206.341 | 2.925 | 3.301 | 4.307 | 5.041 | 6.128 | 6.633 | 7.347 | 7.704 | 7.728 | 7.826 |
| C | $C1\approx C2$ | [0.5,0.5] | 236039.589 | 2.890 | 3.237 | 4.118 | 4.790 | 6.018 | 6.266 | 7.333 | 7.344 | 7.707 | 7.718 |
| D | $C1<C2$ | [0.3,0.7] | 197471.304 | 2.842 | 3.197 | 4.014 | 4.706 | 5.942 | 5.949 | 7.191 | 7.320 | 7.606 | 7.713 |
| E | $C1<<C2$ | [0.1,0.9] | 139509.614 | 2.679 | 3.134 | 3.841 | 4.559 | 5.340 | 5.747 | 6.836 | 7.226 | 7.428 | 7.659 |
| | | | **NSGA-II** | | | | | | | | | | |
| | | | *Vol* | *Fr 1* | *Fr 2* | *Fr 3* | *Fr 4* | *Fr 5* | *Fr 6* | *Fr 6* | *Fr 7* | *Fr 8* | *Fr 10* |
| A | $C1>>C2$ | [0.9,0.1] | 345978.388 | 2.929 | 3.311 | 4.344 | 5.188 | 6.174 | 6.703 | 7.360 | 7.718 | 7.797 | 7.884 |
| B | $C1>C2$ | [0.7,0.3] | 333857.067 | 2.927 | 3.307 | 4.328 | 5.115 | 6.149 | 6.719 | 7.357 | 7.706 | 7.754 | 7.874 |
| C | $C1\approx C2$ | [0.5,0.5] | 250685.935 | 2.904 | 3.257 | 4.200 | 4.896 | 6.070 | 6.438 | 7.341 | 7.534 | 7.716 | 7.781 |
| D | $C1<C2$ | [0.3,0.7] | 156724.220 | 2.794 | 3.122 | 3.850 | 4.604 | 5.662 | 5.865 | 7.039 | 7.274 | 7.478 | 7.682 |
| E | $C1<<C2$ | [0.1,0.9] | 122419.420 | 2.659 | 3.058 | 3.706 | 4.499 | 5.191 | 5.726 | 6.841 | 7.124 | 7.384 | 7.587 |
| | | | **MoPSO** | | | | | | | | | | |
| | | | *Vol* | *Fr 1* | *Fr 2* | *Fr 3* | *Fr 4* | *Fr 5* | *Fr 6* | *Fr 6* | *Fr 7* | *Fr 8* | *Fr 10* |
| A | $C1>>C2$ | [0.9,0.1] | 330761.861 | 2.925 | 3.302 | 4.316 | 5.103 | 6.148 | 6.718 | 7.353 | 7.684 | 7.752 | 7.868 |
| B | $C1>C2$ | [0.7,0.3] | 296607.241 | 2.920 | 3.290 | 4.278 | 5.024 | 6.125 | 6.586 | 7.348 | 7.642 | 7.729 | 7.823 |
| C | $C1\approx C2$ | [0.5,0.5] | 250204.361 | 2.901 | 3.245 | 4.145 | 4.823 | 6.045 | 6.423 | 7.336 | 7.416 | 7.707 | 7.748 |
| D | $C1<C2$ | [0.3,0.7] | 173154.045 | 2.812 | 3.157 | 3.901 | 4.624 | 5.753 | 5.889 | 7.043 | 7.292 | 7.513 | 7.698 |
| E | $C1<<C2$ | [0.1,0.9] | 144606.744 | 2.676 | 3.128 | 3.782 | 4.533 | 5.222 | 5.747 | 6.748 | 7.206 | 7.421 | 7.648 |

Table 10: Different possible scenarios for problem P2 corresponding solutions

| Scenario | Importance of criteria | Possible priority weights | Selected solution by MTDM | | | | | | | | | | |
|---|---|---|---|---|---|---|---|---|---|---|---|---|---|
| | | | **MoCSS** | | | | | | | | | | |
| | | | ***Vol*** | ***Fr 1*** | ***Fr 2*** | ***Fr 3*** | ***Fr 4*** | ***Fr 5*** | ***Fr 6*** | ***Fr 6*** | ***Fr 7*** | ***Fr 8*** | ***Fr 10*** |
| A | $C1 \gg C2$ | [0.9,0.1] | **326933.129** | **2.923** | **3.301** | **4.304** | **5.051** | **6.141** | **6.712** | **7.354** | **7.692** | **7.746** | **7.860** |
| B | $C1 > C2$ | [0.7,0.3] | **269391.844** | **2.907** | **3.260** | **4.179** | **4.844** | **6.066** | **6.482** | **7.336** | **7.511** | **7.717** | **7.773** |
| C | $C1 \approx C2$ | [0.5,0.5] | **226550.518** | **2.880** | **3.225** | **4.076** | **4.769** | **5.997** | **6.018** | **7.307** | **7.322** | **7.618** | **7.723** |
| D | $C1 < C2$ | [0.3,0.7] | **192550.177** | **2.828** | **3.172** | **3.930** | **4.665** | **5.816** | **5.928** | **7.066** | **7.302** | **7.547** | **7.713** |
| E | $C1 \ll C2$ | [0.1,0.9] | **145201.120** | **2.710** | **3.129** | **3.815** | **4.549** | **5.287** | **5.742** | **6.757** | **7.196** | **7.421** | **7.629** |
| | | | **NSGA-II** | | | | | | | | | | |
| | | | ***Vol*** | ***Fr 1*** | ***Fr 2*** | ***Fr 3*** | ***Fr 4*** | ***Fr 5*** | ***Fr 6*** | ***Fr 6*** | ***Fr 7*** | ***Fr 8*** | ***Fr 10*** |
| A | $C1 \gg C2$ | [0.9,0.1] | **327223.085** | **2.921** | **3.301** | **4.301** | **5.035** | **6.139** | **6.735** | **7.351** | **7.675** | **7.739** | **7.858** |
| B | $C1 > C2$ | [0.7,0.3] | **283797.653** | **2.909** | **3.263** | **4.176** | **4.852** | **6.093** | **6.524** | **7.330** | **7.587** | **7.714** | **7.761** |
| C | $C1 \approx C2$ | [0.5,0.5] | **232374.178** | **2.876** | **3.225** | **4.064** | **4.754** | **5.988** | **6.026** | **7.323** | **7.345** | **7.621** | **7.728** |
| D | $C1 < C2$ | [0.3,0.7] | **200606.877** | **2.841** | **3.171** | **3.922** | **4.660** | **5.747** | **5.951** | **7.110** | **7.299** | **7.534** | **7.715** |
| E | $C1 \ll C2$ | [0.1,0.9] | **148505.290** | **2.709** | **3.131** | **3.804** | **4.541** | **5.284** | **5.745** | **6.756** | **7.189** | **7.416** | **7.643** |
| | | | **MoPSO** | | | | | | | | | | |
| | | | ***Vol*** | ***Fr 1*** | ***Fr 2*** | ***Fr 3*** | ***Fr 4*** | ***Fr 5*** | ***Fr 6*** | ***Fr 6*** | ***Fr 7*** | ***Fr 8*** | ***Fr 10*** |
| A | $C1 \gg C2$ | [0.9,0.1] | **333535.467** | **2.926** | **3.304** | **4.323** | **5.115** | **6.155** | **6.701** | **7.355** | **7.711** | **7.772** | **7.868** |
| B | $C1 > C2$ | [0.7,0.3] | **277474.385** | **2.905** | **3.266** | **4.209** | **4.907** | **6.073** | **6.438** | **7.349** | **7.520** | **7.725** | **7.807** |
| C | $C1 \approx C2$ | [0.5,0.5] | **228311.990** | **2.879** | **3.224** | **4.072** | **4.759** | **6.014** | **6.141** | **7.328** | **7.335** | **7.671** | **7.722** |
| D | $C1 < C2$ | [0.3,0.7] | **171550.040** | **2.793** | **3.162** | **3.906** | **4.616** | **5.689** | **5.858** | **6.953** | **7.284** | **7.486** | **7.702** |
| E | $C1 \ll C2$ | [0.1,0.9] | **138294.760** | **2.667** | **3.128** | **3.820** | **4.552** | **5.249** | **5.733** | **6.771** | **7.202** | **7.418** | **7.650** |

Table 11: Different possible scenarios for problem P3 corresponding solutions

| Scenario | Importance of criteria | Possible priority weights | Selected solution by MTDM | | | | | | | | | | |
|---|---|---|---|---|---|---|---|---|---|---|---|---|---|
| | | | **MoCSS** | | | | | | | | | | |
| | | | ***Vol*** | ***Fr 1*** | ***Fr 2*** | ***Fr 3*** | ***Fr 4*** | ***Fr 5*** | ***Fr 6*** | ***Fr 6*** | ***Fr 7*** | ***Fr 8*** | ***Fr 10*** |
| A | $C1 \gg C2$ | [0.9,0.1] | **306505.148** | **2.91** | **3.29** | **4.30** | **5.12** | **6.14** | **6.61** | **7.35** | **7.69** | **7.76** | **7.86** |
| B | $C1 > C2$ | [0.7,0.3] | **264968.249** | **2.90** | **3.26** | **4.24** | **5.00** | **6.09** | **6.52** | **7.35** | **7.61** | **7.73** | **7.84** |
| C | $C1 \approx C2$ | [0.5,0.5] | **198943.337** | **2.84** | **3.21** | **4.10** | **4.80** | **5.97** | **6.12** | **7.31** | **7.34** | **7.72** | **7.72** |
| D | $C1 < C2$ | [0.3,0.7] | **150206.901** | **2.72** | **3.15** | **3.93** | **4.62** | **5.65** | **5.80** | **6.98** | **7.31** | **7.55** | **7.67** |
| E | $C1 \ll C2$ | [0.1,0.9] | **109942.578** | **2.571** | **3.10** | **3.69** | **4.43** | **5.01** | **5.69** | **6.70** | **7.00** | **7.36** | **7.46** |
| | | | **NSGA-II** | | | | | | | | | | |
| | | | ***Vol*** | ***Fr 1*** | ***Fr 2*** | ***Fr 3*** | ***Fr 4*** | ***Fr 5*** | ***Fr 6*** | ***Fr 6*** | ***Fr 7*** | ***Fr 8*** | ***Fr 10*** |
| A | $C1 \gg C2$ | [0.9,0.1] | **306734.745** | **2.918** | **3.292** | **4.302** | **5.122** | **6.145** | **6.603** | **7.360** | **7.690** | **7.759** | **7.870** |
| B | $C1 > C2$ | [0.7,0.3] | **269714.224** | **2.904** | **3.273** | **4.248** | **5.005** | **6.101** | **6.481** | **7.353** | **7.620** | **7.740** | **7.830** |
| C | $C1 \approx C2$ | [0.5,0.5] | **226966.502** | **2.881** | **3.237** | **4.156** | **4.858** | **6.025** | **6.324** | **7.344** | **7.418** | **7.724** | **7.774** |
| D | $C1 < C2$ | [0.3,0.7] | **171357.449** | **2.779** | **3.180** | **4.025** | **4.729** | **5.885** | **5.908** | **7.155** | **7.339** | **7.694** | **7.708** |
| E | $C1 \ll C2$ | [0.1,0.9] | **132956.454** | **2.643** | **3.132** | **3.852** | **4.548** | **5.326** | **5.716** | **6.838** | **7.246** | **7.435** | **7.633** |
| | | | **MoPSO** | | | | | | | | | | |
| | | | ***Vol*** | ***Fr 1*** | ***Fr 2*** | ***Fr 3*** | ***Fr 4*** | ***Fr 5*** | ***Fr 6*** | ***Fr 6*** | ***Fr 7*** | ***Fr 8*** | ***Fr 10*** |
| A | $C1 \gg C2$ | [0.9,0.1] | **314937.556** | **2.920** | **3.296** | **4.311** | **5.156** | **6.165** | **6.616** | **7.362** | **7.701** | **7.778** | **7.877** |
| B | $C1 > C2$ | [0.7,0.3] | **271749.572** | **2.902** | **3.275** | **4.259** | **5.040** | **6.109** | **6.503** | **7.356** | **7.636** | **7.748** | **7.848** |
| C | $C1 \approx C2$ | [0.5,0.5] | **227666.621** | **2.886** | **3.240** | **4.160** | **4.864** | **6.031** | **6.354** | **7.342** | **7.425** | **7.722** | **7.767** |
| D | $C1 < C2$ | [0.3,0.7] | **181474.149** | **2.824** | **3.187** | **4.033** | **4.747** | **5.929** | **5.971** | **7.191** | **7.334** | **7.675** | **7.713** |
| E | $C1 \ll C2$ | [0.1,0.9] | **150504.796** | **2.754** | **3.143** | **3.908** | **4.610** | **5.678** | **5.813** | **6.977** | **7.306** | **7.533** | **7.671** |